\documentclass[11pt]{article}

\usepackage[a4paper,margin=32mm]{geometry}
\usepackage[T1]{fontenc}
\usepackage{lmodern}
\usepackage{microtype}
\usepackage{mathtools,amssymb,mathrsfs}
\usepackage{amsthm}
\usepackage[hidelinks]{hyperref}

\newtheorem{proposition}{Proposition}[section]
\newtheorem{theorem}[proposition]{Theorem}
\newtheorem{remark}[proposition]{Remark}
\newtheorem{corollary}[proposition]{Corollary}
\newtheorem{lemma}[proposition]{Lemma}
\newtheorem{conjecture}[proposition]{Conjecture}
\numberwithin{equation}{section}

\title{\bfseries Determinant Identities and Branch Selection\\
in the Stationary Multi-Horizon Balance Problem}
\author{Daheng Min\thanks{Funded by the Deutsche Forschungsgemeinschaft (DFG,
German Research Foundation) under Germany's Excellence Strategy
EXC 2044-390685587, Mathematics M\"unster: Dynamics--Geometry--Structure.}}
\date{}

\begin{document}

\maketitle

\begin{abstract}
We study the Hennig--Neugebauer family of stationary and axisymmetric vacuum
multi-horizon candidates. The regularity conditions on the symmetry axis lead to
quadratic determinant relations which leave an apparent ambiguity between two sign
branches. Under explicit non-degeneracy assumptions, we prove that one of these branches
is impossible. This branch selection yields a system of explicit determinant
equations characterizing equilibrium within the Hennig--Neugebauer family.
\end{abstract}

\section{Introduction}

The question whether several aligned black holes can form a stationary equilibrium
configuration is a long-standing problem in general relativity. In an asymptotically
flat, stationary and axisymmetric vacuum spacetime, such an equilibrium would require
all ordinary portions of the symmetry axis to be regular. In particular, the
elementary-flatness condition must hold on every axis interval separating the horizon
components; otherwise the metric contains conical singularities that may be
interpreted as strings or struts supporting the configuration. Soliton methods reduce
the corresponding boundary-value problem to the study of finite-parameter families of
exact candidate solutions. The remaining problem is therefore to determine whether
their axis regularity conditions can be satisfied and whether the resulting candidates
possess the additional global and geometric properties required of a multi-black-hole
spacetime.

The cases of one and two horizon components have been studied in considerable detail.
For a single connected horizon, the soliton construction reproduces the Kerr family
and provides a constructive version of the corresponding uniqueness result. For two
components, the associated boundary-value problem leads to the double-Kerr--NUT
family. Tomimatsu and Kihara carried out an explicit analysis of the metric functions
on the symmetry axis and derived the algebraic conditions required for the absence of
axis singularities \cite{TomimatsuKihara1982}. An important step in their argument is
the exclusion of one of the two sign choices left undetermined by the
elementary-flatness condition, using a sum of non-negative terms. Hennig and Neugebauer
subsequently showed that every sufficiently regular stationary two-black-hole
configuration in vacuum must belong to the double-Kerr--NUT family and, by combining
its explicit parametrization with the angular-momentum--area inequality, proved that
two separated sub-extremal black holes cannot be in stationary equilibrium
\cite{NeugebauerHennig2009}.

Until recently, the balance problem for three or more components remained open.
General multi-soliton candidates were obtained and studied by Neugebauer, Kramer and
Yamazaki, who also derived algebraic formulas for the associated metric functions
\cite{NeugebauerBT,Kramer1986,Yamazaki1983}. Varzugin subsequently obtained
corresponding expressions for the interaction forces between neighbouring components
\cite{Varzugin1997,Varzugin1998}. More recently, Hennig derived the general rational
form of the axis potentials compatible with the boundary conditions for an arbitrary
number of aligned horizon components, thereby reducing the problem to the analysis of
finite-parameter families of candidates \cite{Hennig2020,Hennig2025,Hennig2026}.

Complementary to the soliton approach, the stationary axisymmetric vacuum equations
also admit a harmonic-map formulation. Recently, Han, Khuri, Weinstein and Xiong
proved that, under their prescribed regularity and asymptotic hypotheses, every bounded
axis rod in a multi-horizon configuration has strictly negative logarithmic angle
defect; the angular momentum is additionally assumed to be nonzero for every degenerate
horizon component. Consequently, no regular asymptotically flat equilibrium
configuration with more than one horizon component exists in this class. Together with
the standard rigidity and connected-horizon uniqueness theorems, this yields Kerr
uniqueness in the corresponding analytic rotating setting
\cite{HanKhuriWeinsteinXiong2026}.

From the viewpoint of the present determinant construction, this result shifts the
remaining question: one should determine which hypothesis of Theorem~1.1 of
\cite{HanKhuriWeinsteinXiong2026}, or equivalently which part of its harmonic-map
regularity and asymptotic framework, fails for the equilibrium candidates produced by
the Hennig--Neugebauer determinants. Although the determinant identities established
below do not by themselves settle this regularity question, they isolate the admissible
equilibrium branch and provide explicit boundary formulas needed for a systematic
comparison with those hypotheses.

Yamazaki's determinant formulas are essentially equivalent to the axis representation
used in the present work. In formulating the elementary-flatness conditions, however,
he selected one of the two sign branches allowed by the corresponding quadratic
relation, without excluding the alternative for general~$N$. For $N=2$, this choice is
justified by the non-negative identity of Tomimatsu and Kihara. Our main result supplies
the missing branch-selection argument for arbitrary~$N$.

The Hennig--Neugebauer candidate considered here is determined by ordered real
parameters
\begin{equation*}
 K_1>\cdots>K_{2N}
\end{equation*}
and unit complex parameters
\begin{equation*}
 \alpha_1,\ldots,\alpha_{2N},
 \qquad |\alpha_i|=1.
\end{equation*}
For every axis interval $I$, the determinant construction associates two explicit
algebraic quantities
\begin{equation*}
 H_N^I=H_N^I(K_1,\ldots,K_{2N};\alpha_1,\ldots,\alpha_{2N}),
 \qquad
 D_N^I=D_N^I(K_1,\ldots,K_{2N};\alpha_1,\ldots,\alpha_{2N}).
\end{equation*}
The equilibrium problem is therefore reduced to determining which parameter values
satisfy the determinant equations expressing regularity of all ordinary axis intervals.

\begin{theorem}[Equilibrium equations for the Hennig--Neugebauer candidates]
\label{thm:introductionMainResult}
Under the normalization and non-degeneracy assumptions of
Theorem~\ref{thm:branchSelection}, the Hennig--Neugebauer candidate determined by
$(K_i,\alpha_i)_{i=1}^{2N}$ satisfies the regular-axis equilibrium conditions if and
only if
\begin{equation*}
 D_N^{\mathcal A^j}=(-1)^jD_N^+,
 \qquad j=1,\ldots,N,
\end{equation*}
and
\begin{equation*}
 H_N^{\mathcal A^j}=(-1)^jH_N^+,
 \qquad j=1,\ldots,N-1.
\end{equation*}
\end{theorem}

Thus the possible equilibrium parameters are selected by $2N-1$ explicit determinant
equations in $K_i$ and $\alpha_i$. The $D_N^I$-equations encode equality of the
boundary values of the rotational potential $a$, whereas the $H_N^I$-equations encode
elementary flatness. The detailed statement and proof are given in
Theorem~\ref{thm:branchSelection}.

The paper is organized as follows. Section~2 reviews the Weyl--Papapetrou reduction,
the Ernst equation and its linear problem, and establishes the Darboux criterion for
polynomial transformation matrices. Section~3 introduces the Hennig--Neugebauer
determinant family, explains its origin from the integration of the linear problem along
the boundary, and verifies that it satisfies the Darboux criterion. In Section~4, we
derive explicit determinant formulas for the boundary values of the Ernst potential and
the metric functions $e^{2U}$, $e^{2k}$, and $a$; the required determinant identities
are obtained principally from Pl\"ucker relations. Section~5 translates regularity of
the ordinary axis intervals into equilibrium equations and proves the branch-selection
theorem stated above. Section~6 records the resulting formulas for the candidate horizon
quantities. Section~7 examines the cases $N=1,2,3$, while Section~8 collects several
algebraic conjectures and discusses their possible relation to the global regularity
problem.

All conclusions concerning equilibrium in this paper are algebraic conditions on the
determinant candidates and do not, by themselves, establish global regularity of the
reconstructed spacetime.

\section{The Ernst equation and its linear problem}

\subsection{Weyl--Papapetrou coordinates}

Let
\[
 \mathscr O\subset(0,\infty)\times\mathbb R
\]
be a connected open set.  We use coordinates with ranges
\begin{equation}
 (\rho,\zeta)\in\mathscr O,
 \qquad
 \varphi\in\mathbb R/(2\pi\mathbb Z),
 \qquad
 t\in\mathbb R.
 \label{eq:WPranges}
\end{equation}
Thus the coordinate domain under consideration is
\(\mathscr O\times S^1\times\mathbb R\).  The coordinate \(\varphi\)
parametrizes the orbits of the axial symmetry and is understood modulo
\(2\pi\), whereas \(t\) parametrizes the stationary symmetry.  Since
\(\rho>0\) on \(\mathscr O\), the symmetry axis is not contained in this
coordinate domain; it will later be approached through the boundary limit
\(\rho\to0\).

Let
\[
 F,E,a:\mathscr O\longrightarrow\mathbb R
\]
be sufficiently regular real-valued functions, and assume
\begin{equation}
 F\ne0\qquad\text{on }\mathscr O.
 \label{eq:Fnonzero}
\end{equation}
Following the conventional notation, we write
\begin{equation}
 F=e^{2U},\qquad E=e^{2k}.
 \label{eq:FEnotation}
\end{equation}
These are names for real-valued functions and do not assert positivity.  In
particular, \(e^{-2U}\) means \(1/F\), which explains the necessity of
\eqref{eq:Fnonzero}.

In Weyl--Papapetrou coordinates, consider the stationary and axisymmetric
line element
\begin{equation}
 ds^2=e^{-2U}\bigl[e^{2k}(d\rho^2+d\zeta^2)+\rho^2d\varphi^2\bigr]
      -e^{2U}(dt+a\,d\varphi)^2 .
 \label{eq:WPmetric}
\end{equation}
Equivalently,
\begin{equation}
 ds^2=\frac EF(d\rho^2+d\zeta^2)
      +\frac{\rho^2}{F}d\varphi^2-F(dt+a\,d\varphi)^2.
 \label{eq:WPmetricFE}
\end{equation}
At this stage, \eqref{eq:WPmetric} is regarded as a symmetric quadratic
form.  It defines a Lorentzian metric only where it is non-degenerate and
has the required signature.

The restriction to the \((t,\varphi)\)-directions has matrix
\begin{equation}
 G_{t\varphi}=
 \begin{pmatrix}
 -F&-Fa\\[1mm]
 -Fa&\displaystyle\frac{\rho^2}{F}-Fa^2
 \end{pmatrix}.
 \label{eq:tphiblock}
\end{equation}
Its determinant is
\begin{equation}
 \det G_{t\varphi}=-\rho^2.
 \label{eq:tphideterminant}
\end{equation}
Since \(\rho>0\), this block is non-degenerate and has one positive and one
negative eigenvalue, independently of the signs of \(F\) and \(E\).  The
restriction to the \((\rho,\zeta)\)-directions is
\begin{equation}
 \frac EF(d\rho^2+d\zeta^2).
 \label{eq:rhozetaquadratic}
\end{equation}
If
\begin{equation}
 \frac EF>0,
 \label{eq:positiveRatio}
\end{equation}
the complete quadratic form has signature \((-+++)\).  If instead
\begin{equation}
 \frac EF<0,
 \label{eq:negativeRatio}
\end{equation}
it has signature \((+---)\).  Consequently, with our convention,
\eqref{eq:WPmetric} defines a Lorentzian metric of signature \((-+++)\)
precisely where
\begin{equation}
 E\ne0,
 \qquad
 \frac EF>0.
 \label{eq:Lorentzcondition}
\end{equation}
Since \(F\ne0\) has already been assumed, this is equivalent to requiring
that \(E\) be nonzero and have the same sign as \(F\).

If \(E=0\), the quadratic form \eqref{eq:WPmetric} is defined but
degenerate.  If \(F=0\), the Weyl--Papapetrou expression itself is not
defined because it contains \(F^{-1}\).  In either case, failure of this
particular expression does not by itself prove that the underlying spacetime
geometry is singular; a different coordinate representation may still
extend smoothly across the corresponding locus.

\subsection{The Ernst potential}

Introduce a real twist potential \(b\), locally and up to an additive
constant, by
\begin{equation}
 a_{,\rho}=\rho e^{-4U}b_{,\zeta},
 \qquad
 a_{,\zeta}=-\rho e^{-4U}b_{,\rho}.
 \label{eq:twist}
\end{equation}
The complex Ernst potential is
\begin{equation}
 f=e^{2U}+ib.
 \label{eq:Ernstdef}
\end{equation}
The vacuum field equations for \eqref{eq:WPmetric} imply
\begin{equation}
 (\Re f)\left(f_{,\rho\rho}+f_{,\zeta\zeta}
                    +\frac1\rho f_{,\rho}\right)
 =f_{,\rho}^2+f_{,\zeta}^2,
 \label{eq:Ernst}
\end{equation}
wherever \(\Re f\ne0\).

Conversely, let \(f\) be a sufficiently regular solution of
\eqref{eq:Ernst} on a simply connected region on which \(\Re f\ne0\).  We
first define
\begin{equation}
 e^{2U}:=\Re f,
 \qquad b:=\Im f.
 \label{eq:recoverUb}
\end{equation}
The function \(a\) is then obtained, up to an additive constant, from
\begin{equation}
 a_{,\rho}=\frac{\rho}{(\Re f)^2}(\Im f)_{,\zeta},
 \qquad
 a_{,\zeta}=-\frac{\rho}{(\Re f)^2}(\Im f)_{,\rho}.
 \label{eq:recovera}
\end{equation}
The Ernst equation is the integrability condition for this system.

The remaining first-order equations are
\begin{align}
 k_{,\rho}
 &=\frac{\rho}{4(\Re f)^2}
 \left(f_{,\rho}\overline{f_{,\rho}}
       -f_{,\zeta}\overline{f_{,\zeta}}\right),
 \label{eq:krhof}\\
 k_{,\zeta}
 &=\frac{\rho}{4(\Re f)^2}
 \left(f_{,\rho}\overline{f_{,\zeta}}
       +f_{,\zeta}\overline{f_{,\rho}}\right).
 \label{eq:kzetaf}
\end{align}
Equivalently, in terms of \(U\) and \(b\),
\begin{align}
 k_{,\rho}
 &=\rho\left[U_{,\rho}^2-U_{,\zeta}^2
 +\frac14e^{-4U}(b_{,\rho}^2-b_{,\zeta}^2)\right],
 \label{eq:krhoUb}\\
 k_{,\zeta}
 &=2\rho\left[U_{,\rho}U_{,\zeta}
 +\frac14e^{-4U}b_{,\rho}b_{,\zeta}\right].
 \label{eq:kzetaUb}
\end{align}
The Ernst equation also implies the integrability condition for
\eqref{eq:krhof}--\eqref{eq:kzetaf}.

There is an important point concerning the notation \(k\).  The line
element contains the real-valued function \(e^{2k}\), whereas the quantity
differentiated in \eqref{eq:krhof}--\eqref{eq:kzetaf} is, on a connected
region where \(e^{2k}\ne0\),
\begin{equation}
 k:=\frac12\log|e^{2k}|.
 \label{eq:kabsolute}
\end{equation}
Consequently these equations determine \(|e^{2k}|\), up to a positive
multiplicative constant, but do not determine its sign.  On each connected
region where \(e^{2k}\ne0\), one may write
\begin{equation}
 e^{2k}=\varepsilon_k e^{2k_{\mathrm{abs}}},
 \qquad \varepsilon_k\in\{+1,-1\},
 \qquad k_{\mathrm{abs}}=\frac12\log|e^{2k}|.
 \label{eq:ksign}
\end{equation}
The sign \(\varepsilon_k\) is additional discrete information.  If the
line element is required to have signature \((-+++)\), it must be chosen so
that
\begin{equation}
 \frac{e^{2k}}{\Re f}>0.
 \label{eq:ksignaturechoice}
\end{equation}

Minkowski space corresponds to
\begin{equation}
 f_0=1,
 \qquad e^{2U_0}=e^{2k_0}=1,
 \qquad a_0=b_0=0.
 \label{eq:MinkowskiErnst}
\end{equation}

\subsection{The spectral covering and the linear problem}

Put
\[
 z=\rho+i\zeta,
 \qquad \bar z=\rho-i\zeta,
\]
and introduce the spectral relation
\begin{equation}
 \lambda^2=\frac{K-i\bar z}{K+iz}.
 \label{eq:spectralrelation}
\end{equation}
The branch values \(K=i\bar z\) and \(K=-iz\) are excluded at first.  The
total spectral space over \(\mathscr O\) is
\begin{equation}
 \mathscr S_{\mathscr O}:=
 \left\{(\rho,\zeta,K,\lambda):
 \begin{array}{l}
 (\rho,\zeta)\in\mathscr O,\ K\in\mathbb C,\ \lambda\in\mathbb C^\times,\\
 K\ne i\bar z,\ K\ne-iz,\quad
 \lambda^2=(K-i\bar z)/(K+iz)
 \end{array}\right\}.
 \label{eq:totalS}
\end{equation}
Both \(K\) and \(\lambda\) are single-valued functions on this space.

For fixed \(K\), define
\[
 \mathscr O_K:=\{(\rho,\zeta)\in\mathscr O:
                 K\ne i\bar z,\ K\ne-iz\}
\]
and
\begin{equation}
 \mathscr S_K:=\{(\rho,\zeta,\lambda):
 (\rho,\zeta)\in\mathscr O_K,\ 
 \lambda^2=(K-i\bar z)/(K+iz)\}.
 \label{eq:fixedK}
\end{equation}
The projection \(\mathscr S_K\to\mathscr O_K\) is an unbranched double
cover.  Consequently the vector fields \(\partial_z,\partial_{\bar z}\) on
the base have unique lifts to \(\mathscr S_K\); these lifted derivatives are
used whenever a function on \(\mathscr S_{\mathscr O}\) is differentiated.
Differentiating \eqref{eq:spectralrelation} at fixed \(K\) gives
\begin{equation}
 \lambda_{,z}=\frac{\lambda}{4\rho}(\lambda^2-1),
 \qquad
 \lambda_{,\bar z}=\frac{\lambda^2-1}{4\rho\lambda}.
 \label{eq:lambdaderivatives}
\end{equation}

For fixed \(x=(\rho,\zeta)\), the spectral fibre is
\[
 \mathscr S_x:=\{(K,\lambda):(x,K,\lambda)\in\mathscr S_{\mathscr O}\}.
\]
After adjoining the omitted points, its compactification is a Riemann sphere,
\begin{equation}
 \widehat{\mathscr S}_x\simeq\mathbb {CP}^1_\lambda,
 \qquad
 K=\frac{i\bar z+iz\lambda^2}{1-\lambda^2}.
 \label{eq:fibrecompactification}
\end{equation}
Here \(\lambda\) is a global meromorphic parameter on the compactified
fibre.  The points \(\lambda=0,\infty\) lie above the two branch values,
whereas \(\lambda=\pm1\) are the two points above \(K=\infty\).

For a compact Riemann surface \(X\), let \(\mathcal M(X)\) denote its field
of meromorphic functions.  Since \(\rho>0\) and the two branch points have
been removed, the total spectral space \(\mathscr S_{\mathscr O}\) is a
smooth manifold.  Its smooth structure may be described by the local
fixed-\(K\) covers, or equivalently by using \(\lambda\) as a local
coordinate along the fibres.

Let \(r\in\mathbb N_0\).  A scalar family
\[
 F:\mathscr S_{\mathscr O}\longrightarrow\mathbb {CP}^1
\]
is called a \(C^r\) fibrewise meromorphic family if it is a \(C^r\) map on
the total spectral space and, for every \(x\in\mathscr O\), its restriction
to the fibre extends meromorphically to the compactified fibre:
\begin{equation}
 F(x,\cdot)\in\mathcal M(\widehat{\mathscr S}_x).
 \label{eq:fibremeromorphic}
\end{equation}
A matrix family \(F=(F_{ab})\) is called \(C^r\) fibrewise meromorphic if
every entry \(F_{ab}\) has this property.  Equivalently,
\[
 F(x,\cdot)\in
 \operatorname{Mat}_{2\times2}\bigl(\mathcal M(\widehat{\mathscr S}_x)\bigr)
\]
for every \(x\in\mathscr O\), and its entries depend jointly
\(C^r\)-smoothly on the points of the total spectral space, with a pole
interpreted as the value \(\infty\in\mathbb {CP}^1\).  Thus this definition
contains two separate requirements: meromorphic dependence along every
compactified spectral fibre and joint \(C^r\) dependence on the total
spectral space.

For a \(C^2\) Ernst potential with \(f+\bar f\ne0\), set
\begin{equation}
 M=\frac{f_{,z}}{f+\bar f},
 \qquad
 N=\frac{\bar f_{,z}}{f+\bar f}.
 \label{eq:MN}
\end{equation}
The associated linear problem is
\begin{align}
 \Phi_{,z}&=\mathscr A(\lambda)\Phi,
 &\mathscr A(\lambda)&=
 \begin{pmatrix}N&0\\0&M\end{pmatrix}
 +\lambda\begin{pmatrix}0&N\\M&0\end{pmatrix},
 \label{eq:LPz}\\
 \Phi_{,\bar z}&=\mathscr B(\lambda)\Phi,
 &\mathscr B(\lambda)&=
 \begin{pmatrix}\bar M&0\\0&\bar N\end{pmatrix}
 +\lambda^{-1}\begin{pmatrix}0&\bar M\\\bar N&0\end{pmatrix}.
 \label{eq:LPzb}
\end{align}
For every fixed \(x\), the coefficient matrices are fibrewise meromorphic:
\begin{equation}
 \mathscr A(x,\cdot),\mathscr B(x,\cdot)
 \in\operatorname{Mat}_{2\times2}
 \bigl(\mathcal M(\widehat{\mathscr S}_x)\bigr).
 \label{eq:ABmeromorphic}
\end{equation}
Indeed, \(\mathscr A\) is affine in \(\lambda\), whereas \(\mathscr B\) is
affine in \(\lambda^{-1}\).  Thus the LP is naturally an equation between
fibrewise meromorphic matrix families.

A solution of the LP is a \(C^2\) fibrewise meromorphic matrix family
\(\Phi\) satisfying \eqref{eq:LPz}--\eqref{eq:LPzb} away from its polar
set.  The derivatives are taken with respect to the lifted vector fields on
the fixed-\(K\) covers.  Since both sides extend fibrewise meromorphically,
the two equations are understood as meromorphic identities on every
compactified fibre.

The LP is called integrable when
\begin{equation}
 \mathscr A_{,\bar z}-\mathscr B_{,z}+[\mathscr A,\mathscr B]=0.
 \label{eq:zerocurvature}
\end{equation}
Comparison of the coefficients of the different powers of \(\lambda\),
with the derivatives taken at fixed \(K\), gives
\begin{align}
 N_{,\bar z}+\frac{N+\bar M}{4\rho}
       +N(\bar N-\bar M)&=0,
 \label{eq:coefficientN}\\
 M_{,\bar z}+\frac{M+\bar N}{4\rho}
       +M(\bar M-\bar N)&=0,
 \label{eq:coefficientM}
\end{align}
together with their conjugates and the constant-term relations
\begin{align}
 N_{,\bar z}-\bar M_{,z}+N\bar N-M\bar M&=0,
 \label{eq:constantN}\\
 M_{,\bar z}-\bar N_{,z}+M\bar M-N\bar N&=0.
 \label{eq:constantM}
\end{align}
Substituting \eqref{eq:MN} into these relations, and using their conjugates,
reduces the zero-curvature equation to the Ernst equation
\eqref{eq:Ernst}.  Conversely, an Ernst potential gives coefficient matrices
that satisfy \eqref{eq:zerocurvature}.  Thus the Ernst equation is
equivalent to the integrability of the LP.

\medskip
\noindent\textit{Symmetries.}
Let
\begin{equation}
 \sigma_1=\begin{pmatrix}0&1\\1&0\end{pmatrix},
 \qquad
 \sigma_3=\begin{pmatrix}1&0\\0&-1\end{pmatrix}.
 \label{eq:Pauli}
\end{equation}
The coefficient matrices satisfy
\begin{equation}
 \mathscr A(-\lambda)=\sigma_3\mathscr A(\lambda)\sigma_3,
 \qquad
 \mathscr B(-\lambda)=\sigma_3\mathscr B(\lambda)\sigma_3.
 \label{eq:coefficientSheetSymmetry}
\end{equation}
If
\begin{equation}
 Y(\lambda)=\begin{pmatrix}\psi(\lambda)\\\chi(\lambda)\end{pmatrix}
 \label{eq:columnY}
\end{equation}
is a column solution, then \(\sigma_3Y(-\lambda)\) is another column
solution.  Hence one may form
\begin{equation}
 \Phi(\lambda)=
 \bigl(Y(\lambda),\sigma_3Y(-\lambda)\bigr)
 =\begin{pmatrix}
 \psi(\lambda)&\psi(-\lambda)\\
 \chi(\lambda)&-\chi(-\lambda)
 \end{pmatrix}.
 \label{eq:PhiColumns}
\end{equation}
By construction,
\begin{equation}
 \Phi(-\lambda)=\sigma_3\Phi(\lambda)\sigma_1.
 \label{eq:sheetsymmetry}
\end{equation}
Thus the sheet symmetry follows from the coefficient symmetry and the chosen
construction of the second column.

There is also an anti-linear symmetry.  For a column solution \(Y\), define
\begin{equation}
 (\mathcal RY)(\rho,\zeta,K,\lambda)
 :=\sigma_1\overline{Y\left(\rho,\zeta,\bar K,1/\bar\lambda\right)}.
 \label{eq:realityoperator}
\end{equation}
The spectral relation implies that the argument on the right again lies on
the corresponding spectral covering.  Direct substitution in the LP shows
that \(\mathcal RY\) is a column solution whenever \(Y\) is, and
\begin{equation}
 \mathcal R^2Y=Y.
 \label{eq:Rinvolution}
\end{equation}
We may therefore choose a reality-adapted column satisfying
\(\mathcal RY=Y\).  In components this means
\begin{align}
 \chi(\rho,\zeta,K,\lambda)
 &=\overline{\psi\left(\rho,\zeta,\bar K,1/\bar\lambda\right)},
 \label{eq:realitypsi}\\
 \psi(\rho,\zeta,K,\lambda)
 &=\overline{\chi\left(\rho,\zeta,\bar K,1/\bar\lambda\right)}.
 \label{eq:realitychi}
\end{align}

Right multiplication by a matrix depending on \(K\) but not on
\((\rho,\zeta)\) preserves the LP.  Subject to compatibility with the two
symmetries, this right-gauge freedom can be used to impose the following
normalization whenever the required spectral-infinity limits are
non-degenerate.

A \(C^r\) fibrewise meromorphic matrix family is called normalized if it
satisfies the sheet symmetry \eqref{eq:sheetsymmetry}, the reality symmetry
\eqref{eq:realitypsi}--\eqref{eq:realitychi}, and the lower
spectral-infinity normalization below.  A normalized LP solution is a
normalized matrix that also satisfies the LP.

\medskip
\noindent\textit{Spectral infinity and recovery of the Ernst potential.}
All spectral-infinity limits are taken with \((\rho,\zeta)\) fixed.  The
lower spectral infinity is \(K\to\infty,\lambda\to-1\), whereas the upper
spectral infinity is \(K\to\infty,\lambda\to1\).  The lower normalization
is
\begin{equation}
 \lim_{\substack{K\to\infty\\\lambda\to-1}}\psi=1,
 \qquad
 \lim_{\substack{K\to\infty\\\lambda\to-1}}\chi=1.
 \label{eq:lowernormalization}
\end{equation}
The reality relation makes the two equalities equivalent; both are displayed
to make the limiting first column explicit:
\begin{equation}
 \lim_{\substack{K\to\infty\\\lambda\to-1}}
 \Phi(\rho,\zeta,K,\lambda)
 \begin{pmatrix}1\\0\end{pmatrix}
 =\begin{pmatrix}1\\1\end{pmatrix}.
 \label{eq:lowerfirstcolumn}
\end{equation}

Suppose that a normalized LP solution has an upper spectral-infinity limit,
and define
\begin{equation}
 \widehat f:=
 \lim_{\substack{K\to\infty\\\lambda\to1}}
 \chi(\rho,\zeta,K,\lambda).
 \label{eq:recoverf}
\end{equation}
The notation \(\chi(1)\) will occasionally abbreviate this fibrewise limit.
Taking the upper limit in the LP gives
\begin{equation}
 \widehat f_{,z}=M(\widehat f+\overline{\widehat f}),
 \qquad
 \overline{\widehat f}_{,z}=N(\widehat f+\overline{\widehat f}),
 \label{eq:fhatLP}
\end{equation}
together with the conjugate \(\bar z\)-relations.

Suppose first that the LP coefficients were constructed from an Ernst
potential \(f\) by \eqref{eq:MN}, and assume that both real parts are
nonzero.  Put
\begin{equation}
 r:=\frac{\widehat f+\overline{\widehat f}}{f+\bar f}
   =\frac{\Re\widehat f}{\Re f}.
 \label{eq:rdefinition}
\end{equation}
Equations \eqref{eq:MN} and \eqref{eq:fhatLP}, including their conjugates,
give
\begin{equation}
 d\widehat f=r\,df.
 \label{eq:dfhat}
\end{equation}
Taking real parts yields
\[
 d(\Re\widehat f)=r\,d(\Re f).
\]
On the other hand, \(\Re\widehat f=r\Re f\), so
\[
 d(\Re\widehat f)=r\,d(\Re f)+(\Re f)\,dr.
\]
Since \(\Re f\ne0\), comparison gives \(dr=0\).  Thus \(r=c\) is a
nonzero real constant, and \eqref{eq:dfhat} gives
\(d(\widehat f-cf)=0\).  The real part of this constant vanishes, hence
\begin{equation}
 \widehat f=cf+iC,
 \qquad c\in\mathbb R^\times,\quad C\in\mathbb R.
 \label{eq:fgauge}
\end{equation}
This is precisely the freedom \(f\mapsto cf+iC\), which leaves \(M,N\)
unchanged.

Conversely, suppose an integrable LP of the form
\eqref{eq:LPz}--\eqref{eq:LPzb} has a normalized solution for which
\eqref{eq:recoverf} exists and has nonzero real part.  Equation
\eqref{eq:fhatLP} then gives
\begin{equation}
 M=\frac{\widehat f_{,z}}
        {\widehat f+\overline{\widehat f}},
 \qquad
 N=\frac{\overline{\widehat f}_{,z}}
        {\widehat f+\overline{\widehat f}}.
 \label{eq:recoverMN}
\end{equation}
Substitution in the zero-curvature relations shows that \(\widehat f\)
satisfies the Ernst equation.

Finally, \(f_0=1\) corresponds to
\begin{equation}
 \Phi_0=\begin{pmatrix}1&1\\1&-1\end{pmatrix}.
 \label{eq:Minkowskiseed}
\end{equation}

\subsection{Polynomial transformation matrices}

Let \(\Phi_0\) be a normalized solution of a seed LP with coefficient
matrices \(\mathscr A_0,\mathscr B_0\).  Following Kramer
\cite{Kramer1986}, fix an integer \(N\geq1\) and consider
\begin{equation}
 T(\rho,\zeta,K,\lambda)=
 \left(\frac{K+iz}{K}\right)^N
 \sum_{\nu=0}^{2N}Y_\nu(\rho,\zeta)\lambda^\nu,
 \label{eq:Trepresentation}
\end{equation}
where
\begin{equation}
 Y_\nu\in C^2\bigl(\mathscr O,
 \operatorname{Mat}_{2\times2}(\mathbb C)\bigr),
 \qquad \nu=0,\ldots,2N.
 \label{eq:Yregularity}
\end{equation}
For fixed \(x=(\rho,\zeta)\), the scalar prefactor in
\eqref{eq:Trepresentation} is meromorphic on
\(\widehat{\mathscr S}_x\), while the remaining factor is a polynomial in
the global meromorphic parameter \(\lambda\).  Consequently
\begin{equation}
 T(x,\cdot)\in
 \operatorname{Mat}_{2\times2}
 \bigl(\mathcal M(\widehat{\mathscr S}_x)\bigr),
 \label{eq:Tmeromorphic}
\end{equation}
so \(T\) is a \(C^2\) fibrewise meromorphic matrix family.

We put
\begin{equation}
 \Phi=T\Phi_0.
 \label{eq:transformedPhi}
\end{equation}
No Darboux equation is assumed at this stage.

For a constant \(K_i\in\mathbb C\), define
\[
 \mathscr O_i:=\{(\rho,\zeta)\in\mathscr O:
 K_i\ne i\bar z,\ K_i\ne-iz\}
\]
and the fixed-\(K_i\) slice
\begin{equation}
 \mathscr S_i:=\mathscr S_{\mathscr O}\cap\{K=K_i\}
 =\left\{(\rho,\zeta,\lambda):
 (\rho,\zeta)\in\mathscr O_i,\quad
 \lambda^2=\frac{K_i-i\bar z}{K_i+iz}\right\}.
 \label{eq:Si}
\end{equation}
It is an unbranched double cover of \(\mathscr O_i\); its connectedness is
determined by the square-root monodromy around the branch values removed
from \(\mathscr O\).

\begin{proposition}[Darboux criterion]
\label{prop:Darbouxcriterion}
Suppose that \(\Phi_0\) is a normalized \(C^2\) fibrewise meromorphic
solution of the seed LP.  Let \(T\) be given by
\eqref{eq:Trepresentation}, and suppose that:
\begin{enumerate}
 \item \(\Phi=T\Phi_0\) is a normalized \(C^2\) fibrewise meromorphic
 matrix family;
 \item on every compactified fibre, each zero \(q\) of \(\det T\) away from
 \(\lambda=0,\infty\) is simple and satisfies \(\operatorname{rank}T(q)=1\);
 \item there exist an integer \(r\geq1\), constants
 \[
 K_i\in\mathbb C,
 \qquad i\in\{1,\ldots,r\},
 \]
 connected components \(\mathcal Z_i\) of the zero set of \(\det T\)
 satisfying
 \[
 \mathcal Z_i\subset\mathscr S_i,
 \qquad i\in\{1,\ldots,r\},
 \]
 and constant vectors
 \[
 b_i\in\mathbb C^2\setminus\{0\},
 \qquad i\in\{1,\ldots,r\},
 \]
 such that
 \[
 T\Phi_0b_i=0
 \quad\text{on }\mathcal Z_i,
 \qquad i\in\{1,\ldots,r\};
 \]
 moreover, \(\mathcal Z_1,\ldots,\mathcal Z_r\) exhaust all zeros of
 \(\det T\) away from \(\lambda=0,\infty\);
 \item the fibrewise endpoint matrices
 \[
 T_0=\lim_{\lambda\to0}T,
 \qquad T_\infty=\lim_{\lambda\to\infty}T
 \]
 are invertible.
\end{enumerate}
Then there exist \(M,N\in C^1(\mathscr O)\) such that \(\Phi\) satisfies
the Ernst LP \eqref{eq:LPz}--\eqref{eq:LPzb}.  Consequently \(T\) obeys
the corresponding Darboux equations.  If the upper spectral-infinity limit
of \(\chi\) exists and has nonzero real part, it defines an Ernst potential.
\end{proposition}

\begin{proof}
Where \(T\) is invertible, define
\begin{equation}
 \mathscr A_T=(T_{,z}+T\mathscr A_0)T^{-1},
 \qquad
 \mathscr B_T=(T_{,\bar z}+T\mathscr B_0)T^{-1}.
 \label{eq:logderivatives}
\end{equation}
Then \(\Phi_{,z}=\mathscr A_T\Phi\) and
\(\Phi_{,\bar z}=\mathscr B_T\Phi\).

Let \(q\in\mathcal Z_i\) be a zero of \(\det T\) away from
\(\lambda=0,\infty\), and set \(B_i=\Phi_0b_i\).  Since \(b_i\) is
constant and \(\Phi_0\) solves the seed LP,
\[
 (B_i)_{,z}=\mathscr A_0B_i,
 \qquad (B_i)_{,\bar z}=\mathscr B_0B_i.
\]
Differentiating \(TB_i=0\) along \(\mathscr S_i\) gives
\begin{equation}
 (T_{,z}+T\mathscr A_0)B_i=0,
 \qquad
 (T_{,\bar z}+T\mathscr B_0)B_i=0.
 \label{eq:numeratorannihilates}
\end{equation}

We now show that the apparent pole in \eqref{eq:logderivatives} is
removable.  Let \(w^T\) span the left kernel of \(T(q)\).  Since
\(\operatorname{rank}T(q)=1\), the adjugate has the form
\[
 \operatorname{adj}T(q)=\gamma B_i(q)w^T,
 \qquad \gamma\ne0.
\]
With \(\tau=\lambda-\lambda(q)\), simplicity of the zero gives
\[
 \det T=c\tau+O(\tau^2),\qquad c\ne0.
\]
Using \(T^{-1}=\operatorname{adj}T/\det T\), we obtain
\[
 T^{-1}=\frac{\gamma}{c}\frac{B_i(q)w^T}{\tau}+O(1).
\]
Equation \eqref{eq:numeratorannihilates} therefore yields
\[
 \mathscr A_T=(T_{,z}+T\mathscr A_0)T^{-1}
 =\frac{\gamma}{c}
   \frac{(T_{,z}+T\mathscr A_0)(q)B_i(q)w^T}{\tau}
   +O(1)=O(1).
\]
The same argument applies to \(\mathscr B_T\).  Hence all poles arising
from the prescribed zeros of \(\det T\) are removable.

It remains to account for the common scalar factor.  Write
\[
 T=s\widehat T,
 \qquad s=\left(\frac{K+iz}{K}\right)^N.
\]
At fixed \(K\),
\[
 \frac{s_{,z}}s=\frac{Ni}{K+iz},
 \qquad \frac{s_{,\bar z}}s=0.
\]
Thus the factor \(K^{-N}\) does not create a pole of the logarithmic
derivatives at \(K=0\); possible poles arising from
\(\widehat T^{-1}\) have already been covered by the preceding zero
analysis.

From \eqref{eq:fibrecompactification},
\begin{equation}
 \frac{K+iz}{K}=\frac{2\rho}{\bar z+z\lambda^2}.
 \label{eq:scalarfactorlambda}
\end{equation}
Consequently \eqref{eq:Trepresentation} gives automatically
\begin{align}
 T&=T_0+T_1\lambda+O(\lambda^2),
 &T_0&=\left(\frac{2\rho}{\bar z}\right)^NY_0,
 &&\lambda\to0,
 \label{eq:Tzeroexpansion}\\
 T&=T_\infty+T_{-1}\lambda^{-1}+O(\lambda^{-2}),
 &T_\infty&=\left(\frac{2\rho}{z}\right)^NY_{2N},
 &&\lambda\to\infty.
 \label{eq:Tinfexpansion}
\end{align}
Using \eqref{eq:lambdaderivatives} to differentiate these expansions at
fixed \(K\), one obtains
\[
 \begin{array}{c|cc}
 &\lambda\to0&\lambda\to\infty\\ \hline
 T_{,z}&O(1)&O(\lambda)\\
 T_{,\bar z}&O(\lambda^{-1})&O(1).
 \end{array}
\]
Since the seed coefficients have the corresponding growth
\[
 \mathscr A_0=O(1),\quad \mathscr B_0=O(\lambda^{-1})
 \quad(\lambda\to0),
\]
and
\[
 \mathscr A_0=O(\lambda),\quad \mathscr B_0=O(1)
 \quad(\lambda\to\infty),
\]
invertibility of \(T_0,T_\infty\) gives
\begin{align}
 \mathscr A_T&=O(1),&\mathscr B_T&=O(\lambda^{-1})
 &&(\lambda\to0),\label{eq:growthzero}\\
 \mathscr A_T&=O(\lambda),&\mathscr B_T&=O(1)
 &&(\lambda\to\infty).\label{eq:growthinf}
\end{align}

Thus, on each compactified fibre, \(\mathscr A_T\) has no finite poles and
has at most a simple pole at infinity.  Likewise \(\mathscr B_T\) has only
a possible simple pole at zero.  The standard meromorphic-function argument
on \(\mathbb {CP}^1\) therefore gives
\begin{equation}
 \mathscr A_T=P+\lambda Q,
 \qquad
 \mathscr B_T=\widetilde P+\lambda^{-1}\widetilde Q,
 \label{eq:spectraldegree}
\end{equation}
with \(C^1\) coefficient matrices on \(\mathscr O\).

The sheet symmetry of \(\Phi\) implies
\[
 \mathscr A_T(-\lambda)=\sigma_3\mathscr A_T(\lambda)\sigma_3.
\]
Hence \(P\) is diagonal and \(Q\) is off-diagonal.  Write
\[
 P=\begin{pmatrix}N&0\\0&M\end{pmatrix},
 \qquad Q=\begin{pmatrix}0&q\\p&0\end{pmatrix}.
\]
Let \(e=(1,1)^T\).  The lower spectral-infinity normalization gives a
constant limiting first column, and therefore
\[
 \mathscr A_T(-1)e=0.
\]
It follows that \(q=N\) and \(p=M\), so
\[
 \mathscr A_T=
 \begin{pmatrix}N&0\\0&M\end{pmatrix}
 +\lambda\begin{pmatrix}0&N\\M&0\end{pmatrix}.
\]
Applying the reality symmetry gives
\[
 \mathscr B_T=
 \begin{pmatrix}\bar M&0\\0&\bar N\end{pmatrix}
 +\lambda^{-1}\begin{pmatrix}0&\bar M\\\bar N&0\end{pmatrix}.
\]
This proves the two LP equations.  Finally, \(\Phi=T\Phi_0\) implies
\[
 T_{,z}=\mathscr A_TT-T\mathscr A_0,
 \qquad
 T_{,\bar z}=\mathscr B_TT-T\mathscr B_0,
\]
which are the Darboux equations.  The assertion concerning the recovered
Ernst potential follows from the preceding subsection.
\end{proof}

\section{The multi-horizon determinant family}

\subsection{The axis configuration}

Fix an integer
\[
 N\geq1.
\]
Let
\[
 K_j\in\mathbb R,
 \qquad j\in\{1,\ldots,2N\},
\]
be \(2N\) distinct real constants ordered by
\begin{equation}
 K_1>K_2>\cdots>K_{2N}.
 \label{eq:Kordering}
\end{equation}
We consider the boundary \(\rho=0\) of the Weyl--Papapetrou half-plane,
parametrized by \(\zeta\in\mathbb R\).  The constants \(K_j\) divide this
boundary into \(N\) finite candidate horizon intervals and \(N+1\) open
axis intervals.

For every \(j\in\{1,\ldots,N\}\), define the \(j\)-th candidate horizon
interval by
\begin{equation}
 \mathcal H^j
 :=\{(0,\zeta):K_{2j}<\zeta<K_{2j-1}\}.
 \label{eq:horizonintervals}
\end{equation}
Its two endpoints are \((0,K_{2j-1})\) and \((0,K_{2j})\), and its
coordinate length is
\begin{equation}
 \ell_j:=K_{2j-1}-K_{2j}>0,
 \qquad j\in\{1,\ldots,N\}.
 \label{eq:horizonlengths}
\end{equation}

The upper semi-infinite axis interval is
\begin{equation}
 \mathcal A^0:=\{(0,\zeta):\zeta>K_1\}.
 \label{eq:upperaxis}
\end{equation}
For every \(j\in\{1,\ldots,N-1\}\), the finite axis interval between
\(\mathcal H^j\) and \(\mathcal H^{j+1}\) is
\begin{equation}
 \mathcal A^j:=\{(0,\zeta):K_{2j+1}<\zeta<K_{2j}\}.
 \label{eq:middleaxes}
\end{equation}
Finally, the lower semi-infinite axis interval is
\begin{equation}
 \mathcal A^N:=\{(0,\zeta):\zeta<K_{2N}\}.
 \label{eq:loweraxis}
\end{equation}
In accordance with the notation used in the two-horizon problem, we also
write
\begin{equation}
 \mathcal A^0=\mathcal A^+,
 \qquad
 \mathcal A^N=\mathcal A^-.
 \label{eq:outeraxisnotation}
\end{equation}

At this stage, the terminology ``horizon interval'' designates only the
prescribed boundary decomposition.  We do not yet claim that a set
\(\mathcal H^j\) extends to a regular Killing horizon of a reconstructed
spacetime.  Likewise, the constants \(K_j\) are presently algebraic
boundary data; their geometric interpretation as horizon endpoints requires
additional regularity assumptions.

We shall also need a spectral space extending the construction of
Section~2.3 to the closed orbit half-plane.  Put
\[
 \overline{\mathscr P}:=[0,\infty)\times\mathbb R
\]
and let \(\mathscr O\subset\overline{\mathscr P}\) be relatively open.  Thus
\(\mathscr O\) carries the subspace topology inherited from \(\mathbb R^2\).
With \(z=\rho+i\zeta\), define
\begin{equation}
 \begin{split}
 \mathscr S_{\mathscr O}:=
 \biggl\{(\rho,\zeta,K,\lambda)
 &\in\mathscr O\times\mathbb C\times\mathbb C^\times:
 K\ne i\bar z,\quad K\ne-iz,\\
 &\lambda^2=\frac{K-i\bar z}{K+iz}\biggr\},
 \end{split}
 \label{eq:closedtotalspectralspace}
\end{equation}
where \(\mathbb C^\times:=\mathbb C\setminus\{0\}\).  If
\(\mathscr O\subset(0,\infty)\times\mathbb R\), this is exactly the space
introduced in Section~2.3.  The exclusion of the branch values ensures that
\(\lambda\) is finite and nonzero.

On \(\mathscr S_{\mathscr O}\), introduce
\begin{equation}
 \mu:=\lambda(K+iz)=\lambda(K-\zeta+i\rho).
 \label{eq:defmu}
\end{equation}
The spectral relation gives
\begin{equation}
 \mu^2=(K-i\bar z)(K+iz)=(K-\zeta)^2+\rho^2.
 \label{eq:muspectralrelation}
\end{equation}
Let \([X:Y:W]\) denote homogeneous coordinates on \(\mathbb {CP}^2\).  We
define the projective completion of the total spectral space by
\begin{equation}
 \widehat{\mathscr S}_{\mathscr O}
 :=\left\{(\rho,\zeta,[X:Y:W])\in
 \mathscr O\times\mathbb {CP}^2:
 Y^2=(X-\zeta W)^2+\rho^2W^2\right\}.
 \label{eq:projectivetotalspace}
\end{equation}
It is given the subspace topology inherited from
\(\mathscr O\times\mathbb {CP}^2\).  The map
\begin{equation}
 (\rho,\zeta,K,\lambda)
 \longmapsto
 (\rho,\zeta,[K:\mu:1])
 \label{eq:spectralembedding}
\end{equation}
embeds \(\mathscr S_{\mathscr O}\) into
\(\widehat{\mathscr S}_{\mathscr O}\).  Its image is the part on which
\(W\ne0\) and \(Y\ne0\).  Thus the projective model restores the branch
points \(Y=0\), adds the two points above \(K=\infty\), and, on the boundary
\(\rho=0\), also adds the point at which the two limiting components meet.

For a fixed base point \(x=(\rho,\zeta)\in\mathscr O\), define
\begin{equation}
 \begin{split}
 \mathscr S_x=\mathscr S_{(\rho,\zeta)}:=
 \biggl\{(K,\lambda)&\in\mathbb C\times\mathbb C^\times:
 K\ne i\bar z,\quad K\ne-iz,\\
 &\lambda^2=\frac{K-i\bar z}{K+iz}\biggr\}
 \end{split}
 \label{eq:closedbasefibre}
\end{equation}
and
\begin{equation}
 \widehat{\mathscr S}_x
 =\widehat{\mathscr S}_{(\rho,\zeta)}
 :=\left\{[X:Y:W]\in\mathbb {CP}^2:
 Y^2=(X-\zeta W)^2+\rho^2W^2\right\}.
 \label{eq:projectivebasefibre}
\end{equation}
For \(\rho>0\), the latter is a smooth projective conic and hence isomorphic
to \(\mathbb {CP}^1\).  At \(W=0\), its two points at spectral infinity are
\begin{equation}
 \mathfrak p_+=[1:1:0],
 \qquad
 \mathfrak p_-=[1:-1:0].
 \label{eq:projectiveinfinities}
\end{equation}

For \(\rho=0\), the fibre equation factors as
\begin{equation}
 \bigl(Y-(X-\zeta W)\bigr)
 \bigl(Y+(X-\zeta W)\bigr)=0.
 \label{eq:nodalfactorization}
\end{equation}
Consequently
\begin{equation}
 \widehat{\mathscr S}_{(0,\zeta)}
 =L_+^\zeta\cup L_-^\zeta,
 \qquad
 L_\pm^\zeta:=\{Y=\pm(X-\zeta W)\}.
 \label{eq:nodalfibre}
\end{equation}
The two projective lines meet at
\begin{equation}
 \mathfrak n_\zeta:=[\zeta:0:1].
 \label{eq:node}
\end{equation}
This point is a node: in the affine coordinates
\(u=K-\zeta\), \(y=\mu\), the boundary fibre has equation
\((y-u)(y+u)=0\), so two smooth local branches with distinct tangent
directions meet there.

For a fixed \(K_0\in\mathbb C\), put
\begin{equation}
 \mathscr O_{K_0}:=\left\{(\rho,\zeta)\in\mathscr O:
 K_0\ne i\bar z,\ K_0\ne-iz\right\}
 \label{eq:closedfixedKbase}
\end{equation}
and define the fixed-\(K_0\) space
\begin{equation}
 \mathscr S_{K_0}:=
 \left\{(\rho,\zeta,\lambda):
 (\rho,\zeta)\in\mathscr O_{K_0},\quad
 \lambda^2=\frac{K_0-i\bar z}{K_0+iz}\right\}.
 \label{eq:closedfixedKspace}
\end{equation}
The projection \(\mathscr S_{K_0}\to\mathscr O_{K_0}\) is an unbranched
double cover.  We shall not introduce a projective completion of this
fixed-\(K_0\) space; the hat notation will be reserved for the total space
and the fixed-base spectral fibres.

The space \(\mathscr S_{\mathscr O}\) remains a smooth manifold with boundary.
Indeed, at \(\rho=0\) its defining relation has \(\lambda=\pm1\), while
\(K=\zeta\) has been excluded, so the implicit-function description remains
regular.  By contrast, the projective total space
\(\widehat{\mathscr S}_{\mathscr O}\) is not naturally a smooth manifold at
the nodal locus
\begin{equation}
 \mathscr N:=\{(0,\zeta,[\zeta:0:1]):(0,\zeta)\in\mathscr O\}.
 \label{eq:nodallocus}
\end{equation}
Near such a point its affine equation is
\(\mu^2-(K-\zeta)^2-\rho^2=0\), whose differential vanishes at the node.
Nevertheless, \(\widehat{\mathscr S}_{\mathscr O}\) is a well-defined
topological space and is smooth away from \(\mathscr N\).

\begin{remark}
For \(\rho>0\), the function \(\lambda\) extends meromorphically to the
smooth compactified fibre \(\widehat{\mathscr S}_{(\rho,\zeta)}\).  In
homogeneous coordinates it may be written as
\begin{equation}
 \lambda
 =\frac{Y}{X-\zeta W+i\rho W}
 =\frac{X-\zeta W-i\rho W}{Y},
 \label{eq:lambdaprojective}
\end{equation}
where the two expressions agree by the conic equation.  Its values at
\(\mathfrak p_+\) and \(\mathfrak p_-\) are \(+1\) and \(-1\), respectively,
and the two branch points correspond to \(\lambda=0\) and \(\lambda=\infty\).

This statement does not extend unchanged to \(\rho=0\).  Away from the node,
\begin{equation}
 \lambda|_{L_+^\zeta}=+1,
 \qquad
 \lambda|_{L_-^\zeta}=-1.
 \label{eq:lambdanodalbranches}
\end{equation}
The two branches therefore give different limiting values at
\(\mathfrak n_\zeta\).  Hence \(\lambda\) does not descend to a continuous,
and therefore not to a regular, function on the entire nodal fibre.  On its
normalization it may instead be regarded as the pair \((1,-1)\).  This is
why \(\mu\), rather than \(\lambda\), is used in the unified projective
model.
\end{remark}

In particular, on the boundary \(\rho=0\), one has
\begin{equation}
 z=i\zeta,
 \qquad
 \bar z=-i\zeta,
 \label{eq:zaxis}
\end{equation}
and hence
\begin{equation}
 \lambda^2=1
 \qquad(K\ne\zeta).
 \label{eq:lambdaaxis}
\end{equation}
Thus the two sheets reduce on every open boundary interval to
\begin{equation}
 \lambda=+1,
 \qquad
 \lambda=-1.
 \label{eq:axissheets}
\end{equation}

\subsection{Integration of the linear problem along the boundary}

Following Neugebauer and Hennig \cite{NeugebauerHennig2009}, we begin with a
hypothetical stationary, axisymmetric, asymptotically flat vacuum spacetime.
More precisely, we assume that the Weyl--Papapetrou expression \(ds^2\)
defines a Lorentzian metric \(g\) of signature \((-+++)\), satisfying
\begin{equation}
 \operatorname{Ric}(g)=0,
 \label{eq:vacuumassumption}
\end{equation}
and possessing the regularity and boundary properties required in the
axis-integration argument.  Let \(f\) denote its Ernst potential, and let
\(\mathscr O\subset\overline{\mathscr P}\) be a relatively open neighbourhood
of the relevant part of the boundary \(\rho=0\).  Define
\begin{equation}
 \mathscr O^\circ:=\mathscr O\cap\{\rho>0\}
 \label{eq:Ointerior}
\end{equation}
and
\begin{equation}
 \mathscr O_*:=
 \mathscr O\setminus
 \left\{(0,K_i):i\in\{1,\ldots,2N\}\right\}.
 \label{eq:Ostar}
\end{equation}

In the interior \(\mathscr O^\circ\), we assume that \(f\) induces an
integrable LP admitting a normalized \(C^2\) fibrewise meromorphic solution
\(\Phi\) on \(\mathscr S_{\mathscr O^\circ}\).  We further assume that
\(\Phi\) extends to a normalized \(C^2\) fibrewise meromorphic matrix family
on the manifold with boundary
\begin{equation}
 \mathscr S_{\mathscr O_*}.
 \label{eq:Phiboundarydomain}
\end{equation}
Here the terminology has the same meaning as in Section~2.3: the entries of
\(\Phi\) depend jointly \(C^2\)-smoothly on the total spectral space.  On a
boundary fibre, fibrewise meromorphicity means meromorphic dependence on
\(K\) separately on the two components \(\lambda=+1\) and \(\lambda=-1\).
Thus differentiability is assumed up to every open axis and candidate
horizon interval, but no differentiability at the endpoints
\((0,K_i)\) is included in this hypothesis.

Finally, we impose the separate endpoint-matching assumption that \(\Phi\)
admits a continuous extension to the projective spectral space
\begin{equation}
 \widehat{\mathscr S}_{\mathscr O}.
 \label{eq:Phicontinuousextension}
\end{equation}
In particular, this extension assigns a single matrix value to every nodal
point lying over an endpoint \((0,K_i)\).  Equivalently, the limits obtained
by approaching that point through the two components of the corresponding
boundary fibre are required to agree.  This is an additional assumption; it
does not follow merely from meromorphic dependence on the two separate
boundary branches.

On every open boundary interval, the equality \(\lambda^2=1\) reduces the LP
to two systems corresponding to \(\lambda=+1\) and \(\lambda=-1\).
Hennig and Neugebauer integrate these reduced systems separately on the axis
intervals
\[
 \mathcal A^0,\ldots,\mathcal A^N
\]
and the candidate horizon intervals
\[
 \mathcal H^1,\ldots,\mathcal H^N.
\]
The resulting integration matrices depend on \(K\) but are constant along
each individual boundary interval.

Further assumptions are imposed on the hypothetical spacetime, including
the prescribed asymptotic behaviour, regularity of the metric potentials on
the open axis and horizon intervals, the horizon boundary conditions, and
sufficient continuity at the endpoints
\[
 (0,K_i),
 \qquad i\in\{1,\ldots,2N\}.
\]
These conditions relate the integration matrices belonging to adjacent
intervals.  Eliminating the intermediate matrices gives algebraic relations
for the spectral data on the upper axis.

The conclusion of this boundary integration is that the upper-axis Ernst
potential must be a rational function
\begin{equation}
 f^+(\zeta):=
 f(0,\zeta)\big|_{\mathcal A^+}
 =\frac{n_N(\zeta)}{d_N(\zeta)},
 \label{eq:rationalaxisdatum}
\end{equation}
where
\begin{align}
 n_N(\zeta)
 &=\zeta^N+n_{N-1}\zeta^{N-1}+\cdots+n_0,
 \label{eq:nN}\\
 d_N(\zeta)
 &=\zeta^N+d_{N-1}\zeta^{N-1}+\cdots+d_0
 \label{eq:dN}
\end{align}
are monic polynomials of degree \(N\).  Their common leading coefficient is
normalized to one by the asymptotic condition
\begin{equation}
 f^+(\zeta)\longrightarrow1
 \qquad(\zeta\to+\infty).
 \label{eq:axisasymptotics}
\end{equation}

Thus the boundary integration does not yet construct an Ernst potential in
the half-plane.  Rather, it gives a necessary algebraic form for the
upper-axis value of any hypothetical solution satisfying the preceding
regularity and matching assumptions.  The subsequent B\"acklund
construction starts from this rational axis datum and produces an explicit
determinant candidate in the interior.  It remains a separate task to
determine which of the properties assumed in deriving
\eqref{eq:rationalaxisdatum} are actually satisfied by that candidate.

\subsection{The Hennig--Neugebauer determinant family}

Hennig and Neugebauer found that the rational upper-axis datum obtained in
the preceding subsection can be reproduced by a special family of
\(2N\)-fold B\"acklund transformations of the Minkowski seed
\cite{NeugebauerHennig2009}.  We shall not review the general
B\"acklund-transformation formalism.  Instead, we write down explicitly the
particular determinant family relevant to the multi-horizon problem.

Let
\begin{equation}
 \alpha_i\in\mathbb C,
 \qquad
 |\alpha_i|=1,
 \qquad
 i\in\{1,\ldots,2N\}.
 \label{eq:alphadata}
\end{equation}
For each \(i\in\{1,\ldots,2N\}\), define
\begin{equation}
 r_i=r_i(\rho,\zeta):=
 \sqrt{(K_i-\zeta)^2+\rho^2}.
 \label{eq:ri}
\end{equation}

Recall the extended spectral variable
\[
 \mu=\lambda(K+iz),
 \qquad
 \mu^2=(K-\zeta)^2+\rho^2.
\]
Define the free spectral column
\begin{equation}
 V_N(\rho,\zeta;K,\lambda)
 :=
 \begin{pmatrix}
 K^N\\
 \mu K^{N-1}\\
 K^{N-1}\\
 \mu K^{N-2}\\
 K^{N-2}\\
 \vdots\\
 \mu K\\
 K\\
 \mu\\
 1
 \end{pmatrix}
 \in\mathbb C^{2N+1}.
 \label{eq:freecolumn}
\end{equation}
Thus the rows alternate between
\[
 K^N,K^{N-1},\ldots,K,1
\]
and
\[
 \mu K^{N-1},\mu K^{N-2},\ldots,\mu K,\mu.
\]
For each \(i\in\{1,\ldots,2N\}\), define the fixed column
\begin{equation}
 C_i(\rho,\zeta)
 :=
 \begin{pmatrix}
 K_i^N\\
 \alpha_i r_iK_i^{N-1}\\
 K_i^{N-1}\\
 \alpha_i r_iK_i^{N-2}\\
 K_i^{N-2}\\
 \vdots\\
 \alpha_i r_iK_i\\
 K_i\\
 \alpha_i r_i\\
 1
 \end{pmatrix}
 \in\mathbb C^{2N+1}.
 \label{eq:fixedcolumns}
\end{equation}
The \((2N+1)\times(2N+1)\) determinant matrix is
\begin{equation}
 \mathcal M_N(\rho,\zeta;K,\lambda)
 :=
 \bigl(V_N(\rho,\zeta;K,\lambda),C_1,\ldots,C_{2N}\bigr),
 \label{eq:detmatrixMN}
\end{equation}
where parentheses denote the matrix having the displayed vectors as its
columns.

Let \(\mathcal M_N^{(-)}\) have the same final \(2N\) columns as
\(\mathcal M_N\), but first column
\begin{equation}
 e_-:=
 \begin{pmatrix}
 1&-1&0&\cdots&0
 \end{pmatrix}^{\mathsf T}
 \in\mathbb C^{2N+1}.
 \label{eq:eminus}
\end{equation}
Set
\begin{equation}
 \Delta_N(\rho,\zeta)
 :=
 \det\mathcal M_N^{(-)}(\rho,\zeta).
 \label{eq:DeltaN}
\end{equation}
At every point at which \(\Delta_N(\rho,\zeta)\ne0\), define
\begin{equation}
 \chi_N(\rho,\zeta;K,\lambda)
 :=
 \frac1{K^N}
 \frac{\det\mathcal M_N(\rho,\zeta;K,\lambda)}
      {\Delta_N(\rho,\zeta)}.
 \label{eq:chiN}
\end{equation}
The factor \(K^{-N}\) is part of the normalization of this special
B\"acklund family.

Define
\begin{equation}
 \mathscr U:=
 \left\{(\rho,\zeta)\in\mathscr O_*:
 \Delta_N(\rho,\zeta)\ne0\right\}.
 \label{eq:Udet}
\end{equation}
On \(\mathscr O_*\), the functions \(r_i\) are smooth.  For a fixed base
point, the numerator of \(\chi_N\) is polynomial in \(K\) and \(\mu\), while
\(K^{-N}\) is meromorphic on the compactified spectral fibre.  Consequently
\begin{equation}
 \chi_N:\mathscr S_{\mathscr U}\longrightarrow\mathbb {CP}^1
 \label{eq:chimeromorphic}
\end{equation}
is a \(C^\infty\) fibrewise meromorphic function.

For a scalar spectral function \(F\), define
\begin{equation}
 F^\sharp(\rho,\zeta,K,\lambda)
 :=
 \overline{
 F\left(\rho,\zeta,\bar K,\frac1{\bar\lambda}\right)}.
 \label{eq:sharp}
\end{equation}
Set
\begin{equation}
 \psi_N:=\chi_N^\sharp
 \label{eq:psiN}
\end{equation}
and define the matrix candidate
\begin{equation}
 \Phi_N(\rho,\zeta;K,\lambda)
 :=
 \begin{pmatrix}
 \psi_N(\lambda)&\psi_N(-\lambda)\\
 \chi_N(\lambda)&-\chi_N(-\lambda)
 \end{pmatrix}.
 \label{eq:PhiN}
\end{equation}
At this stage, these formulas are taken as algebraic definitions.  In
particular, we have not yet proved that \(\Phi_N\) satisfies an LP.

The two spectral-infinity limits explain the normalization.  Dividing the
free column by \(K^N\) gives
\begin{equation}
 \lim_{\substack{K\to\infty\\ \lambda\to-1}}
 K^{-N}V_N=e_-.
 \label{eq:minusfreelimit}
\end{equation}
Consequently,
\begin{equation}
 \lim_{\substack{K\to\infty\\ \lambda\to-1}}
 \chi_N=1.
 \label{eq:minuschi}
\end{equation}
By construction, \(\Phi_N\) satisfies the sheet and reality symmetries of
Section~2.3, while \eqref{eq:minuschi} gives the lower spectral-infinity
normalization.  Hence \(\Phi_N\) is a normalized \(C^\infty\) fibrewise
meromorphic matrix family on \(\mathscr S_{\mathscr U}\).  This assertion
does not yet imply that \(\Phi_N\) satisfies an LP.
Similarly, let \(\mathcal M_N^{(+)}\) have the same final \(2N\) columns
and first column
\begin{equation}
 e_+:=
 \begin{pmatrix}
 1&1&0&\cdots&0
 \end{pmatrix}^{\mathsf T}.
 \label{eq:eplus}
\end{equation}
Then
\begin{equation}
 \lim_{\substack{K\to\infty\\ \lambda\to+1}}
 K^{-N}V_N=e_+,
 \label{eq:plusfreelimit}
\end{equation}
and the associated Ernst candidate is
\begin{equation}
 f(\rho,\zeta)
 :=
 \lim_{\substack{K\to\infty\\ \lambda\to+1}}
 \chi_N(\rho,\zeta;K,\lambda)
 =
 \frac{\det\mathcal M_N^{(+)}(\rho,\zeta)}
      {\det\mathcal M_N^{(-)}(\rho,\zeta)}.
 \label{eq:Ernstdeterminant}
\end{equation}
Thus the same determinant construction gives the matrix candidate
\(\Phi_N\) and the associated Ernst candidate \(f\).

The determinant family was introduced by Hennig and Neugebauer in order to
reproduce the rational upper-axis datum obtained from the boundary
integration.  At the present stage, however, we take the determinant
expressions only as algebraic definitions.  In particular, we have not yet
proved either that
\begin{equation}
 f(0,\zeta)\big|_{\mathcal A^+}
 =
 \frac{n_N(\zeta)}{d_N(\zeta)}
 \label{eq:axisreproductionpending}
\end{equation}
or that \(\Phi_N\) satisfies an Ernst LP.  In the following subsections we
identify explicit non-degeneracy conditions under which the determinant
matrix satisfies the hypotheses of Proposition~2.1.  On that parameter
domain, \(\Phi_N\) is an LP solution and \(f\) is its associated Ernst
potential.

\subsection{Verification of the Darboux criterion}

We now verify directly that the determinant family of the preceding
subsection satisfies Proposition~\ref{prop:Darbouxcriterion} on an explicit
non-degenerate interior domain.  The proof also identifies the precise role
of the scalar factor \(K^{-N}\), the prescribed kernel data, and the two
branch-point matrices.

Let \(q_\ell=q_\ell(\rho,\zeta)\),
\(\ell\in\{1,\ldots,2N+1\}\), denote the signed cofactor of the entry in
row \(\ell\) of the first column of \(\mathcal M_N\).  Thus
\begin{equation}
 \det\mathcal M_N
 =\sum_{\ell=1}^{2N+1}q_\ell(V_N)_\ell,
 \qquad
 \Delta_N=q_1-q_2.
 \label{eq:qcofactorsexpansion}
\end{equation}
For \(\rho>0\), put
\begin{equation}
 A(\lambda):=\frac{K}{K+iz}
 =\frac{\bar z+z\lambda^2}{2\rho},
 \qquad
 B(\lambda):=\frac1{K+iz}
 =\frac{1-\lambda^2}{2i\rho}.
 \label{eq:ABlambda}
\end{equation}
Define the odd and even polynomials
\begin{align}
 \gamma(\lambda)
 &:=\frac{1}{\Delta_N}
 \sum_{j=1}^{N}q_{2j}\lambda
 A(\lambda)^{N-j}B(\lambda)^{j-1}
 =\sum_{j=1}^{N}\gamma_{2j-1}\lambda^{2j-1},
 \label{eq:gammadefinition}\\
 \delta(\lambda)
 &:=\frac{1}{\Delta_N}
 \sum_{j=0}^{N}q_{2j+1}
 A(\lambda)^{N-j}B(\lambda)^j
 =\sum_{j=0}^{N}\delta_{2j}\lambda^{2j}.
 \label{eq:deltadefinition}
\end{align}
The coefficients \(\gamma_{2j-1}\) and \(\delta_{2j}\) are smooth
functions of \((\rho,\zeta)\) wherever \(\rho>0\) and \(\Delta_N\ne0\).

We shall also use the following boundary notation.  Let \(I\) be any one
of the open intervals introduced in Section~3.1.  On \(I\), choose the
signs \(\lambda_i^I\in\{1,-1\}\) so that
\begin{equation}
 r_i(0,\zeta)=\lambda_i^I(K_i-\zeta),
 \qquad
 u_i^I:=\alpha_i\lambda_i^I,
 \qquad i\in\{1,\ldots,2N\}.
 \label{eq:axisuibeforeprop}
\end{equation}

\begin{remark}
The signs \(\lambda_i^I\) in \eqref{eq:axisuibeforeprop} arise from a
single, consistent choice on the real-\(K\) part of the spectral covering.
To make this explicit, consider
\begin{equation}
 \mathscr S'
 :=
 \left\{
 (\rho,\zeta,K,\lambda)\in\mathscr S_{\mathscr O^\circ}:
 K\in\mathbb R
 \right\}.
 \label{eq:realKspectralset}
\end{equation}
For \(\rho>0\) and \(K\in\mathbb R\), put
\begin{equation}
 r(\rho,\zeta,K)
 :=\sqrt{(K-\zeta)^2+\rho^2}.
 \label{eq:realKradialfunction}
\end{equation}
The spectral relation implies that every point of \(\mathscr S'\) belongs
to exactly one of the two disjoint subsets
\begin{equation}
 \mathscr S'_+
 :=
 \left\{
 \lambda=\frac{K-\zeta-i\rho}{r(\rho,\zeta,K)}
 \right\}
 \label{eq:realKplussheet}
\end{equation}
and
\begin{equation}
 \mathscr S'_-
 :=
 \left\{
 \lambda=-\frac{K-\zeta-i\rho}{r(\rho,\zeta,K)}
 \right\}.
 \label{eq:realKminussheet}
\end{equation}
Indeed, the continuous function
\begin{equation}
 \varepsilon
 :=
 \lambda\frac{r(\rho,\zeta,K)}{K-\zeta-i\rho}
 \label{eq:realKsheetindicator}
\end{equation}
takes only the values \(+1\) and \(-1\), and
\begin{equation}
 \mathscr S'_\pm=\varepsilon^{-1}(\pm1).
\end{equation}
Thus \(\mathscr S'_+\) and \(\mathscr S'_-\) are disjoint unions of
connected components of \(\mathscr S'\), and the sheet involution
\begin{equation}
 (K,\lambda)\longmapsto(K,-\lambda)
\end{equation}
interchanges them.

For the boundary calculations in this paper, we work on \(\mathscr S'_+\).
Restricting its defining relation to \(K=K_i\) gives
\begin{equation}
 \lambda_i(\rho,\zeta)
 =
 \frac{K_i-\zeta-i\rho}{r_i(\rho,\zeta)}.
 \label{eq:lambdaiOnRealPlusSheet}
\end{equation}
Hence, on an open boundary interval \(I\),
\begin{equation}
 \lambda_i^I
 =
 \lim_{\substack{\rho\downarrow0\\ \mathscr S'_+}}
 \lambda_i(\rho,\zeta)
 =
 \frac{K_i-\zeta}{|K_i-\zeta|}
 \in\{1,-1\}.
 \label{eq:lambdaiIexplained}
\end{equation}
This limit is independent of the running coordinate \(\zeta\in I\), and
it gives
\begin{equation}
 r_i(0,\zeta)
 =|K_i-\zeta|
 =\lambda_i^I(K_i-\zeta),
\end{equation}
which explains \eqref{eq:axisuibeforeprop}.

Choosing \(\mathscr S'_-\) instead changes all signs simultaneously:
\begin{equation}
 \lambda_i^I\longmapsto-\lambda_i^I,
 \qquad
 u_i^I\longmapsto-u_i^I,
 \qquad i\in\{1,\ldots,2N\}.
\end{equation}
The choice of \(\mathscr S'_+\) is therefore an auxiliary convention used
to perform the boundary determinant calculations.  The final quantities
defined on the \((\rho,\zeta)\)-plane will be independent of this
simultaneous sheet choice.
\end{remark}

Put
\begin{equation}
 c^I=(u_1^I,\ldots,u_{2N}^I),
 \qquad d=(1,\ldots,1),
 \qquad K=\operatorname{diag}(K_1,\ldots,K_{2N}),
 \label{eq:axiscdbeforeprop}
\end{equation}
and define
\begin{equation}
 H_N^I:=
 [c^IK^{N-1},dK^{N-1},\ldots,c^IK,dK,c^I,d],
 \label{eq:HNaxisdeterminant}
\end{equation}
where square brackets denote the determinant having the displayed row
vectors as its rows.

\begin{proposition}[The Hennig--Neugebauer determinant is an LP solution]
\label{prop:HNdeterminantLP}
Assume the data of Section~3.1.  Let \(\mathscr V\) be an open subset of
the interior \(\rho>0\) on which
\begin{equation}
 \Delta_N\ne0
 \label{eq:DeltaNondegenerate}
\end{equation}
and
\begin{equation}
 \delta_0\delta_{2N}\ne0.
 \label{eq:deltaNondegenerate}
\end{equation}
Then the matrix \(\Phi_N\) defined by \eqref{eq:PhiN} satisfies all the
hypotheses of Proposition~\ref{prop:Darbouxcriterion}, with the Minkowski
seed
\begin{equation}
 \Phi_0=\begin{pmatrix}1&1\\1&-1\end{pmatrix}.
 \label{eq:Minkowskiseedagain}
\end{equation}
Consequently, \(\Phi_N\) satisfies an Ernst LP on \(\mathscr V\).  If in
addition the function \(f\) in \eqref{eq:Ernstdeterminant} has nonzero real
part, then it is the Ernst potential associated with that LP.

More explicitly, let \(I\) be one of the open boundary intervals and let
\(\zeta_0\in I\) satisfy \(\zeta_0\ne0\) and
\begin{equation}
 \Delta_N(0,\zeta_0)\ne0.
 \label{eq:axisDeltaNonzero}
\end{equation}
If
\begin{equation}
 H_N^I\ne0,
 \label{eq:axisHNonzero}
\end{equation}
then \eqref{eq:deltaNondegenerate} holds in a sufficiently small
half-neighbourhood of \((0,\zeta_0)\).  Thus \(H_N^I\ne0\) is a sufficient,
but not necessary, boundary criterion for the branch-point hypothesis of
Proposition~\ref{prop:Darbouxcriterion}.
\end{proposition}

\begin{proof}
We divide the verification into five steps.

\medskip
\noindent\emph{Step 1: the transformation matrix.}
Expanding the determinant in \eqref{eq:chiN} along its first column and
using \eqref{eq:ABlambda} gives
\begin{equation}
 \chi_N(\lambda)
 =g_N(\lambda)\bigl(\gamma(\lambda)+\delta(\lambda)\bigr),
 \qquad
 g_N(\lambda):=\left(\frac{K+iz}{K}\right)^N.
 \label{eq:chigammadelta}
\end{equation}
Since \(\gamma\) is odd and \(\delta\) is even, direct multiplication by
\(\Phi_0^{-1}\) yields
\begin{equation}
 T:=\Phi_N\Phi_0^{-1}
 =g_N
 \begin{pmatrix}
  \lambda^{2N}\overline{\delta(\bar\lambda^{-1})}&
  \lambda^{2N}\overline{\gamma(\bar\lambda^{-1})}\\
  \gamma(\lambda)&\delta(\lambda)
 \end{pmatrix}.
 \label{eq:THNexplicit}
\end{equation}
Hence
\begin{equation}
 T=g_N\widehat T,
 \qquad
 \widehat T=\sum_{\nu=0}^{2N}Y_\nu\lambda^\nu,
 \label{eq:Thatpolynomial}
\end{equation}
where
\begin{align}
 Y_{2j}
 &=\begin{pmatrix}
   \overline{\delta_{2N-2j}}&0\\0&\delta_{2j}
  \end{pmatrix},
 &&j\in\{0,\ldots,N\},
 \label{eq:YevenHN}\\
 Y_{2j-1}
 &=\begin{pmatrix}
   0&\overline{\gamma_{2N-2j+1}}\\
   \gamma_{2j-1}&0
  \end{pmatrix},
 &&j\in\{1,\ldots,N\}.
 \label{eq:YoddHN}
\end{align}
This is precisely the representation required in
\eqref{eq:Trepresentation}.  Section~3.3 already established that
\(\Phi_N\) is a normalized \(C^\infty\) fibrewise meromorphic matrix
family on the domain \(\Delta_N\ne0\).  Thus the first hypothesis of
Proposition~\ref{prop:Darbouxcriterion} is satisfied.

\medskip
\noindent\emph{Step 2: the constant kernel vectors.}
For each \(i\in\{1,\ldots,2N\}\), choose on one component of the
fixed-\(K_i\) cover the function \(\lambda_i\) characterized by
\begin{equation}
 \lambda_i(K_i+iz)=r_i.
 \label{eq:lambdairi}
\end{equation}
Let \(V_i^+\) and \(V_i^-\) denote the free column \(V_N\) evaluated at
\(K=K_i\) and \(\lambda=\lambda_i\) and \(-\lambda_i\), respectively.
Only the rows containing \(\mu\) change sign.  Consequently
\begin{equation}
 C_i=\frac{1+\alpha_i}{2}V_i^+
     +\frac{1-\alpha_i}{2}V_i^-.
 \label{eq:Cilinearcombination}
\end{equation}
Replacing the \(i\)-th fixed column by the right-hand side and using
multilinearity of the determinant gives
\begin{equation}
 (1+\alpha_i)\chi_N(\lambda_i)
 +(1-\alpha_i)\chi_N(-\lambda_i)=0.
 \label{eq:chikernelrelation}
\end{equation}
To obtain the corresponding equation for the first row, recall the
definition
\begin{equation}
 \psi_N(\rho,\zeta,K,\lambda)
 =\overline{\chi_N\left(\rho,\zeta,\bar K,
             \frac1{\bar\lambda}\right)}.
 \label{eq:sharprecalled}
\end{equation}
For real \(K_i\), the spectral relation gives
\begin{equation}
 |\lambda_i|^2=|\lambda_i^2|
 =\left|\frac{K_i-\zeta-i\rho}{K_i-\zeta+i\rho}\right|=1.
 \label{eq:lambdaimodulus}
\end{equation}
Hence \(1/\bar\lambda_i=\lambda_i\).  Taking the complex conjugate of
\eqref{eq:chikernelrelation}, using
\(\bar\alpha_i=\alpha_i^{-1}\), and multiplying by \(\alpha_i\), gives
\begin{equation}
 (1+\alpha_i)\psi_N(\lambda_i)
 +(\alpha_i-1)\psi_N(-\lambda_i)=0.
 \label{eq:psikernelrelation}
\end{equation}
Equations \eqref{eq:chikernelrelation} and
\eqref{eq:psikernelrelation} are respectively the second and first rows of
\begin{equation}
 \Phi_N(\lambda_i)b_i^+=0,
 \qquad
 b_i^+:=\frac12
 \begin{pmatrix}1+\alpha_i\\ \alpha_i-1\end{pmatrix}.
 \label{eq:biplus}
\end{equation}
The sheet symmetry gives the second kernel condition
\begin{equation}
 \Phi_N(-\lambda_i)b_i^-=0,
 \qquad
 b_i^-:=\sigma_1b_i^+.
 \label{eq:biminus}
\end{equation}
Both vectors are constant and nonzero.  Since
\(T\Phi_0=\Phi_N\), equations \eqref{eq:biplus}--\eqref{eq:biminus} are
exactly the constant-kernel conditions required in
Proposition~\ref{prop:Darbouxcriterion}.

\medskip
\noindent\emph{Step 3: exhaustion, simplicity, and rank.}
Since \(T=g_N\widehat T\), the matrix kernel relations
\eqref{eq:biplus}--\eqref{eq:biminus} allow the common scalar factor
\(g_N\) to be cancelled as a meromorphic identity.  They therefore imply
directly that
\begin{equation}
 \widehat T(\lambda_i)\Phi_0b_i^+=0,
 \qquad
 \widehat T(-\lambda_i)\Phi_0b_i^-=0.
 \label{eq:Thatkernelconditions}
\end{equation}
For a fixed base point \((\rho,\zeta)\in\mathscr V\), the polynomial
\(\det\widehat T(\lambda)\) has degree at most \(4N\).  The strict
ordering of the real constants \(K_i\) implies that the \(4N\) points
\begin{equation}
 \lambda_1,-\lambda_1,\ldots,
 \lambda_{2N},-\lambda_{2N}
 \label{eq:fourNroots}
\end{equation}
are pairwise distinct and lie away from \(0\) and \(\infty\) when
\(\rho>0\).  Equation \eqref{eq:Thatkernelconditions} shows that all of
them are zeros of \(\det\widehat T\).

Moreover,
\begin{equation}
 \widehat T(0)=Y_0
 =\begin{pmatrix}\overline{\delta_{2N}}&0\\0&\delta_0\end{pmatrix}.
 \label{eq:Thatatzero}
\end{equation}
Condition \eqref{eq:deltaNondegenerate} therefore implies
\(\det\widehat T\not\equiv0\).  A nonzero polynomial of degree at most
\(4N\) with the \(4N\) distinct zeros in \eqref{eq:fourNroots} has degree
exactly \(4N\); these are all its zeros, and every one of them is simple.
At such a zero \(q\), the matrix \(\widehat T(q)\) cannot be the zero
matrix: otherwise every entry would vanish to positive order in a local
coordinate and \(\det\widehat T\) would have a zero of order at least two.
Thus
\begin{equation}
 \operatorname{rank}\widehat T(q)=1.
 \label{eq:rankoneThat}
\end{equation}

The scalar factor \(g_N\) has no finite zeros.  Hence every finite zero of
\(\det T=g_N^2\det\widehat T\) is one of the zeros just listed.  At a pole
of \(g_N\), the simple zero of \(\det\widehat T\) cannot produce a zero of
\(\det T\); at every zero of \(\det T\), the factor \(g_N\) is therefore
finite and nonzero.  The zero remains simple and the matrix has rank one.
Thus all finite zeros required in Proposition~\ref{prop:Darbouxcriterion}
have been exhausted.

\medskip
\noindent\emph{Step 4: the branch-point matrices.}
Equations \eqref{eq:THNexplicit} and \eqref{eq:scalarfactorlambda} give
\begin{align}
 T_0
 &:=\lim_{\lambda\to0}T
 =\left(\frac{2\rho}{\bar z}\right)^N
 \begin{pmatrix}\overline{\delta_{2N}}&0\\0&\delta_0\end{pmatrix},
 \label{eq:T0HN}\\
 T_\infty
 &:=\lim_{\lambda\to\infty}T
 =\left(\frac{2\rho}{z}\right)^N
 \begin{pmatrix}\overline{\delta_0}&0\\0&\delta_{2N}\end{pmatrix}.
 \label{eq:TinfHN}
\end{align}
Since \(\rho>0\), the scalar prefactors are nonzero wherever these
expressions are used.  Thus condition \eqref{eq:deltaNondegenerate} is
equivalent to invertibility of both \(T_0\) and \(T_\infty\).  This proves
the last hypothesis of Proposition~\ref{prop:Darbouxcriterion}.  The first
assertion of the proposition now follows from that criterion.

\medskip
\noindent\emph{Step 5: a boundary criterion for Step 4.}
Fix an open boundary interval \(I\), use the notation
\eqref{eq:axisuibeforeprop}--\eqref{eq:HNaxisdeterminant}, and put
\begin{equation}
 X:=\operatorname{diag}(K_1-\zeta,\ldots,K_{2N}-\zeta).
 \label{eq:Xaxiscriterion}
\end{equation}

From \eqref{eq:deltadefinition},
\begin{align}
 (2\rho)^N\delta_0
 &=\frac1{\Delta_N}\sum_{j=0}^{N}
 q_{2j+1}\bar z^{N-j}i^{-j},
 \label{eq:delta0scaled}\\
 (2\rho)^N\delta_{2N}
 &=\frac1{\Delta_N}\sum_{j=0}^{N}
 q_{2j+1}z^{N-j}i^j.
 \label{eq:delta2Nscaled}
\end{align}
Consequently, for \(\zeta\in I\) and \(\zeta\ne0\),
\begin{equation}
 \lim_{\rho\downarrow0}(T_0)_{22}
 =\frac{1}{\zeta^N\Delta_N(0,\zeta)}
 \sum_{j=0}^{N}q_{2j+1}(0,\zeta)\zeta^{N-j}.
 \label{eq:T0axisqsum}
\end{equation}
The same scalar occurs in the lower-right entry of \(T_\infty\), while the
upper-left entries are its complex conjugate.

It remains to simplify the numerator in \eqref{eq:T0axisqsum}.  Let
\begin{equation}
 W(\zeta)=(\zeta^N,0,\zeta^{N-1},0,\ldots,\zeta,0,1)^{\mathsf T}.
 \label{eq:Wzetadefinition}
\end{equation}
Expansion along the first column gives
\begin{equation}
 \sum_{j=0}^{N}q_{2j+1}(0,\zeta)\zeta^{N-j}
 =[W(\zeta),C_1^I,\ldots,C_{2N}^I].
 \label{eq:qsumdeterminant}
\end{equation}
Subtract \(W(\zeta)\) from every fixed column.  The \(i\)-th resulting
column has the common factor \(K_i-\zeta\).  Extracting all these factors
gives \(\det X\).  Expansion along the final row, followed by unit-triangular
row operations based on
\begin{equation}
 \frac{K_i^{m+1}-\zeta^{m+1}}{K_i-\zeta}
 =K_i^m+\zeta K_i^{m-1}+\cdots+\zeta^m,
 \label{eq:divideddifferencepowers}
\end{equation}
reduces the remaining determinant to
\begin{equation}
 [dK^{N-1},c^IK^{N-1},\ldots,dK,c^IK,d,c^I].
 \label{eq:qsumrows}
\end{equation}
Interchanging the two rows in each of the \(N\) pairs yields
\begin{equation}
 \sum_{j=0}^{N}q_{2j+1}(0,\zeta)\zeta^{N-j}
 =(-1)^N\det(X)H_N^I.
 \label{eq:qsumHidentity}
\end{equation}
Hence
\begin{equation}
 T_0(0,\zeta)=T_\infty(0,\zeta)
 =\begin{pmatrix}\overline{\Theta_N^I(\zeta)}&0\\
                  0&\Theta_N^I(\zeta)\end{pmatrix},
 \qquad
 \Theta_N^I(\zeta)
 =(-1)^N\frac{\det(X)H_N^I}
 {\zeta^N\Delta_N(0,\zeta)}.
 \label{eq:T0axisHcriterion}
\end{equation}
At an interior point of \(I\), \(\det X\ne0\).  Moreover, the explicit
expressions \eqref{eq:delta0scaled}--\eqref{eq:delta2Nscaled} show that
\(T_0\) and \(T_\infty\) extend smoothly to any such axis point at which
\(\zeta\ne0\) and \(\Delta_N\ne0\).  If \(H_N^I\ne0\), their boundary
values are invertible; continuity then makes both matrices invertible in a
sufficiently small interior half-neighbourhood.  This proves the final
assertion.
\end{proof}

\section{Boundary values}

\subsection{The Ernst potential on the open boundary intervals}

We begin by calculating the restriction of the determinant Ernst potential
to every open boundary interval.  This calculation uses only the determinant
formula of Section~3 and does not assume continuity at the endpoints \(K_i\).

For notational convenience, put
\[
 K_0:=+\infty,
 \qquad
 K_{2N+1}:=-\infty,
\]
and define
\begin{equation}
 I_m:=\{(0,\zeta):K_{m+1}<\zeta<K_m\},
 \qquad m\in\{0,\ldots,2N\}.
 \label{eq:Imdefinition}
\end{equation}
In the notation of Section~3.1,
\begin{equation}
 I_{2j-1}=\mathcal H^j,
 \qquad j\in\{1,\ldots,N\},
 \label{eq:Imhorizons}
\end{equation}
and
\begin{equation}
 I_{2j}=\mathcal A^j,
 \qquad j\in\{0,\ldots,N\}.
 \label{eq:Imaxes}
\end{equation}
In particular,
\begin{equation}
 I_0=\mathcal A^+,
 \qquad
 I_{2N}=\mathcal A^-.
 \label{eq:Imouteraxes}
\end{equation}

Fix \(m\in\{0,\ldots,2N\}\).  For \(\zeta\in I_m\), define
\begin{equation}
 \lambda_i^{I_m}:=
 \begin{cases}
  +1,&i\leq m,\\
  -1,&i>m,
 \end{cases}
 \qquad i\in\{1,\ldots,2N\}.
 \label{eq:lambdaIm}
\end{equation}
Indeed, the ordering \(K_1>\cdots>K_{2N}\) implies
\begin{equation}
 r_i(0,\zeta)=|K_i-\zeta|
 =\lambda_i^{I_m}(K_i-\zeta).
 \label{eq:riIm}
\end{equation}
Set
\begin{equation}
 u_i^{I_m}:=\alpha_i\lambda_i^{I_m},
 \qquad i\in\{1,\ldots,2N\},
 \label{eq:uiIm}
\end{equation}
and introduce
\begin{equation}
 c^{I_m}:=(u_1^{I_m},\ldots,u_{2N}^{I_m}),
 \qquad d:=(1,\ldots,1).
 \label{eq:cdIm}
\end{equation}
Let
\begin{equation}
 X=X(\zeta):=
 \operatorname{diag}(K_1-\zeta,\ldots,K_{2N}-\zeta)
 \label{eq:Xboundary}
\end{equation}
and define
\begin{equation}
 c_j^{I_m}:=c^{I_m}X^j,
 \qquad d_j:=dX^j,
 \qquad j\in\{0,\ldots,N\}.
 \label{eq:cjdjboundary}
\end{equation}
Thus \(c_0^{I_m}=c^{I_m}\) and \(d_0=d\).  In the proof below we fix one
interval \(I=I_m\) and suppress the subscript \(m\).

\begin{proposition}[Boundary value of the Ernst potential]
\label{prop:boundaryErnst}
Let \(I\) be one of the open boundary intervals.  Define
\begin{align}
 \widetilde P_N^I(\zeta)
 &:=[c_N^I-d_N,d_{N-1},c_{N-1}^I,
 d_{N-2},c_{N-2}^I,\ldots,d_1,c_1^I,d],
 \label{eq:PtildeI}\\
 \widetilde Q_N^I(\zeta)
 &:=[c_N^I+d_N,d_{N-1},c_{N-1}^I,
 d_{N-2},c_{N-2}^I,\ldots,d_1,c_1^I,d].
 \label{eq:QtildeI}
\end{align}
Then \(\widetilde P_N^I,\widetilde Q_N^I\in\mathbb C[\zeta]\) have degree
at most \(N\).  Their coefficients of \(\zeta^N\) coincide and satisfy
\begin{equation}
 [\zeta^N]\widetilde P_N^I
 =[\zeta^N]\widetilde Q_N^I
 =(-1)^NH_N^I,
 \label{eq:commonleadingcoefficient}
\end{equation}
where
\begin{equation}
 H_N^I=[c^IK^{N-1},dK^{N-1},\ldots,c^IK,dK,c^I,d].
 \label{eq:HNboundaryagain}
\end{equation}
At every \(\zeta\in I\) satisfying
\begin{equation}
 \Delta_N^I(\zeta):=\Delta_N(0,\zeta)\ne0,
 \label{eq:DeltaIboundarynonzero}
\end{equation}
the determinant Ernst potential has the boundary value
\begin{equation}
 f^I(\zeta):=f(0,\zeta)|_I
 =\frac{\widetilde P_N^I(\zeta)}
 {\widetilde Q_N^I(\zeta)}.
 \label{eq:fIboundary}
\end{equation}
Moreover,
\begin{equation}
 \widetilde Q_N^I(\zeta)=\Delta_N^I(\zeta).
 \label{eq:QtildeDelta}
\end{equation}
If \(H_N^I\ne0\), then both polynomials have degree exactly \(N\).  In
that case define
\begin{equation}
 P_N^I:=\frac{\widetilde P_N^I}{(-1)^NH_N^I},
 \qquad
 Q_N^I:=\frac{\widetilde Q_N^I}{(-1)^NH_N^I}.
 \label{eq:monicPIQI}
\end{equation}
Then \(P_N^I\) and \(Q_N^I\) are monic polynomials of degree \(N\), and
\begin{equation}
 f^I=\frac{P_N^I}{Q_N^I}.
 \label{eq:fImonic}
\end{equation}
\end{proposition}

\begin{proof}
Let \(q_\ell^I\), \(\ell\in\{1,\ldots,2N+1\}\), denote the boundary
values on \(I\) of the signed cofactors introduced in Section~3.4.  The
first columns of \(\mathcal M_N^{(+)}\) and \(\mathcal M_N^{(-)}\) are
respectively
\[
 e_+=(1,1,0,\ldots,0)^{\mathsf T},
 \qquad
 e_-=(1,-1,0,\ldots,0)^{\mathsf T}.
\]
Expansion along the first column gives
\begin{equation}
 \det\mathcal M_N^{(+)}(0,\zeta)=q_1^I+q_2^I
 \label{eq:Mplusq12}
\end{equation}
and
\begin{equation}
 \det\mathcal M_N^{(-)}(0,\zeta)
 =q_1^I-q_2^I=\Delta_N^I(\zeta).
 \label{eq:Mminusq12}
\end{equation}
The notation \(\widetilde Q_N^I\) therefore does not define a new
function: it is the same determinant \(\Delta_N^I\), now regarded
explicitly as a polynomial in the running boundary coordinate \(\zeta\).
We retain both notations because they emphasize two different roles of the
same object.

On \(I\), the fixed columns are
\begin{equation}
 C_i^I=
 \begin{pmatrix}
  K_i^N\\
  u_i^I(K_i-\zeta)K_i^{N-1}\\
  K_i^{N-1}\\
  u_i^I(K_i-\zeta)K_i^{N-2}\\
  K_i^{N-2}\\
  \vdots\\
  K_i\\
  u_i^I(K_i-\zeta)\\
  1
 \end{pmatrix},
 \qquad i\in\{1,\ldots,2N\}.
 \label{eq:CiIboundary}
\end{equation}
Define
\begin{align}
 A_N^I
 &:=[c_1^IK^{N-1},dK^{N-1},c_1^IK^{N-2},
 dK^{N-2},\ldots,c_1^I,d],
 \label{eq:ANIboundary}\\
 B_N^I
 &:=[dK^N,dK^{N-1},c_1^IK^{N-2},
 dK^{N-2},\ldots,c_1^I,d].
 \label{eq:BNIboundary}
\end{align}
The sign in the second cofactor is important.  Direct expansion gives
\begin{equation}
 q_1^I=A_N^I,
 \qquad
 q_2^I=-B_N^I.
 \label{eq:q1q2AB}
\end{equation}
Consequently,
\begin{equation}
 q_1^I+q_2^I=A_N^I-B_N^I,
 \qquad
 q_1^I-q_2^I=A_N^I+B_N^I.
 \label{eq:qplusminusAB}
\end{equation}

Using \(K=X+\zeta I_{2N}\), unit-triangular row operations give
\begin{align}
 A_N^I
 &=[c_N^I,d_{N-1},c_{N-1}^I,\ldots,d_1,c_1^I,d],
 \label{eq:ANIafterX}\\
 B_N^I
 &=[d_N,d_{N-1},c_{N-1}^I,\ldots,d_1,c_1^I,d].
 \label{eq:BNIafterX}
\end{align}
These determinants have the same final \(2N-1\) rows.  Linearity in the
first row therefore gives
\begin{equation}
 q_1^I+q_2^I=\widetilde P_N^I,
 \qquad
 q_1^I-q_2^I=\widetilde Q_N^I.
 \label{eq:q12PQtilde}
\end{equation}
Equations \eqref{eq:Mplusq12}--\eqref{eq:Mminusq12} prove
\eqref{eq:fIboundary} and \eqref{eq:QtildeDelta}.

It remains to determine the degrees and leading coefficients.  In the
original fixed columns \eqref{eq:CiIboundary}, dependence on \(\zeta\)
occurs only in the \(N\) rows of the form
\begin{equation}
 u_i^I(K_i-\zeta)K_i^j
 =u_i^IK_i^{j+1}-\zeta u_i^IK_i^j.
 \label{eq:zetarowsboundary}
\end{equation}
Every determinant term therefore contains at most \(N\) factors of
\(\zeta\), and hence
\begin{equation}
 \deg\widetilde P_N^I\leq N,
 \qquad
 \deg\widetilde Q_N^I\leq N.
 \label{eq:degreePQtilde}
\end{equation}

To obtain the coefficient of \(\zeta^N\), one must choose the term
\(-\zeta u_i^IK_i^j\) from every one of these \(N\) rows.  In the Laplace
expansion along \(e_\pm\), the contribution from the second entry of the
first column deletes one of these rows and therefore has degree at most
\(N-1\).  Thus the degree-\(N\) term comes only from the first entry of
\(e_\pm\), and
\begin{equation}
 [\zeta^N]\widetilde P_N^I
 =[\zeta^N]\widetilde Q_N^I
 =(-1)^N[c^IK^{N-1},dK^{N-1},\ldots,c^I,d].
 \label{eq:leadingPQtildeproof}
\end{equation}
The determinant on the right is \(H_N^I\), proving
\eqref{eq:commonleadingcoefficient}.  The final assertion follows by
dividing both polynomials by their common leading coefficient.
\end{proof}

The monic denominator polynomial on the upper axis also recovers the
parameters \(\alpha_i\) directly.

\begin{corollary}[Recovery of the parameters \(\alpha_i\)]
\label{cor:recoveralpha}
Assume that \(H_N^+\ne0\), and let \(Q_N^+\) be the monic polynomial
defined in Proposition~\ref{prop:boundaryErnst}.  For every
\(i\in\{1,\ldots,2N\}\) such that \(Q_N^+(K_i)\ne0\),
\begin{equation}
 \alpha_i=
 \frac{\overline{Q_N^+(K_i)}}{Q_N^+(K_i)}.
 \label{eq:recoveralpha}
\end{equation}
\end{corollary}

\begin{proof}
On the upper axis \(\mathcal A^+\), one has
\begin{equation}
 \lambda_k^+=-1,
 \qquad
 u_k^+=-\alpha_k,
 \qquad k\in\{1,\ldots,2N\}.
 \label{eq:upperaxisuk}
\end{equation}
For brevity, write \(u_k:=u_k^+\), \(c:=(u_1,\ldots,u_{2N})\), and
retain the notation \(c_j=cX^j\), \(d_j=dX^j\).  Since
\(|\alpha_k|=1\),
\begin{equation}
 u_k\overline{u_k}=1.
 \label{eq:ukumodone}
\end{equation}
By Proposition~\ref{prop:boundaryErnst},
\begin{equation}
 \widetilde Q_N^+(\zeta)
 =[c_N+d_N,d_{N-1},c_{N-1},\ldots,d_1,c_1,d].
 \label{eq:Qupperdet}
\end{equation}

Fix \(i\in\{1,\ldots,2N\}\).  At \(\zeta=K_i\), the \(i\)-th
diagonal entry of \(X\) vanishes.  Hence the \(i\)-th entry of every
row \(c_j\) and \(d_j\), \(j\ge1\), is zero.  The \(i\)-th column of
the determinant in \eqref{eq:Qupperdet} is therefore
\((0,\ldots,0,1)^{\mathsf T}\).  Laplace expansion along this column
gives
\begin{equation}
 \widetilde Q_N^+(K_i)=(-1)^{2N+i}\mathcal Q_i=(-1)^i\mathcal Q_i,
 \label{eq:QiLaplace}
\end{equation}
where
\begin{equation}
 \mathcal Q_i
 :=[c_N+d_N,d_{N-1},c_{N-1},\ldots,d_1,c_1]_{\widehat i}.
 \label{eq:calQi}
\end{equation}
The subscript \(\widehat i\) means that the \(i\)-th column has been
deleted.

Conjugate \eqref{eq:calQi} and multiply every remaining column indexed
by \(k\ne i\) by \(u_k\).  At \(\zeta=K_i\), all entries of \(X\) are
real, and \eqref{eq:ukumodone} yields
\begin{equation}
 \overline{c_j}\,\operatorname{diag}(u_k)_{k\ne i}=d_j,
 \qquad
 \overline{d_j}\,\operatorname{diag}(u_k)_{k\ne i}=c_j.
 \label{eq:cdconjugationsubmatrix}
\end{equation}
Consequently,
\begin{equation}
 \left(\prod_{k\ne i}u_k\right)\overline{\mathcal Q_i}
 =[c_N+d_N,c_{N-1},d_{N-1},\ldots,c_1,d_1]_{\widehat i}.
 \label{eq:QiConjugated}
\end{equation}
Interchanging the two rows in each of the \(N-1\) pairs
\((c_j,d_j)\), \(j=1,\ldots,N-1\), restores the row order in
\eqref{eq:calQi}.  Thus
\begin{equation}
 \left(\prod_{k\ne i}u_k\right)\overline{\mathcal Q_i}
 =(-1)^{N-1}\mathcal Q_i.
 \label{eq:QiConjugationIdentity}
\end{equation}
The Laplace sign in \eqref{eq:QiLaplace} is real, and hence the same
identity holds with \(\mathcal Q_i\) replaced by
\(\widetilde Q_N^+(K_i)\):
\begin{equation}
 \left(\prod_{k\ne i}u_k\right)
 \overline{\widetilde Q_N^+(K_i)}
 =(-1)^{N-1}\widetilde Q_N^+(K_i).
 \label{eq:QtildeConjugationAtKi}
\end{equation}

Put \(A:=\prod_{k=1}^{2N}\alpha_k\).  Since \(u_k=-\alpha_k\),
\begin{equation}
 \prod_{k\ne i}u_k=-\frac{A}{\alpha_i}.
 \label{eq:productuexcepti}
\end{equation}
Provided that \(\widetilde Q_N^+(K_i)\ne0\), equations
\eqref{eq:QtildeConjugationAtKi} and \eqref{eq:productuexcepti} imply
\begin{equation}
 \frac{\overline{\widetilde Q_N^+(K_i)}}
      {\widetilde Q_N^+(K_i)}
 =(-1)^N\frac{\alpha_i}{A}.
 \label{eq:QtildeRatioAtKi}
\end{equation}

It remains to account for the monic normalization.  By
Proposition~\ref{prop:boundaryErnst}, the leading coefficient of
\(\widetilde Q_N^+\) is
\begin{equation}
 L^+:=[\zeta^N]\widetilde Q_N^+=(-1)^NH_N^+,
 \label{eq:LplusHplus}
\end{equation}
where
\begin{equation}
 H_N^+=[cK^{N-1},dK^{N-1},\ldots,cK,dK,c,d].
 \label{eq:Hplusdetrecoveralpha}
\end{equation}
Conjugating this determinant and multiplying its \(k\)-th column by
\(u_k\), for every \(k\), interchanges the two rows in each of its
\(N\) row pairs.  Since \(\prod_{k=1}^{2N}u_k=A\), one obtains
\begin{equation}
 A\,\overline{H_N^+}=(-1)^NH_N^+.
 \label{eq:HplusConjugation}
\end{equation}
The assumption \(H_N^+\ne0\) therefore gives
\begin{equation}
 \frac{L^+}{\overline{L^+}}
 =\frac{H_N^+}{\overline{H_N^+}}
 =(-1)^NA.
 \label{eq:LplusPhase}
\end{equation}
Finally, \(Q_N^+=\widetilde Q_N^+/L^+\), and hence
\begin{align}
 \frac{\overline{Q_N^+(K_i)}}{Q_N^+(K_i)}
 &=\frac{L^+}{\overline{L^+}}
   \frac{\overline{\widetilde Q_N^+(K_i)}}
        {\widetilde Q_N^+(K_i)} \notag\\
 &=\left((-1)^NA\right)
   \left((-1)^N\frac{\alpha_i}{A}\right)
 =\alpha_i.
 \label{eq:recoveralphaend}
\end{align}
This proves \eqref{eq:recoveralpha}.
\end{proof}

\medskip
\noindent\textit{Endpoint determinants and continuity.}
We next compare the denominator polynomials belonging to two adjacent
boundary intervals.  Fix \(k\in\{1,\ldots,2N\}\).  On the interval
\(I_{k-1}\), immediately above \(K_k\), the parameters are
\begin{equation}
 u_i^{I_{k-1}}=
 \begin{cases}
  \alpha_i,&i<k,\\
  -\alpha_i,&i\ge k.
 \end{cases}
 \label{eq:ukbeforeendpoint}
\end{equation}
Delete the \(k\)-th entry and define
\begin{equation}
 X_k:=\operatorname{diag}(K_i-K_k)_{i\ne k},
 \label{eq:Xkendpoint}
\end{equation}
\begin{equation}
 c^{(k)}:=(\alpha_1,\ldots,\alpha_{k-1},
           -\alpha_{k+1},\ldots,-\alpha_{2N}),
 \qquad
 d^{(k)}:=(1,\ldots,1),
 \qquad
 v^{(k)}:=c^{(k)}+d^{(k)}.
 \label{eq:cdvkendpoint}
\end{equation}
For \(j\ge1\), put
\begin{equation}
 c_j^{(k)}:=c^{(k)}X_k^j,
 \qquad
 d_j^{(k)}:=d^{(k)}X_k^j,
 \qquad
 v_j^{(k)}:=v^{(k)}X_k^j.
 \label{eq:cdvjkendpoint}
\end{equation}
For \(N\ge2\), define the reduced endpoint determinant
\begin{equation}
 E_k:=[v_{N-1}^{(k)},\ldots,v_1^{(k)},
       c_{N-2}^{(k)},\ldots,c_1^{(k)},c^{(k)},d^{(k)}].
 \label{eq:Ekendpoint}
\end{equation}
When \(N=1\), the corresponding definition is
\begin{equation}
 E_k:=[v^{(k)}].
 \label{eq:EkendpointNone}
\end{equation}
Thus \(E_k\) is a determinant of order \(2N-1\) in both cases.

\begin{proposition}[Endpoint denominator]
\label{prop:endpointdenominator}
Define
\begin{equation}
 \Delta_k:=\widetilde Q_N^{I_{k-1}}(K_k).
 \label{eq:Deltakdefinition}
\end{equation}
Then
\begin{equation}
 \widetilde Q_N^{I_{k-1}}(K_k)
 =\widetilde Q_N^{I_k}(K_k),
 \label{eq:adjacentQagree}
\end{equation}
so \(\Delta_k\) is the common endpoint value of the two adjacent
denominator polynomials.  Moreover,
\begin{equation}
 \Delta_k
 =(-1)^{\,k+1+\frac{(N-2)(N-3)}2}\det(X_k)E_k.
 \label{eq:DeltakdetXEk}
\end{equation}
Equivalently, if
\begin{equation}
 \Lambda_k:=|\det X_k|
 =\prod_{i<k}(K_i-K_k)\prod_{i>k}(K_k-K_i)>0,
 \label{eq:Lambdakendpoint}
\end{equation}
then
\begin{equation}
 \Delta_k
 =(-1)^{\,1+\frac{(N-2)(N-3)}2}\Lambda_kE_k.
 \label{eq:DeltakLambdaEk}
\end{equation}
\end{proposition}

\begin{proof}
When one crosses \(K_k\), only \(u_k^I\) changes sign.  At
\(\zeta=K_k\), however,
\begin{equation}
 (c_j^I)_k=u_k^I(K_k-\zeta)^j=0,
 \qquad
 (d_j)_k=(K_k-\zeta)^j=0,
 \qquad j\ge1.
 \label{eq:endpointcolumnvanishing}
\end{equation}
Every occurrence of \(u_k^I\) in \(\widetilde Q_N^I\) is therefore
annihilated, proving \eqref{eq:adjacentQagree}.

At \(\zeta=K_k\), the \(k\)-th column of the determinant
\eqref{eq:QtildeI} is \((0,\ldots,0,1)^{\mathsf T}\).  Since this
determinant has \(2N\) rows, Laplace expansion along that column gives
\begin{align}
 \Delta_k
 &=(-1)^{2N+k}
 [c_N+d_N,d_{N-1},c_{N-1},\ldots,d_1,c_1]_{\widehat k}\notag\\
 &=(-1)^k
 [c_N+d_N,d_{N-1},c_{N-1},\ldots,d_1,c_1]_{\widehat k}.
 \label{eq:DeltakLaplace}
\end{align}
All row vectors in the determinant on the right are now restricted to
the \(2N-1\) columns different from \(k\).  On those columns the matrix
\(X\) becomes the invertible diagonal matrix \(X_k\).  Right
multiplication by \(X_k^{-1}\) therefore gives
\begin{align}
 &[c_N+d_N,d_{N-1},c_{N-1},\ldots,d_1,c_1]_{\widehat k}
 \notag\\
 &\quad=\det(X_k)
 [v_{N-1}^{(k)},d_{N-2}^{(k)},c_{N-2}^{(k)},\ldots,
   d_1^{(k)},c_1^{(k)},d^{(k)},c^{(k)}].
 \label{eq:endpointfactorXk}
\end{align}

For every \(j\in\{1,\ldots,N-2\}\), apply the row operation
\begin{equation}
 d_j^{(k)}\longmapsto d_j^{(k)}+c_j^{(k)}=v_j^{(k)}.
 \label{eq:endpointdjtoj}
\end{equation}
The resulting ordered row group is
\begin{equation}
 \bigl(v_{N-1}^{(k)},v_{N-2}^{(k)},c_{N-2}^{(k)},
 v_{N-3}^{(k)},c_{N-3}^{(k)},\ldots,
 v_1^{(k)},c_1^{(k)},d^{(k)},c^{(k)}\bigr).
 \label{eq:endpointinterleavedrows}
\end{equation}
Moving all the \(v\)-rows in front of the \(c\)-rows requires
\begin{equation}
 0+1+\cdots+(N-3)=\frac{(N-2)(N-3)}2
 \label{eq:endpointpermutationcount}
\end{equation}
row interchanges.  One final interchange changes the last two rows
from \((d^{(k)},c^{(k)})\) to \((c^{(k)},d^{(k)})\).  Thus
\begin{align}
 &[c_N+d_N,d_{N-1},c_{N-1},\ldots,d_1,c_1]_{\widehat k}
 \notag\\
 &\quad=(-1)^{\,1+\frac{(N-2)(N-3)}2}\det(X_k)E_k.
 \label{eq:endpointminorEk}
\end{align}
Together with \eqref{eq:DeltakLaplace}, this proves
\eqref{eq:DeltakdetXEk} for \(N\ge2\).  For \(N=1\), the same formula
follows directly by expanding the defining \(2\times2\) determinant.

Finally, the ordering \(K_1>\cdots>K_{2N}\) implies
\begin{equation}
 \det(X_k)=(-1)^k\Lambda_k.
 \label{eq:detXkLambdakendpoint}
\end{equation}
Substitution in \eqref{eq:DeltakdetXEk} proves
\eqref{eq:DeltakLambdaEk}.
\end{proof}

\begin{corollary}[Continuity at an endpoint]
\label{cor:endpointcontinuity}
For every \(k\in\{1,\ldots,2N\}\),
\begin{equation}
 \Delta_k\ne0\quad\Longleftrightarrow\quad E_k\ne0.
 \label{eq:DeltakEkNonzero}
\end{equation}
If this condition holds, then the boundary Ernst potential extends
continuously across \(K_k\).
\end{corollary}

\begin{proof}
The equivalence \eqref{eq:DeltakEkNonzero} follows immediately from
\eqref{eq:DeltakLambdaEk} and \(\Lambda_k>0\).  The same argument used
in \eqref{eq:endpointcolumnvanishing} for the numerator gives
\begin{equation}
 G_k:=\widetilde P_N^{I_{k-1}}(K_k)
     =\widetilde P_N^{I_k}(K_k).
 \label{eq:Gkdefinition}
\end{equation}
If \(E_k\ne0\), then \(\Delta_k\ne0\), and hence
\begin{equation}
 \lim_{\substack{\zeta\to K_k\\\zeta\in I_{k-1}}}f(0,\zeta)
 =\frac{G_k}{\Delta_k}
 =\lim_{\substack{\zeta\to K_k\\\zeta\in I_k}}f(0,\zeta).
 \label{eq:fendpointlimits}
\end{equation}
Defining
\begin{equation}
 f(0,K_k):=\frac{G_k}{\Delta_k}
 \label{eq:fendpointdefinition}
\end{equation}
therefore makes the boundary restriction continuous at \(K_k\).
Moreover, the original numerator and denominator determinants are
continuous functions of \((\rho,\zeta)\).  Their denominator is
nonzero at \((0,K_k)\), and hence remains nonzero in a sufficiently
small closed half-neighbourhood.  Thus \(f\) itself is continuous in
such a half-neighbourhood.
\end{proof}

\subsection{Boundary values of \texorpdfstring{\(e^{2U}\)}{exp(2U)}}
\label{subsec:boundarye2U}

The calculation of the real part of the boundary Ernst potential will use
the following elementary determinant identity.  We state it
separately so that it can also be cited in the later calculations of the
remaining metric functions.

\begin{lemma}[Pl\"ucker relation]
\label{lem:plucker}
Let \(V\) be an \(n\)-dimensional complex vector space, let
\(\mathcal R=(r_1,\ldots,r_{n-2})\) be an ordered family in \(V\), and fix
a volume form on \(V\).  For \(\xi,\eta\in V\), write
\([mathcal R,\xi,\eta]\) for the value of this volume form on the ordered
\(n\)-tuple \((r_1,\ldots,r_{n-2},\xi,\eta)\).  Then, for arbitrary
\(\xi,\eta,\mu,\nu\in V\),
\begin{align}
 &[\mathcal R,\xi,\eta][\mathcal R,\mu,\nu]
 -[\mathcal R,\xi,\nu][\mathcal R,\mu,\eta]
 \notag\\
 &\qquad=[\mathcal R,\xi,\mu][\mathcal R,\eta,\nu].
 \label{eq:pluckerrelation}
\end{align}
\end{lemma}

\begin{proof}
If the vectors in \(\mathcal R\) are linearly dependent, every determinant
in \eqref{eq:pluckerrelation} vanishes.  We may therefore assume that
\begin{equation}
 W:=\operatorname{span}\{r_1,\ldots,r_{n-2}\}
 \label{eq:pluckerW}
\end{equation}
has dimension \(n-2\).  On the two-dimensional quotient space
\(Q:=V/W\), define
\begin{equation}
 \omega([x],[y]):=[\mathcal R,x,y],
 \label{eq:quotientomega}
\end{equation}
where \([x]\) denotes the class of \(x\) modulo \(W\).  This is a
well-defined alternating bilinear form: adding an element of \(W\) to
either argument produces a determinant whose rows are linearly dependent.

Assume first that
\begin{equation}
 D:=[\mathcal R,\xi,\mu]\ne0.
 \label{eq:pluckerD}
\end{equation}
Then \([\xi],[\mu]\) form a basis of \(Q\), so there exist unique scalars
\(a,b,c,d\in\mathbb C\) such that
\begin{equation}
 [\eta]=a[\xi]+b[\mu],
 \qquad
 [\nu]=c[\xi]+d[\mu].
 \label{eq:quotientexpansions}
\end{equation}
Alternation and bilinearity give
\begin{equation}
 [\mathcal R,\xi,\eta]=bD,
 \qquad
 [\mathcal R,\mu,\nu]=-cD,
 \label{eq:pluckerfourone}
\end{equation}
and
\begin{equation}
 [\mathcal R,\xi,\nu]=dD,
 \qquad
 [\mathcal R,\mu,\eta]=-aD.
 \label{eq:pluckerfourtwo}
\end{equation}
Consequently, the left-hand side of \eqref{eq:pluckerrelation} equals
\begin{equation}
 (bD)(-cD)-(dD)(-aD)=(ad-bc)D^2.
 \label{eq:pluckerLHS}
\end{equation}
On the other hand,
\begin{align}
 [\mathcal R,\eta,\nu]
 &=\omega(a[\xi]+b[\mu],c[\xi]+d[\mu])\notag\\
 &=(ad-bc)D.
 \label{eq:pluckerlastbracket}
\end{align}
Multiplication by \([\mathcal R,\xi,\mu]=D\) proves
\eqref{eq:pluckerrelation} whenever \(D\ne0\).

Finally, both sides of \eqref{eq:pluckerrelation} are polynomial
expressions in the components of all the displayed vectors.  Since they
agree on the open set \(D\ne0\), they agree identically.  Thus the relation
also holds when \(D=0\).
\end{proof}

\begin{proposition}[Boundary values of \(e^{2U}\)]
\label{prop:boundarye2U}
Fix an open boundary interval
\begin{equation}
 I=I_m,
 \qquad m\in\{0,\ldots,2N\}.
 \label{eq:fixImU}
\end{equation}
At every \(\zeta\in I\) for which \(\Delta_N^I(\zeta)\ne0\),
\begin{equation}
 e^{2U^I}(\zeta)
 =\frac{\det X(\zeta)\,|H_N^I|^2}
        {|\Delta_N^I(\zeta)|^2},
 \label{eq:boundarye2U}
\end{equation}
where \(X\) is defined in \eqref{eq:Xboundary}.  If, in addition,
\(H_N^I\ne0\), then, wherever \(\Delta_N^I\ne0\),
\begin{equation}
 \operatorname{sgn}e^{2U^I}
 =\operatorname{sgn}\det X=(-1)^m.
 \label{eq:signe2UI}
\end{equation}
Consequently,
\begin{equation}
 e^{2U^{\mathcal A^j}}>0,
 \qquad j\in\{0,\ldots,N\},
 \label{eq:e2Uaxispositive}
\end{equation}
on the ordinary axis intervals, whereas
\begin{equation}
 e^{2U^{\mathcal H^j}}<0,
 \qquad j\in\{1,\ldots,N\},
 \label{eq:e2Uhorizonnegative}
\end{equation}
on the candidate horizon intervals, provided that the corresponding
\(H_N^I\) and \(\Delta_N^I\) do not vanish.  Moreover, if \(E_k\ne0\),
then \(f\) extends continuously to \(K_k\) by
Corollary~\ref{cor:endpointcontinuity}, and
\begin{equation}
 e^{2U}(0,K_k)=0.
 \label{eq:e2Uendpointzero}
\end{equation}
\end{proposition}

\begin{proof}
We retain the boundary notation \(u_i^I,c^I,d,K,H_N^I\) introduced in
\eqref{eq:axisuibeforeprop}--\eqref{eq:HNaxisdeterminant}; the associated
matrix \(X\) and rows \(c_j^I,d_j\) are those of
\eqref{eq:Xboundary}--\eqref{eq:cjdjboundary}.  We also use the signed
cofactors \(q_\ell^I\), \(\ell\in\{1,\ldots,2N+1\}\), obtained as the
boundary values on \(I\) of the cofactors \(q_\ell\) defined in
\eqref{eq:qcofactorsexpansion}.

By Proposition~\ref{prop:boundaryErnst},
\begin{equation}
 f^I=\frac{q_1^I+q_2^I}{q_1^I-q_2^I},
 \qquad
 \Delta_N^I=q_1^I-q_2^I.
 \label{eq:fIq12fore2U}
\end{equation}
It follows that
\begin{align}
 e^{2U^I}=\Re f^I
 &=\frac12\left(
   \frac{q_1^I+q_2^I}{q_1^I-q_2^I}
  +\frac{\overline{q_1^I}+\overline{q_2^I}}
         {\overline{q_1^I}-\overline{q_2^I}}
  \right)\notag\\
 &=\frac{q_1^I\overline{q_1^I}-q_2^I\overline{q_2^I}}
         {|\Delta_N^I|^2}.
 \label{eq:realfq12}
\end{align}

Put
\begin{equation}
 U_0^I:=\prod_{i=1}^{2N}u_i^I
 \label{eq:U0I}
\end{equation}
and introduce the ordered family of \(2N-2\) common rows
\begin{equation}
 \mathcal R^I:=(c_{N-1}^I,d_{N-1},\ldots,c_1^I,d_1).
 \label{eq:RIe2U}
\end{equation}
Reordering the cofactor determinants gives
\begin{equation}
 q_1^I=(-1)^{N-1}[\mathcal R^I,c_N^I,d],
 \qquad
 q_2^I=(-1)^N[\mathcal R^I,d_N,d].
 \label{eq:q12RI}
\end{equation}
Since \(|u_i^I|=1\), conjugating a determinant and then multiplying its
\(i\)-th column by \(u_i^I\) interchanges \(c_j^I\) and \(d_j\).
Consequently,
\begin{equation}
 U_0^I\overline{q_1^I}=[\mathcal R^I,d_N,c^I],
 \qquad
 U_0^I\overline{q_2^I}=-[\mathcal R^I,c_N^I,c^I].
 \label{eq:conjugateq12RI}
\end{equation}
Thus
\begin{align}
 &U_0^I\left(q_1^I\overline{q_1^I}
                    -q_2^I\overline{q_2^I}\right)
 \notag\\
 &\quad=(-1)^{N-1}\Bigl(
 [\mathcal R^I,c_N^I,d][\mathcal R^I,d_N,c^I]
 -[\mathcal R^I,c_N^I,c^I][\mathcal R^I,d_N,d]
 \Bigr).
 \label{eq:e2UbeforePlucker}
\end{align}
Apply Lemma~\ref{lem:plucker} with
\begin{equation}
 \xi=c_N^I,
 \qquad \eta=d,
 \qquad \mu=d_N,
 \qquad \nu=c^I.
 \label{eq:e2UPluckersubstitution}
\end{equation}
The expression in parentheses in \eqref{eq:e2UbeforePlucker} becomes
\begin{equation}
 [\mathcal R^I,c_N^I,d_N][\mathcal R^I,d,c^I].
 \label{eq:e2UafterPlucker}
\end{equation}
By the definition of \(H_N^I\),
\begin{equation}
 [\mathcal R^I,d,c^I]=-H_N^I.
 \label{eq:secondHfactor}
\end{equation}
Multiplying every column of the matrix defining \(H_N^I\) by the
corresponding diagonal entry of \(X\) gives
\begin{equation}
 [\mathcal R^I,c_N^I,d_N]=\det(X)H_N^I.
 \label{eq:firstHfactor}
\end{equation}
Equations \eqref{eq:e2UbeforePlucker}--\eqref{eq:firstHfactor} yield
\begin{equation}
 U_0^I\left(q_1^I\overline{q_1^I}
                   -q_2^I\overline{q_2^I}\right)
 =(-1)^N\det(X)(H_N^I)^2.
 \label{eq:e2Udetidentity}
\end{equation}

It remains to remove the phase of \(H_N^I\).  Conjugating its defining
determinant and multiplying the \(i\)-th column by \(u_i^I\) interchanges
the two rows in each of its \(N\) row pairs.  Hence
\begin{equation}
 U_0^I\overline{H_N^I}=(-1)^NH_N^I,
 \label{eq:HIconjugation}
\end{equation}
and therefore
\begin{equation}
 \frac{(H_N^I)^2}{U_0^I}=(-1)^N|H_N^I|^2.
 \label{eq:HIphaseRemoval}
\end{equation}
Dividing \eqref{eq:e2Udetidentity} by
\(U_0^I|\Delta_N^I|^2\) and using \eqref{eq:HIphaseRemoval} proves
\eqref{eq:boundarye2U}.

Since \(I=I_m=(K_{m+1},K_m)\), the factors \(K_i-\zeta\) are positive
for \(i\le m\) and negative for \(i>m\).  There are \(2N-m\) negative
factors, and hence
\begin{equation}
 \operatorname{sgn}\det X=(-1)^{2N-m}=(-1)^m.
 \label{eq:signDetXI}
\end{equation}
Together with \(\mathcal A^j=I_{2j}\) and
\(\mathcal H^j=I_{2j-1}\), this proves
\eqref{eq:signe2UI}--\eqref{eq:e2Uhorizonnegative}.

Finally, suppose \(E_k\ne0\).  By
Corollary~\ref{cor:endpointcontinuity}, \(f\), and hence
\(e^{2U}=\Re f\), extends continuously to \(K_k\).  Since
\(\det X(K_k)=0\) and the endpoint denominator is nonzero,
\eqref{eq:boundarye2U} gives \eqref{eq:e2Uendpointzero}.
\end{proof}

\subsection{Boundary values of \texorpdfstring{\(e^{2k}\)}{exp(2k)}}
\label{subsec:boundarye2k}

We begin with a global formula for \(e^{2k}\) obtained by Kramer from the
general B\"acklund transformation.  We shall not reproduce Kramer's general
construction here.  We only record its specialization to the determinant
family introduced in Section~3.

For \(2N\) complex numbers \(\beta_1,\ldots,\beta_{2N}\), define the
Kramer determinant
\begin{equation}
\begin{split}
 &\mathfrak K_N(\beta_1,\ldots,\beta_{2N})\\
 &\quad:=\det
 \begin{pmatrix}
  1&\beta_1\lambda_1&\lambda_1^2&\beta_1\lambda_1^3&
    \cdots&\lambda_1^{2N-2}&\beta_1\lambda_1^{2N-1}\\
  1&\beta_2\lambda_2&\lambda_2^2&\beta_2\lambda_2^3&
    \cdots&\lambda_2^{2N-2}&\beta_2\lambda_2^{2N-1}\\
  \vdots&\vdots&\vdots&\vdots&&\vdots&\vdots\\
  1&\beta_{2N}\lambda_{2N}&\lambda_{2N}^2&
    \beta_{2N}\lambda_{2N}^3&\cdots&\lambda_{2N}^{2N-2}&
    \beta_{2N}\lambda_{2N}^{2N-1}
 \end{pmatrix}.
\end{split}
\label{eq:KramerDeterminant}
\end{equation}
Thus \(\mathfrak K_N\) is a determinant of order \(2N\).  The functions
\(\lambda_i\) on the fixed-\(K_i\) coverings satisfy
\begin{equation}
 \lambda_i^2=
 \frac{K_i-\zeta-i\rho}{K_i-\zeta+i\rho},
 \qquad i\in\{1,\ldots,2N\}.
 \label{eq:lambdaifore2k}
\end{equation}

Kramer's integrated formula, specialized to the present \(2N\)-fold
determinant family with Minkowski seed and parameters \(\alpha_i\), gives
\begin{align}
 e^{2k+C}
 ={}&
 \frac{\displaystyle\prod_{i=1}^{2N}\alpha_i}
      {\displaystyle\prod_{i=1}^{2N}\lambda_i}
 \rho^{2N(N-1)}
 \left[\prod_{i=1}^{2N}(\lambda_i^2-1)\right]^{2-2N}
 \notag\\
 &\times
 \mathfrak K_N(\alpha_1,\ldots,\alpha_{2N})
 \mathfrak K_N(\alpha_1^{-1},\ldots,\alpha_{2N}^{-1}).
 \label{eq:KramerSpecialized}
\end{align}
Formula \eqref{eq:KramerSpecialized} follows from equation~(15) of
Kramer~\cite{Kramer1986}; its normalization and notation have been
translated to those used in the present paper.  It is an identity in the
interior region \(\rho>0\).  Here \(C\in\mathbb R\) is a constant,
independent of \((\rho,\zeta)\), which represents the additive integration
constant in \(k\).  Any later occurrence of the letter \(C\) need not denote
this same constant unless this is stated explicitly.  The individual factors
on its right-hand side
may vanish or diverge as \(\rho\downarrow0\).  In the following proposition
we calculate their combined limit on each open boundary interval \(I\).

\begin{proposition}[Boundary values of \(e^{2k}\)]
\label{prop:boundarye2k}
Retain the boundary notation introduced in
\eqref{eq:axisuibeforeprop}--\eqref{eq:HNaxisdeterminant} and in
Section~4.1, namely
\begin{equation}
 \lambda_i^I,\quad u_i^I,\quad c^I,\quad d,\quad K,\quad X,
 \quad c_j^I,\quad d_j,\quad H_N^I.
 \label{eq:boundaryNotationForK}
\end{equation}
Define
\begin{equation}
 e^{2k^I}:=\lim_{\rho\downarrow0}e^{2k(\rho,\zeta)},
 \qquad \zeta\in I,
 \label{eq:e2kIboundarydefinition}
\end{equation}
whenever this limit exists.  Then, for any two open boundary intervals
\(I\) and \(J\) such that \(H_N^J\ne0\),
\begin{equation}
 \frac{e^{2k^I}}{e^{2k^J}}
 =
 \left(\frac{H_N^I}{H_N^J}\right)^2.
 \label{eq:e2kboundaryratio}
\end{equation}
\end{proposition}

\begin{proof}
Fix an open boundary interval \(I\), and suppress its superscript
temporarily whenever no ambiguity can arise.  Put
\begin{equation}
 x_i:=K_i-\zeta.
\end{equation}
Since \(I\) contains no endpoint, \(x_i\ne0\) for every
\(i\in\{1,\ldots,2N\}\).  By Remark~3.2, on the chosen component of the
real-\(K\) part of the spectral covering,
\begin{equation}
 \lambda_i(\rho,\zeta)
 =\frac{x_i-i\rho}{r_i(\rho,\zeta)}
 \longrightarrow \lambda_i^I
 \qquad(\rho\downarrow0).
 \label{eq:lambdaithetaapproach}
\end{equation}
For sufficiently small \(\rho\), write
\begin{equation}
 \lambda_i=\lambda_i^I e^{-i\theta_i},
 \qquad
 \theta_i=\arctan\left(\frac{\rho}{x_i}\right)
 =\frac{\rho}{x_i}+O(\rho^3),
 \label{eq:thetaiexpansion}
\end{equation}
where the branch of the arctangent is chosen near zero.

Transpose the first Kramer determinant in
\eqref{eq:KramerSpecialized}.  Its rows occur in the order
\begin{equation}
 E_0,O_0,E_1,O_1,\ldots,E_{N-1},O_{N-1},
\end{equation}
where, for \(\ell\in\{0,\ldots,N-1\}\),
\begin{equation}
 E_\ell=(\lambda_1^{2\ell},\ldots,\lambda_{2N}^{2\ell}),
 \qquad
 O_\ell=(\alpha_1\lambda_1^{2\ell+1},\ldots,
          \alpha_{2N}\lambda_{2N}^{2\ell+1}).
 \label{eq:KramerEvenOddRows}
\end{equation}
For \(m\in\{0,\ldots,N-1\}\), define the finite differences
\begin{equation}
 \widetilde E_m
 :=\sum_{\ell=0}^{m}(-1)^{m-\ell}\binom m\ell E_\ell,
 \qquad
 \widetilde O_m
 :=\sum_{\ell=0}^{m}(-1)^{m-\ell}\binom m\ell O_\ell.
 \label{eq:KramerFiniteDifferences}
\end{equation}
These replacements are effected by lower unitriangular row operations
and hence preserve the determinant.

For fixed \(q\), the coefficient of \(\rho^q\) in \(E_\ell\) is a
polynomial in the discrete variable \(\ell\) of degree at most \(q\).
Its degree-\(q\) term is
\begin{equation}
 \frac{(-2i)^q}{q!}\,\ell^q dX^{-q}.
 \label{eq:EvenRowLeadingPolynomial}
\end{equation}
The \(m\)-th finite difference annihilates every polynomial of degree
smaller than \(m\), whereas
\begin{equation}
 \sum_{\ell=0}^{m}(-1)^{m-\ell}\binom m\ell\ell^m=m!.
 \label{eq:FiniteDifferenceMonomial}
\end{equation}
It follows that
\begin{equation}
 \widetilde E_m
 =(-2i)^m\rho^m dX^{-m}+O(\rho^{m+1}).
 \label{eq:EvenFiniteDifferenceAsymptotic}
\end{equation}
Similarly, since
\(
\alpha_i\lambda_i=u_i^I e^{-i\theta_i},
\)
the odd rows satisfy
\begin{equation}
 \widetilde O_m
 =(-2i)^m\rho^m c^IX^{-m}+O(\rho^{m+1}).
 \label{eq:OddFiniteDifferenceAsymptotic}
\end{equation}

The total extracted power of \(\rho\) is
\begin{equation}
 2\sum_{m=0}^{N-1}m=N(N-1).
\end{equation}
After the finite-difference operations, interchange the two rows in each
of the \(N\) pairs.  This contributes the sign \((-1)^N\), and gives
\begin{align}
 \mathfrak K_N(\alpha_1,\ldots,\alpha_{2N})
 ={}&\gamma_N\rho^{N(N-1)}
 [c^I,d,c^IX^{-1},dX^{-1},\ldots,
   c^IX^{-(N-1)},dX^{-(N-1)}]\notag\\
 &+O\!\left(\rho^{N(N-1)+1}\right),
 \label{eq:FirstKramerAxisAsymptotic}
\end{align}
where
\begin{equation}
 \gamma_N:=(-1)^N(-2i)^{N(N-1)}.
\end{equation}
Right multiplication of the matrix in
\eqref{eq:FirstKramerAxisAsymptotic} by \(X^{N-1}\) gives
\begin{equation}
 [c^IX^{N-1},dX^{N-1},\ldots,c^IX,dX,c^I,d].
\end{equation}
Since \(X=K-\zeta\mathbf1\), unitriangular row operations replace the
powers of \(X\) by the corresponding powers of \(K\).  Hence this
determinant is \(H_N^I\).  Since
\begin{equation}
 \det(X^{N-1})
 =\left[\prod_{i=1}^{2N}(K_i-\zeta)\right]^{N-1},
\end{equation}
we obtain
\begin{align}
 \mathfrak K_N(\alpha_1,\ldots,\alpha_{2N})
 ={}&\gamma_N\rho^{N(N-1)}
 \frac{H_N^I}
 {\left[\prod_{i=1}^{2N}(K_i-\zeta)\right]^{N-1}}
 +O\!\left(\rho^{N(N-1)+1}\right).
 \label{eq:FirstKramerReducedAsymptotic}
\end{align}

Because \(|\alpha_i|=1\) and \(\lambda_i^I\in\{1,-1\}\),
\begin{equation}
 \alpha_i^{-1}\lambda_i^I
 =\overline{\alpha_i\lambda_i^I}
 =\overline{u_i^I}.
\end{equation}
The same calculation applied to the second Kramer determinant therefore
gives
\begin{align}
 \mathfrak K_N(\alpha_1^{-1},\ldots,\alpha_{2N}^{-1})
 ={}&\gamma_N\rho^{N(N-1)}
 \frac{\overline{H_N^I}}
 {\left[\prod_{i=1}^{2N}(K_i-\zeta)\right]^{N-1}}
 +O\!\left(\rho^{N(N-1)+1}\right).
 \label{eq:SecondKramerReducedAsymptotic}
\end{align}

The remaining singular factor is governed by the exact identity
\begin{equation}
 \lambda_i^2-1=\frac{-2i\rho}{K_i-\zeta+i\rho}.
 \label{eq:lambdaSquaredMinusOneExact}
\end{equation}
Consequently,
\begin{align}
 \left[\prod_{i=1}^{2N}(\lambda_i^2-1)\right]^{2-2N}
 ={}&(-2i)^{-4N(N-1)}\rho^{-4N(N-1)}
 \left[\prod_{i=1}^{2N}(K_i-\zeta)\right]^{2N-2}
 \bigl(1+O(\rho)\bigr).
 \label{eq:KramerSingularFactorAsymptotic}
\end{align}

Substitution of
\eqref{eq:FirstKramerReducedAsymptotic},
\eqref{eq:SecondKramerReducedAsymptotic}, and
\eqref{eq:KramerSingularFactorAsymptotic} into
\eqref{eq:KramerSpecialized} cancels all powers of \(\rho\) and all
factors \(K_i-\zeta\).  The factor which remains apart from \(H_N^I\)
is independent of the boundary interval \(I\).  The interval-dependent
part is
\begin{equation}
 \frac{\prod_{i=1}^{2N}\alpha_i}
      {\prod_{i=1}^{2N}\lambda_i^I}
 H_N^I\overline{H_N^I}.
 \label{eq:e2kIntervalDependentPart}
\end{equation}

Conjugate the determinant defining \(H_N^I\), and multiply its
\(i\)-th column by \(u_i^I\).  This interchanges \(c^I\) and \(d\) in
each of the \(N\) row pairs.  Therefore
\begin{equation}
 \left(\prod_{i=1}^{2N}u_i^I\right)\overline{H_N^I}
 =(-1)^NH_N^I.
 \label{eq:HNIConjugationForK}
\end{equation}
Since
\begin{equation}
 \prod_{i=1}^{2N}u_i^I
 =\left(\prod_{i=1}^{2N}\alpha_i\right)
  \left(\prod_{i=1}^{2N}\lambda_i^I\right)
\end{equation}
and
\begin{equation}
 \left(\prod_{i=1}^{2N}\lambda_i^I\right)^2=1,
\end{equation}
we have
\begin{equation}
 \frac{\prod_{i=1}^{2N}\alpha_i}
      {\prod_{i=1}^{2N}\lambda_i^I}
 \overline{H_N^I}
 =(-1)^NH_N^I.
 \label{eq:ConjugateHRemovedForK}
\end{equation}
Thus, apart from a factor independent of \(I\), the boundary value
\(e^{2k^I+C}\) is equal to \((H_N^I)^2\).  Applying the same calculation
to \(J\) and taking the quotient cancels both this common factor and the
integration constant \(C\).  This proves \eqref{eq:e2kboundaryratio}.
\end{proof}

\begin{corollary}[Alternating sign]
\label{cor:e2kAlternatingSign}
Assume that
\begin{equation}
 H_N^+\ne0,
 \qquad
 H_N^{I_m}\ne0
\end{equation}
for some \(m\in\{0,\ldots,2N\}\), where
\(I_0=\mathcal A^+\).  Then
\begin{equation}
 \operatorname{sgn}
 \left[
 \left(\frac{H_N^{I_m}}{H_N^+}\right)^2
 \right]
 =(-1)^m.
 \label{eq:HratioAlternatingSign}
\end{equation}
Consequently,
\begin{equation}
 \operatorname{sgn}
 \left(
 \frac{e^{2k^{I_m}}}{e^{2k^+}}
 \right)
 =(-1)^m.
 \label{eq:e2kAlternatingSign}
\end{equation}
\end{corollary}

\begin{proof}
Put
\begin{equation}
 U_0^I:=\prod_{i=1}^{2N}u_i^I.
\end{equation}
When one passes from \(I_{j-1}\) to \(I_j\), exactly one of the factors
\(u_i^I\) changes sign.  Hence
\begin{equation}
 U_0^{I_m}=(-1)^mU_0^+.
 \label{eq:U0AlternatingSign}
\end{equation}
The conjugation identity \eqref{eq:HNIConjugationForK} is equivalent to
\begin{equation}
 (H_N^I)^2=(-1)^NU_0^I|H_N^I|^2.
 \label{eq:HsquareThroughU0}
\end{equation}
Applying this identity to \(I_m\) and \(\mathcal A^+\), and then using
\eqref{eq:U0AlternatingSign}, gives
\begin{equation}
 \left(\frac{H_N^{I_m}}{H_N^+}\right)^2
 =(-1)^m
 \frac{|H_N^{I_m}|^2}{|H_N^+|^2}.
 \label{eq:HratioSignExplicit}
\end{equation}
The quotient on the right is strictly positive, which proves
\eqref{eq:HratioAlternatingSign}.  Equation
\eqref{eq:e2kAlternatingSign} now follows from Proposition~4.7.
\end{proof}

\begin{corollary}[Agreement of the boundary signs]
\label{cor:boundarySignAgreement}
Assume that \(H_N^+\ne0\), and fix the additive constant in \(k\) by
the upper-axis normalization
\begin{equation}
 e^{2k^+}=1.
 \label{eq:upperAxisKNormalization}
\end{equation}
Let \(I=I_m\), where \(m\in\{0,\ldots,2N\}\), and let
\(\zeta\in I_m\) be a point at which
\begin{equation}
 \Delta_N^{I_m}(\zeta)\ne0.
 \label{eq:DeltaNonzeroForSignAgreement}
\end{equation}
Then
\begin{equation}
 \operatorname{sgn}e^{2U^{I_m}}(\zeta)
 =\operatorname{sgn}e^{2k^{I_m}}.
 \label{eq:boundarySignsAgree}
\end{equation}
More precisely,
\begin{equation}
 \operatorname{sgn}e^{2U^{I_m}}(\zeta)
 =\operatorname{sgn}e^{2k^{I_m}}
 =
 \begin{cases}
  (-1)^m,&H_N^{I_m}\ne0,\\[3pt]
  0,&H_N^{I_m}=0.
 \end{cases}
 \label{eq:boundarySignsExplicit}
\end{equation}
Thus, on every open boundary interval on which the Ernst denominator is
nonzero, \(e^{2U}\) and \(e^{2k}\) are simultaneously positive,
simultaneously zero, or simultaneously negative.
\end{corollary}

\begin{proof}
By Proposition~4.6,
\begin{equation}
 e^{2U^{I_m}}(\zeta)
 =
 \frac{\det X(\zeta)\,|H_N^{I_m}|^2}
 {|\Delta_N^{I_m}(\zeta)|^2}.
 \label{eq:e2UForSignAgreement}
\end{equation}
For \(\zeta\in I_m=(K_{m+1},K_m)\), the factors \(K_i-\zeta\) are
positive for \(i\leq m\) and negative for \(i>m\).  Hence
\begin{equation}
 \operatorname{sgn}\det X(\zeta)
 =(-1)^{2N-m}=(-1)^m.
 \label{eq:detXSignForAgreement}
\end{equation}
If \(H_N^{I_m}\ne0\), equations
\eqref{eq:e2UForSignAgreement} and
\eqref{eq:detXSignForAgreement} give
\begin{equation}
 \operatorname{sgn}e^{2U^{I_m}}(\zeta)=(-1)^m.
\end{equation}
On the other hand, Proposition~4.7 and
\eqref{eq:upperAxisKNormalization} give
\begin{equation}
 e^{2k^{I_m}}
 =\left(\frac{H_N^{I_m}}{H_N^+}\right)^2.
 \label{eq:e2kForSignAgreement}
\end{equation}
Corollary~4.8 therefore implies
\begin{equation}
 \operatorname{sgn}e^{2k^{I_m}}=(-1)^m.
\end{equation}
This proves \eqref{eq:boundarySignsAgree} when \(H_N^{I_m}\ne0\).

If \(H_N^{I_m}=0\), Proposition~4.6 gives
\begin{equation}
 e^{2U^{I_m}}(\zeta)=0,
\end{equation}
while \eqref{eq:e2kForSignAgreement} gives
\begin{equation}
 e^{2k^{I_m}}=0.
\end{equation}
Thus the signs also agree in the vanishing case, with the convention
\(\operatorname{sgn}0=0\).  This proves
\eqref{eq:boundarySignsExplicit}.
\end{proof}

\subsection{Boundary values of \texorpdfstring{\(a\)}{a}}

We next determine the boundary values of the metric function \(a\).  As
in the preceding subsection, we begin with a reconstruction formula
obtained from the LP.  Hennig and Neugebauer show that
\begin{equation}
 a+C
 =i e^{-2U}
 \lim_{\substack{K\to\infty\\ \lambda\to+1}}
 K^2\frac{d}{dK}
 \left[\chi(-\lambda)-\psi(-\lambda)\right].
 \label{eq:HNreconstructionOfa}
\end{equation}
This is equation~(30) of Hennig and Neugebauer
\cite{NeugebauerHennig2009}.  Here the derivative is taken on the
compactified spectral fibre, with \((\rho,\zeta)\) fixed, and on the
sheet on which
\begin{equation}
 K\to\infty,
 \qquad
 \lambda\to+1.
\end{equation}
The constant \(C\in\mathbb R\) reflects the additive freedom in \(a\).

Define \(\mathcal L\) by the expansion
\begin{equation}
 \chi_N(-\lambda)
 =1+\frac{\mathcal L}{\Delta_N}\frac1K+O(K^{-2})
 \label{eq:chiMinusLambdaLaurentForA}
\end{equation}
as \(K\to\infty\) on this sheet.  The reality symmetry of the normalized
matrix gives
\begin{equation}
 \psi_N(-\lambda)
 =1+\frac{\overline{\mathcal L}}
          {\overline{\Delta_N}}\frac1K+O(K^{-2}).
 \label{eq:psiMinusLambdaLaurentForA}
\end{equation}
Consequently,
\begin{equation}
 \chi_N(-\lambda)-\psi_N(-\lambda)
 =
 \left(
 \frac{\mathcal L}{\Delta_N}
 -\frac{\overline{\mathcal L}}{\overline{\Delta_N}}
 \right)\frac1K+O(K^{-2}).
 \label{eq:chiPsiDifferenceLaurentForA}
\end{equation}
Since
\begin{equation}
 K^2\frac{d}{dK}\left(\frac1K\right)=-1,
\end{equation}
formula \eqref{eq:HNreconstructionOfa} becomes
\begin{equation}
 a+C
 =-i e^{-2U}
 \left(
 \frac{\mathcal L}{\Delta_N}
 -\frac{\overline{\mathcal L}}{\overline{\Delta_N}}
 \right).
 \label{eq:aThroughLaurentCoefficient}
\end{equation}
Using
\begin{equation}
 e^{2U}
 =\frac{q_1\overline{q_1}-q_2\overline{q_2}}
        {|\Delta_N|^2},
 \label{eq:e2UGlobalCofactorForm}
\end{equation}
we obtain the algebraic reconstruction formula
\begin{equation}
 a+C
 =-i\,
 \frac{
 \mathcal L\,\overline{\Delta_N}
 -\overline{\mathcal L}\,\Delta_N
 }{
 q_1\overline{q_1}-q_2\overline{q_2}
 }.
 \label{eq:aUncancelledReconstruction}
\end{equation}
Formula \eqref{eq:aUncancelledReconstruction} is valid wherever
\begin{equation}
 \Delta_N\ne0,
 \qquad
 q_1\overline{q_1}-q_2\overline{q_2}\ne0.
\end{equation}
It is the starting point for the boundary calculation below.  We shall
first compute the coefficient \(\mathcal L\) from the determinant defining
\(\chi_N\), and then evaluate the numerator and denominator of
\eqref{eq:aUncancelledReconstruction} as \(\rho\downarrow0\) on an open
boundary interval \(I\).

\begin{proposition}[Boundary values of \(a\)]
\label{prop:boundarya}
Assume \(N\geq2\), and let \(I\) be one of the open boundary intervals.
Define
\begin{equation}
 I_{\mathrm{reg}}
 :=\{\zeta\in I:\Delta_N^I(\zeta)\ne0\}.
 \label{eq:IregularForA}
\end{equation}
Assume that
\begin{equation}
 H_N^I\ne0.
 \label{eq:HnonzeroForA}
\end{equation}
Recall that
\begin{equation}
 H_N^I
 =[c^IK^{N-1},dK^{N-1},c^IK^{N-2},dK^{N-2},
   \ldots,c^IK,dK,c^I,d].
 \label{eq:HrecalledForA}
\end{equation}
Put
\begin{equation}
 v^I:=c^I+d
 \label{eq:vIForA}
\end{equation}
and define
\begin{equation}
 D_N^I
 :=[v^IK^N,v^IK^{N-1},\ldots,v^IK,v^I,
    dK^{N-2},\ldots,dK,d].
 \label{eq:DNIForA}
\end{equation}
Then, for every \(\zeta\in I_{\mathrm{reg}}\), the boundary value
\(a^I(\zeta)\) exists and satisfies
\begin{equation}
 a^I(\zeta)+C
 =(-1)^{1+\frac{(N-1)(N-2)}2}\,
 i\,\frac{D_N^I}{H_N^I}.
 \label{eq:aIgeneralFormula}
\end{equation}
The right-hand side is independent of \(\zeta\).  Consequently,
\(a^I\) is constant on \(I_{\mathrm{reg}}\).  In particular, if
\(I_{\mathrm{reg}}=I\), then \(a^I\) is constant on the whole open
interval \(I\).
\end{proposition}

\begin{proof}
Fix \(\zeta\in I_{\mathrm{reg}}\), suppress the superscript \(I\), and
write
\begin{equation}
 X=K-\zeta\mathbf1,
 \qquad c_j=cX^j,
 \qquad d_j=dX^j,
 \qquad v_j=vX^j.
 \label{eq:XRowsForA}
\end{equation}
Since \(K=X+\zeta\mathbf1\), unitriangular row operations give
\begin{equation}
 H_N=[c_{N-1},d_{N-1},\ldots,c_1,d_1,c,d]
 \label{eq:HNinXForA}
\end{equation}
and
\begin{equation}
 D_N=[v_N,v_{N-1},\ldots,v_1,v,d_{N-2},\ldots,d_1,d].
 \label{eq:DNinXForA}
\end{equation}

Let \(q_1,\ldots,q_{2N+1}\) be the signed cofactors of the free column
of \(\mathcal M_N\), numbered from top to bottom, and put
\begin{equation}
 \Delta=q_1-q_2.
\end{equation}
The expansion of the free spectral column at \(K=\infty\) gives the
boundary value of the Laurent coefficient in
\eqref{eq:chiMinusLambdaLaurentForA} as
\begin{equation}
 \mathcal L=q_3+\zeta q_2-q_4.
 \label{eq:LThroughq234}
\end{equation}
On \(I\), the fixed columns of \(\mathcal M_N\) have successive rows
\begin{equation}
 dK^N,\ cXK^{N-1},\ dK^{N-1},\ cXK^{N-2},\ dK^{N-2},
 \ldots,cX,\ d.
 \label{eq:fixedRowsForA}
\end{equation}
Using \(K=X+\zeta\mathbf1\), the first four cofactors become
\begin{align}
 q_1={}&[c_N,d_{N-1},c_{N-1},d_{N-2},c_{N-2},\ldots,c_1,d],
 \label{eq:q1ForA}\\
 q_2={}&-[d_N,d_{N-1},c_{N-1},d_{N-2},c_{N-2},\ldots,c_1,d],
 \label{eq:q2ForA}\\
 q_3={}&[d_N+N\zeta d_{N-1},c_N,c_{N-1},d_{N-2},c_{N-2},
          \ldots,c_1,d],
 \label{eq:q3ForA}\\
 q_4={}&-[d_N,c_N+(N-1)\zeta c_{N-1},d_{N-1},d_{N-2},c_{N-2},
          \ldots,c_1,d].
 \label{eq:q4ForA}
\end{align}
Multilinearity gives
\begin{equation}
 \Delta=[v_N,d_{N-1},c_{N-1},d_{N-2},c_{N-2},\ldots,c_1,d]
 \label{eq:DeltaForA}
\end{equation}
and
\begin{equation}
 \mathcal L=\mathcal L_0-N\zeta\Delta,
 \qquad
 \mathcal L_0
 :=[d_N,c_N,v_{N-1},d_{N-2},c_{N-2},\ldots,c_1,d].
 \label{eq:L0ForA}
\end{equation}
Since \(\zeta\) is real,
\begin{equation}
 \mathcal L\overline\Delta-\overline{\mathcal L}\Delta
 =\mathcal L_0\overline\Delta-\overline{\mathcal L_0}\Delta.
 \label{eq:zetaCancelsForA}
\end{equation}

The denominator in \eqref{eq:aUncancelledReconstruction} has already
been evaluated in the proof of Proposition~4.6.  With
\begin{equation}
 U_0:=\prod_{i=1}^{2N}u_i,
\end{equation}
equation \eqref{eq:e2Udetidentity} reads
\begin{equation}
 U_0(q_1\overline{q_1}-q_2\overline{q_2})
 =(-1)^N\det(X)H_N^2.
 \label{eq:denominatorRecalledForA}
\end{equation}
Since \(\zeta\in I\), one has \(\det X\ne0\); moreover \(U_0\ne0\).
The assumption \(H_N\ne0\) therefore implies
\begin{equation}
 q_1\overline{q_1}-q_2\overline{q_2}\ne0.
\end{equation}
Thus the un-cancelled reconstruction formula
\eqref{eq:aUncancelledReconstruction} is well defined.

It remains to evaluate its numerator.  Introduce the ordered row
families
\begin{equation}
 \mathcal P=(c_{N-2},d_{N-2},\ldots,c_1,d_1)
 \label{eq:PfamilyForA}
\end{equation}
and
\begin{equation}
 \mathcal Q=(v_N,v_{N-1},\mathcal P).
 \label{eq:QfamilyForA}
\end{equation}
Using \eqref{eq:zetaCancelsForA}, the conjugation rule obtained by
multiplying the \(i\)-th column by \(u_i\), and \(v_j=c_j+d_j\), direct
expansion gives
\begin{align}
 &U_0(\mathcal L\overline\Delta-\overline{\mathcal L}\Delta)
 \notag\\
 &\quad=(-1)^{N-1}\Bigl(
 [\mathcal Q,c_N,d][\mathcal Q,d_{N-1},c]
 -[\mathcal Q,c_N,c][\mathcal Q,d_{N-1},d]
 \Bigr).
 \label{eq:numeratorBeforePluckerForA}
\end{align}
The four determinants have the same \(2N-2\) rows \(\mathcal Q\).
Apply Lemma~\ref{lem:plucker} with
\begin{equation}
 \xi=c_N,
 \qquad \eta=d,
 \qquad \mu=d_{N-1},
 \qquad \nu=c.
\end{equation}
We obtain
\begin{equation}
 U_0(\mathcal L\overline\Delta-\overline{\mathcal L}\Delta)
 =(-1)^N
 [\mathcal Q,c_N,d_{N-1}][\mathcal Q,c,d].
 \label{eq:numeratorAfterPluckerForA}
\end{equation}
Multiplying every column of the determinant defining \(H_N\) by the
corresponding diagonal entry of \(X\) gives
\begin{equation}
 [\mathcal Q,c_N,d_{N-1}]=\det(X)H_N.
 \label{eq:firstNumeratorFactorForA}
\end{equation}
For the second factor, replace \(c_j\) by \(v_j=c_j+d_j\) whenever the
corresponding row \(d_j\) is already present.  These are
determinant-preserving row operations.  Reordering the resulting rows
into the order used in \eqref{eq:DNinXForA} gives
\begin{equation}
 [\mathcal Q,c,d]
 =(-1)^{\frac{(N-1)(N-2)}2}D_N.
 \label{eq:secondNumeratorFactorForA}
\end{equation}
Consequently,
\begin{align}
 U_0(\mathcal L\overline\Delta-\overline{\mathcal L}\Delta)
 ={}&(-1)^{N+\frac{(N-1)(N-2)}2}
 \det(X)H_ND_N.
 \label{eq:numeratorFinalForA}
\end{align}

Divide \eqref{eq:numeratorFinalForA} by
\eqref{eq:denominatorRecalledForA} in
\eqref{eq:aUncancelledReconstruction}.  The factors \(U_0\),
\(\det X\), and \(H_N\) cancel, and one obtains
\begin{equation}
 a^I(\zeta)+C
 =(-1)^{1+\frac{(N-1)(N-2)}2}\,
 i\,\frac{D_N^I}{H_N^I}.
\end{equation}
Both \(D_N^I\) and \(H_N^I\) are expressed solely in terms of the
constant data \(K_i\) and \(u_i^I\), so the right-hand side is independent
of \(\zeta\).  The same formula therefore holds at every point of
\(I_{\mathrm{reg}}\), and \(a^I\) is constant there.
\end{proof}

\begin{remark}[Necessity of the condition \(H_N^I\ne0\)]
\label{rem:vanishingHnonconstantA}
The assumption \(H_N^I\ne0\) in Proposition~\ref{prop:boundarya} is
essential.  If \(H_N^I=0\), the cancelled formula
\[
 a^I+C
 =(-1)^{1+\frac{(N-1)(N-2)}2}
 i\frac{D_N^I}{H_N^I}
\]
is unavailable.  More importantly, even the equilibrium conditions do
not then guarantee that the boundary value of \(a\) is constant on \(I\).

This occurs already in the two-horizon family satisfying the equilibrium
conditions.  We use the
parametrization reviewed in Section~7.2 and choose
\[
 L=1,
 \qquad
 w=2,
 \qquad
 w'=\frac12,
 \qquad
 \varepsilon=1.
\]
Introduce the real Cayley parameter \(t\) by
\[
 \alpha=\frac{1+it}{1-it},
 \qquad |\alpha|=1,
\]
and take
\[
 t=\frac{-26+\sqrt{703}}{27}.
\]
Equivalently, \(t\) is the positive root of
\[
 27t^2+52t-1=0.
\]
The corresponding endpoints are
\[
 K_1=\frac52,
 \qquad
 K_2=1,
 \qquad
 K_3=0,
 \qquad
 K_4=-\frac32,
\]
and the parameters \(\alpha_1,\ldots,\alpha_4\) are obtained from the
equilibrium parametrization in Section~7.2.  Thus these data satisfy the
two-horizon algebraic equilibrium conditions.  Exact computer algebra
gives
\[
 H_2^{\mathcal H^1}=0,
 \qquad
 H_2^+\ne0.
\]
We use the normalization \(a^+=0\).

Since the cancelled formula cannot be used on \(\mathcal H^1\), return
to the un-cancelled reconstruction formula
\begin{equation}
 a+C
 =-i\,
 \frac{
  \mathcal L\overline{\Delta_2}
  -\overline{\mathcal L}\Delta_2
 }{
  q_1\overline{q_1}-q_2\overline{q_2}
 }.
 \label{eq:counterexampleUncancelledA}
\end{equation}
Before using \eqref{eq:counterexampleUncancelledA}, we verify the two
non-vanishing conditions needed below.  First, since \(\Delta_2\) is a
determinant whose entries are
polynomials in the continuous functions
\[
 r_i(\rho,\zeta)=\sqrt{\rho^2+(\zeta-K_i)^2},
\]
it is itself continuous on the closed half-plane \(\rho\geq0\).  At each
of the three rational boundary points
\[
 \zeta_0\in\left\{\frac54,\frac74,\frac94\right\},
\]
exact computer algebra gives
\[
 \Delta_2(0,\zeta_0)\ne0.
\]
By continuity, for each of these points there exists \(\delta>0\) such that
\[
 \Delta_2(\rho,\zeta)\ne0
\]
whenever
\[
 0\le\rho<\delta,
 \qquad
 |\zeta-\zeta_0|<\delta.
\]

Next, put
\[
 \mathcal B
 :=q_1\overline{q_1}-q_2\overline{q_2}.
\]
Its zero-order term vanishes on the degenerate horizon interval, but its
transverse expansion has the form
\[
 \mathcal B(\rho,\zeta)
 =\rho^2B_2(\zeta)+O(\rho^4).
\]
At each of the points \(\zeta_0\) specified above, exact computer algebra
gives \(B_2(\zeta_0)\ne0\).  Hence \(\mathcal B\ne0\)
there for all sufficiently small \(\rho>0\), and
\eqref{eq:counterexampleUncancelledA} may be used along the corresponding
one-sided limiting paths.

If
\[
 \mathcal N
 :=\mathcal L\overline{\Delta_2}
   -\overline{\mathcal L}\Delta_2,
\]
then
\[
 \mathcal N(\rho,\zeta)
 =\rho^2N_2(\zeta)+O(\rho^4),
\]
and therefore
\[
 \left.(a+C)\right|_{\mathcal H^1}(\zeta)
 =-i\frac{N_2(\zeta)}{B_2(\zeta)}
\]
at these points.  The resulting values are
\[
 \left.a\right|_{\mathcal H^1}\!\left(\frac54\right)
 \approx0.6513234978,
\]
\[
 \left.a\right|_{\mathcal H^1}\!\left(\frac74\right)
 \approx-0.6643071267,
\]
and
\[
 \left.a\right|_{\mathcal H^1}\!\left(\frac94\right)
 \approx0.9088050565.
\]
Thus \(a|_{\mathcal H^1}\) is neither zero nor constant.  Consequently,
even for data satisfying the algebraic equilibrium conditions, the
condition \(H_2^{\mathcal H^1}=0\) does not ensure the constant value of
\(a\) required on a candidate Killing-horizon interval.  The assumption
\[
 H_N^{\mathcal H^j}\ne0
\]
must therefore be retained whenever Proposition~\ref{prop:boundarya} is
used to assign a constant value \(a^{\mathcal H^j}\).

The determinant factorizations, the zero and non-zero tests at the three
rational points, and the transverse limits in this example were verified
exactly in Maple; no floating-point specialization was used for those
decisions.  The decimal values above are included only to display the
variation of \(a\) along the interval.
\end{remark}

\begin{remark}[A supplementary condition on the degenerate branches of the
equilibrium family]
In the two-horizon parametrization of Section~7.2, put
\[
 \mathit\Phi:=\cos\phi+\varepsilon\sin\phi,
 \qquad
 \alpha=e^{i\phi}=\frac{1+it}{1-it}.
\]
Exact factorization gives
\begin{equation}
 H_2^{\mathcal H^1}
 =(w'+\mathit\Phi)Q^{(1)}_\varepsilon,
 \label{eq:H1balancedfactorization}
\end{equation}
where, away from the non-zero rational prefactor of
\(Q^{(1)}_\varepsilon\), its additional zero locus is defined by
\begin{align}
 P^{(1)}_\varepsilon(t)
 :={}&w'(1+w')(w+1)^2t^2 \notag\\
 &+2\varepsilon(w^2w'+2w+w')t
   +w'(w'-1)(w-1)^2=0.
 \label{eq:extraH1factor}
\end{align}
The example in Remark~\ref{rem:vanishingHnonconstantA} was obtained by
choosing a root of \eqref{eq:extraH1factor} while keeping
\(w'+\mathit\Phi\ne0\).  For the displayed values
\(w=2\), \(w'=1/2\), and \(\varepsilon=1\), equation
\eqref{eq:extraH1factor} is, after multiplication by (4), precisely
\[
 27t^2+52t-1=0.
\]

The other factor in \eqref{eq:H1balancedfactorization} behaves
differently.  The transverse expansion of the un-cancelled formula
\eqref{eq:counterexampleUncancelledA} shows that, under the normalization
\(a^+=0\),
\begin{equation}
 w'+\mathit\Phi=0
 \quad\Longrightarrow\quad
 a|_{\mathcal H^1}=0.
 \label{eq:H1supplementarycondition}
\end{equation}
This suggests the following small supplement to the equilibrium
parametrization: whenever \(H_2^{\mathcal H^1}=0\), require
\(w'+\mathit\Phi=0\).  It is a sufficient condition for the desired
constant value of (a), but it is not asserted to be necessary.

The second horizon has a completely analogous factorization.  Exact
computer algebra gives
\begin{equation}
 H_2^{\mathcal H^2}
 =(1-w\mathit\Phi)Q^{(2)}_\varepsilon.
 \label{eq:H2balancedfactorization}
\end{equation}
The analogous supplement is therefore: whenever
\(H_2^{\mathcal H^2}=0\), require
\begin{equation}
 1-w\mathit\Phi=0.
 \label{eq:H2supplementarycondition}
\end{equation}
The same transverse calculation, with the roles of the two horizon
intervals interchanged, gives
\begin{equation}
 1-w\mathit\Phi=0
 \quad\Longrightarrow\quad
 a|_{\mathcal H^2}=0
 \label{eq:H2zeroa}
\end{equation}
whenever the non-vanishing conditions needed for the un-cancelled
reconstruction formula hold.  This is used only as a sufficient
condition; necessity is not claimed.  All factorizations
and transverse identities in this remark were verified exactly in Maple.
\end{remark}

\section{Equilibrium conditions and branch selection}

\subsection{Equilibrium conditions}

Throughout this section, we assume that the boundary denominator has no zero on
the line \(\rho=0\):
\begin{equation}
 \Delta_N(0,\zeta)\ne0
 \qquad (\zeta\in\mathbb R).
 \label{eq:boundaryDeltaGlobalNonzero}
\end{equation}
On each open boundary interval \(I\), the expression in
\eqref{eq:boundaryDeltaGlobalNonzero} means the corresponding determinant
\(\Delta_N^I(\zeta)\). Thus none of the boundary expressions having
\(\Delta_N^I\) as denominator develops a pole on an open boundary interval.

We impose the further non-degeneracy condition
\begin{equation}
 H_N^+\ne0,
 \qquad H_N^+:=H_N^{\mathcal A^+}.
 \label{eq:HplusNonzeroEquilibrium}
\end{equation}
By Proposition~\ref{prop:boundaryErnst}, this means that the upper-axis Ernst
potential is genuinely represented as a quotient of two monic polynomials of
degree \(N\):
\begin{equation}
 f^+(\zeta)=\frac{P_N^+(\zeta)}{Q_N^+(\zeta)}.
 \label{eq:upperAxisGenuineNN}
\end{equation}
In this sense, \eqref{eq:HplusNonzeroEquilibrium} excludes a degeneration of
the rational upper-axis datum.

The additive constant in \(k\) is now fixed by the upper-axis normalization
\begin{equation}
 e^{2k^+}=1,
 \label{eq:kUpperAxisNormalization}
\end{equation}
and the additive constant in the reconstruction of \(a\) is fixed by
\begin{equation}
 a^+=0.
 \label{eq:aUpperAxisNormalization}
\end{equation}
On any other open axis interval
\(\mathcal A^j\), \(j\in\{1,\ldots,N\}\), the standard regular-axis
conditions are
\begin{equation}
 a^{\mathcal A^j}=0,
 \qquad
 e^{2k^{\mathcal A^j}}=1.
 \label{eq:regularAxisConditions}
\end{equation}
The second equality is the elementary-flatness condition and excludes a
conical singularity, or strut, on \(\mathcal A^j\); the first excludes the
corresponding torsion singularity.

Propositions~\ref{prop:boundarye2k} and~\ref{prop:boundarya} convert
\eqref{eq:regularAxisConditions} into determinant equations. The equilibrium
conditions are
\begin{equation}
 a^{\mathcal A^j}=a^+,
 \qquad j\in\{1,\ldots,N\},
 \label{eq:aEquilibriumConditions}
\end{equation}
and
\begin{equation}
 \bigl(H_N^{\mathcal A^j}\bigr)^2=(H_N^+)^2,
 \qquad j\in\{1,\ldots,N\}.
 \label{eq:HsquareEquilibriumConditions}
\end{equation}
Equivalently, \eqref{eq:aEquilibriumConditions} can be written as
\begin{equation}
 D_N^{\mathcal A^j}H_N^+
 -D_N^+H_N^{\mathcal A^j}=0,
 \qquad j\in\{1,\ldots,N\}.
 \label{eq:DEquilibriumConditions}
\end{equation}
At this stage it is important to retain the form
\eqref{eq:DEquilibriumConditions}: replacing it prematurely by an equation
involving only sums or differences of the \(D_N^I\) would already amount to
choosing a sign in \eqref{eq:HsquareEquilibriumConditions}.

There is one automatic relation among the elementary-flatness equations. On
the upper axis one has
\begin{equation}
 \lambda_i^+=-1,
 \qquad u_i^+=-\alpha_i,
 \qquad i\in\{1,\ldots,2N\},
\end{equation}
whereas on the lower axis
\begin{equation}
 \lambda_i^-=1,
 \qquad u_i^-=\alpha_i,
 \qquad i\in\{1,\ldots,2N\}.
\end{equation}
Hence \(c^-=-c^+\), while the row \(d\) is unchanged. Since the determinant
defining \(H_N^I\) contains exactly \(N\) rows of \(c\)-type, it follows that
\begin{equation}
 H_N^-=(-1)^N H_N^+.
 \label{eq:HminusHplusRelation}
\end{equation}
Consequently,
\begin{equation}
 (H_N^-)^2=(H_N^+)^2
\end{equation}
holds identically. Thus \eqref{eq:aEquilibriumConditions} supplies \(N\)
equations, whereas at most \(N-1\) of the equations in
\eqref{eq:HsquareEquilibriumConditions} are independent. Altogether, there
are at most \(2N-1\) independent equilibrium equations.

For each \(j\in\{1,\ldots,N\}\), the horizon interval \(\mathcal H^j\) lies
between the adjacent axis intervals \(\mathcal A^{j-1}\) and
\(\mathcal A^j\). Equation \eqref{eq:HsquareEquilibriumConditions} leaves two
apparently possible branches across \(\mathcal H^j\):
\begin{equation}
 H_N^{\mathcal A^j}=H_N^{\mathcal A^{j-1}}
 \label{eq:equalHBranch}
\end{equation}
or
\begin{equation}
 H_N^{\mathcal A^j}=-H_N^{\mathcal A^{j-1}}.
 \label{eq:alternatingHBranch}
\end{equation}
The purpose of the following subsections is to exclude
\eqref{eq:equalHBranch} and thereby establish the alternating branch
\eqref{eq:alternatingHBranch}.

\subsection{Branch selection}

We begin by studying the change of the determinants \(H_N^I\) and
\(D_N^I\) across a single endpoint \(K_k\). For
\(k\in\{1,\ldots,2N\}\), the two open boundary intervals adjacent to
\(K_k\) are \(I_{k-1}\), immediately above \(K_k\), and \(I_k\),
immediately below it. Define the endpoint jump
\begin{equation}
 M_k
 :=
 D_N^{I_k}H_N^{I_{k-1}}
 -
 D_N^{I_{k-1}}H_N^{I_k}.
 \label{eq:Mkdefinition}
\end{equation}

Recall from Section~4.1 the reduced endpoint objects
\(X_k,c^{(k)},d^{(k)},v^{(k)},E_k\), and \(\Lambda_k\), defined in
\eqref{eq:Xkendpoint}--\eqref{eq:EkendpointNone} and
\eqref{eq:Lambdakendpoint}.

\begin{proposition}[Factorization of the endpoint jumps]
\label{prop:endpointJumpFactorization}
For every \(k\in\{1,\ldots,2N\}\),
\begin{equation}
 M_k
 =
 2(-1)^{\frac{(N-1)(N-2)}2}
 \alpha_k\det(X_k)E_k^2.
 \label{eq:MkEkfactorization}
\end{equation}
Let
\begin{equation}
 A:=\prod_{i=1}^{2N}\alpha_i.
 \label{eq:Aallalpha}
\end{equation}
Then there exists a real number \(R_k\), determined up to sign, such that
\begin{equation}
 M_k
 =
 2(-1)^{N+1+\frac{(N-1)(N-2)}2}
 A\Lambda_kR_k^2.
 \label{eq:MkRkfactorization}
\end{equation}
In particular, all the \(M_k\) have the common complex factor
\begin{equation}
 2(-1)^{N+1+\frac{(N-1)(N-2)}2}A,
 \label{eq:Mkcommonphase}
\end{equation}
and, after this factor is removed, the remaining quantities are
nonnegative. Under the assumptions of Section~5.1 they are strictly
positive:
\begin{equation}
 \Lambda_kR_k^2>0.
 \label{eq:Mkstrictlypositive}
\end{equation}
\end{proposition}

\begin{proof}
Only the \(k\)-th component of the row \(c^I\) changes when one passes
from \(I_{k-1}\) to \(I_k\). More precisely,
\begin{equation}
 u_k^{I_{k-1}}=-\alpha_k,
 \qquad
 u_k^{I_k}=\alpha_k,
 \label{eq:ukendpointchange}
\end{equation}
whereas \(u_i^{I_{k-1}}=u_i^{I_k}\) for \(i\ne k\).

Temporarily replace the \(k\)-th entry by an indeterminate \(t\), leaving
all the other entries fixed. Since a determinant is affine in each
column separately, there are coefficients \(h_0,h_1,p_0,p_1\),
independent of \(t\), such that
\begin{equation}
 H(t)=h_0+th_1,
 \qquad
 D(t)=p_0+tp_1.
 \label{eq:HtDtendpoint}
\end{equation}
Here
\begin{align}
 H(-\alpha_k)&=H_N^{I_{k-1}},
 &H(\alpha_k)&=H_N^{I_k},\\
 D(-\alpha_k)&=D_N^{I_{k-1}},
 &D(\alpha_k)&=D_N^{I_k}.
\end{align}
Consequently,
\begin{align}
 M_k
 &=
 (p_0+\alpha_kp_1)(h_0-\alpha_kh_1)
 -(p_0-\alpha_kp_1)(h_0+\alpha_kh_1)
 \notag\\
 &=
 2\alpha_k(p_1h_0-p_0h_1).
 \label{eq:Mkbilinearcofactors}
\end{align}

The coefficients \(h_0,h_1,p_0,p_1\) are obtained by Laplace expansion
in the \(k\)-th column. The position of that column contributes the same
cofactor sign to all four coefficients. Since every term in
\(p_1h_0-p_0h_1\) is a product of two such cofactors, this common sign
occurs twice and cancels.

We may therefore delete the \(k\)-th column, replace powers of \(K\) by
powers of \(K-K_k\mathbf1\) using unitriangular row operations, and work
with the reduced rows defined by \(X_k\). Suppress the superscript
\((k)\) temporarily and introduce the ordered row family
\begin{equation}
 \mathcal W
 =
 (v_{N-1},c_{N-2},d_{N-2},\ldots,c_1,d_1).
 \label{eq:WendpointJump}
\end{equation}
For \(N=2\), this means simply \(\mathcal W=(v_1)\). The case \(N=1\)
follows directly from the corresponding \(1\times1\) determinants; in
the remainder of the proof we assume \(N\ge2\).

Laplace expansion in the deleted column, followed by the identities
\(v_m=c_m+d_m\), gives
\begin{align}
 p_1h_0-p_0h_1
 ={}&
 (-1)^{1+\frac{(N-1)(N-2)}2}
 \Bigl(
 [\mathcal W,v_N,v][\mathcal W,d_{N-1},d]
 \notag\\
 &\hspace{39mm}
 -[\mathcal W,v_N,d][\mathcal W,d_{N-1},v]
 \Bigr).
 \label{eq:Mkfourdeterminants}
\end{align}
The four determinants in parentheses have the same \(2N-3\) rows
\(\mathcal W\). Lemma~\ref{lem:plucker} therefore yields
\begin{equation}
 p_1h_0-p_0h_1
 =
 (-1)^{1+\frac{(N-1)(N-2)}2}
 [\mathcal W,v_N,d_{N-1}]
 [\mathcal W,v,d].
 \label{eq:MkonePlucker}
\end{equation}

We now identify the two factors. Determinant-preserving row operations
based on \(v_m=c_m+d_m\) give
\begin{equation}
 [\mathcal W,v,d]
 =
 (-1)^{\frac{(N-1)(N-2)}2}E_k.
 \label{eq:Mksecondfactor}
\end{equation}
Multiplying every column by the corresponding diagonal entry of \(X_k\),
and then performing the same row operations, gives
\begin{equation}
 [\mathcal W,v_N,d_{N-1}]
 =
 (-1)^{1+\frac{(N-1)(N-2)}2}
 \det(X_k)E_k.
 \label{eq:Mkfirstfactor}
\end{equation}
Substitution into \eqref{eq:MkonePlucker} gives
\begin{equation}
 p_1h_0-p_0h_1
 =
 (-1)^{\frac{(N-1)(N-2)}2}
 \det(X_k)E_k^2.
 \label{eq:Mkcofactorfactorized}
\end{equation}
Together with \eqref{eq:Mkbilinearcofactors}, this proves
\eqref{eq:MkEkfactorization}.

It remains to remove the phase of \(E_k\). Put
\begin{equation}
 U_k
 :=
 \prod_{i\ne k}u_i^{I_{k-1}}.
 \label{eq:UkendpointJump}
\end{equation}
Among the remaining \(2N-1\) entries, precisely \(2N-k\) carry a minus
sign. Since \(2N\) is even,
\begin{equation}
 U_k=(-1)^k\frac{A}{\alpha_k}.
 \label{eq:UkthroughA}
\end{equation}
Columnwise conjugation of \(E_k\), using
\(\overline{\alpha_i}=\alpha_i^{-1}\), gives
\begin{equation}
 U_k\overline{E_k}=(-1)^{N-1}E_k.
 \label{eq:EkconjugationEndpointJump}
\end{equation}
Set
\begin{equation}
 S_k:=(-1)^{N-1}U_k.
 \label{eq:SkendpointJump}
\end{equation}
Then \(S_k\overline{E_k}=E_k\). Since \(|S_k|=1\), choose either square
root of \(S_k\) and put
\begin{equation}
 R_k:=\frac{E_k}{\sqrt{S_k}}.
 \label{eq:Rkdefinition}
\end{equation}
The conjugation relation implies \(\overline{R_k}=R_k\), so
\(R_k\in\mathbb R\). Changing the square root of \(S_k\) changes only
the sign of \(R_k\). Equations \eqref{eq:UkthroughA} and
\eqref{eq:SkendpointJump} give
\begin{equation}
 \alpha_kE_k^2=(-1)^{N-1+k}AR_k^2.
 \label{eq:alphaEkphaseRemoved}
\end{equation}

Moreover, \(X_k\) has precisely \(2N-k\) negative diagonal entries.
Therefore
\begin{equation}
 \det(X_k)=(-1)^k\Lambda_k.
 \label{eq:detXkthroughLambda}
\end{equation}
Substitution of \eqref{eq:alphaEkphaseRemoved} and
\eqref{eq:detXkthroughLambda} into \eqref{eq:MkEkfactorization} proves
\eqref{eq:MkRkfactorization}.

Finally, the standing assumption \eqref{eq:boundaryDeltaGlobalNonzero}
includes the endpoint value \(\Delta_k\ne0\). By
Proposition~\ref{prop:endpointdenominator},
\begin{equation}
 \Delta_k
 =
 (-1)^{1+\frac{(N-2)(N-3)}2}\Lambda_kE_k.
 \label{eq:DeltakForStrictJump}
\end{equation}
Since \(\Lambda_k>0\), it follows that \(E_k\ne0\), hence \(R_k\ne0\).
This proves \eqref{eq:Mkstrictlypositive}.
\end{proof}

\begin{theorem}[Branch selection and determinant equilibrium conditions]
\label{thm:branchSelection}
Assume that
\begin{equation}
 \Delta_N(0,\zeta)\ne0
 \qquad(\zeta\in\mathbb R)
 \label{eq:mainDeltaNonzero}
\end{equation}
and
\begin{equation}
 H_N^+\ne0.
 \label{eq:mainHplusNonzero}
\end{equation}
Fix the normalization
\begin{equation}
 a^+=0,
 \qquad
 e^{2k^+}=1.
 \label{eq:mainUpperNormalization}
\end{equation}
Suppose that the regular-axis equilibrium conditions
\begin{equation}
 a^{\mathcal A^j}=0,
 \qquad
 e^{2k^{\mathcal A^j}}=1,
 \qquad
 j\in\{1,\ldots,N\},
 \label{eq:mainRegularAxisConditions}
\end{equation}
hold. Then, for every \(j\in\{1,\ldots,N\}\), the equal-sign branch
\begin{equation}
 H_N^{\mathcal A^j}=H_N^{\mathcal A^{j-1}}
 \label{eq:mainEqualSignImpossible}
\end{equation}
is impossible. Consequently,
\begin{equation}
 H_N^{\mathcal A^j}
 =
 -H_N^{\mathcal A^{j-1}},
 \qquad
 j\in\{1,\ldots,N\},
 \label{eq:mainAdjacentAlternation}
\end{equation}
and hence
\begin{equation}
 H_N^{\mathcal A^j}
 =
 (-1)^jH_N^+,
 \qquad
 j\in\{1,\ldots,N\}.
 \label{eq:mainGlobalAlternation}
\end{equation}
The last equation in \eqref{eq:mainGlobalAlternation}, corresponding to
\(j=N\), is the automatic identity
\(H_N^-=(-1)^NH_N^+\).

Therefore, under \eqref{eq:mainDeltaNonzero} and
\eqref{eq:mainHplusNonzero}, the regular-axis equilibrium conditions
\eqref{eq:mainRegularAxisConditions} are equivalent to
\begin{equation}
 D_N^{\mathcal A^j}
 =
 (-1)^jD_N^+,
 \qquad
 j\in\{1,\ldots,N\},
 \label{eq:mainDconditions}
\end{equation}
together with
\begin{equation}
 H_N^{\mathcal A^j}
 =
 (-1)^jH_N^+,
 \qquad
 j\in\{1,\ldots,N-1\}.
 \label{eq:mainHconditions}
\end{equation}
\end{theorem}

\begin{proof}
Fix \(j\in\{1,\ldots,N\}\). For readability, write
\begin{align}
 H_{\mathrm{up}}&:=H_N^{\mathcal A^{j-1}},
 &D_{\mathrm{up}}&:=D_N^{\mathcal A^{j-1}},\notag\\
 H_{\mathrm{middle}}&:=H_N^{\mathcal H^j},
 &D_{\mathrm{middle}}&:=D_N^{\mathcal H^j},\notag\\
 H_{\mathrm{down}}&:=H_N^{\mathcal A^j},
 &D_{\mathrm{down}}&:=D_N^{\mathcal A^j}.
 \label{eq:upMiddleDownAbbreviations}
\end{align}
The two endpoint jumps belonging to \(\mathcal H^j\) are
\begin{equation}
 M_{2j-1}
 =
 D_{\mathrm{middle}}H_{\mathrm{up}}
 -
 D_{\mathrm{up}}H_{\mathrm{middle}},
 \label{eq:MupperHorizonEndpoint}
\end{equation}
and
\begin{equation}
 M_{2j}
 =
 D_{\mathrm{down}}H_{\mathrm{middle}}
 -
 D_{\mathrm{middle}}H_{\mathrm{down}}.
 \label{eq:MlowerHorizonEndpoint}
\end{equation}

The equilibrium condition for \(a\) gives
\(a^{\mathcal A^{j-1}}=a^{\mathcal A^j}\). By
Proposition~\ref{prop:boundarya},
\begin{equation}
 \frac{D_{\mathrm{up}}}{H_{\mathrm{up}}}
 =
 \frac{D_{\mathrm{down}}}{H_{\mathrm{down}}}.
 \label{eq:DHoverAdjacentAxes}
\end{equation}
All denominators in \eqref{eq:DHoverAdjacentAxes} are nonzero:
\(H_N^+\ne0\), and the elementary-flatness equations imply
\begin{equation}
 H_{\mathrm{up}}^2
 =
 H_{\mathrm{down}}^2
 =
 (H_N^+)^2.
 \label{eq:adjacentAxisHsquares}
\end{equation}

Suppose, for contradiction, that the equal-sign branch holds:
\begin{equation}
 H_{\mathrm{down}}=H_{\mathrm{up}}.
 \label{eq:equalBranchContradictionAssumption}
\end{equation}
Equation \eqref{eq:DHoverAdjacentAxes} then gives
\begin{equation}
 D_{\mathrm{down}}=D_{\mathrm{up}}.
 \label{eq:equalBranchDsame}
\end{equation}
Adding \eqref{eq:MupperHorizonEndpoint} and
\eqref{eq:MlowerHorizonEndpoint}, and using
\eqref{eq:equalBranchContradictionAssumption}--\eqref{eq:equalBranchDsame},
we obtain
\begin{align}
 M_{2j-1}+M_{2j}
 ={}&
 D_{\mathrm{middle}}H_{\mathrm{up}}
 -D_{\mathrm{up}}H_{\mathrm{middle}}
 \notag\\
 &+
 D_{\mathrm{up}}H_{\mathrm{middle}}
 -D_{\mathrm{middle}}H_{\mathrm{up}}
 =0.
 \label{eq:equalBranchMsumZero}
\end{align}

On the other hand, Proposition~\ref{prop:endpointJumpFactorization}
gives
\begin{align}
 M_{2j-1}+M_{2j}
 ={}&
 2(-1)^{N+1+\frac{(N-1)(N-2)}2}A
 \notag\\
 &\times
 \left(
 \Lambda_{2j-1}R_{2j-1}^2
 +
 \Lambda_{2j}R_{2j}^2
 \right).
 \label{eq:positiveMsum}
\end{align}
The quantity in parentheses is strictly positive by
\eqref{eq:Mkstrictlypositive}, while \(A\ne0\). Hence the right-hand side
of \eqref{eq:positiveMsum} is nonzero, contradicting
\eqref{eq:equalBranchMsumZero}. Thus the equal-sign branch is impossible,
and \eqref{eq:mainAdjacentAlternation} follows. Iterating from
\(\mathcal A^0=\mathcal A^+\) proves
\eqref{eq:mainGlobalAlternation}.

It remains to rewrite the \(a\)-equilibrium conditions. From
\begin{equation}
 \frac{D_N^{\mathcal A^j}}{H_N^{\mathcal A^j}}
 =
 \frac{D_N^+}{H_N^+}
 \label{eq:mainDHoverH}
\end{equation}
and \eqref{eq:mainGlobalAlternation}, we obtain
\eqref{eq:mainDconditions}.

Conversely, assume \eqref{eq:mainDconditions} and
\eqref{eq:mainHconditions}. The missing case \(j=N\) of
\eqref{eq:mainHconditions} is supplied by the automatic relation
\eqref{eq:HminusHplusRelation}. Therefore
\begin{equation}
 \frac{D_N^{\mathcal A^j}}{H_N^{\mathcal A^j}}
 =
 \frac{D_N^+}{H_N^+},
 \qquad
 j\in\{1,\ldots,N\},
\end{equation}
and Proposition~\ref{prop:boundarya} gives
\(a^{\mathcal A^j}=a^+=0\). Likewise,
\begin{equation}
 \left(
 \frac{H_N^{\mathcal A^j}}{H_N^+}
 \right)^2=1,
\end{equation}
so Proposition~\ref{prop:boundarye2k} gives
\(e^{2k^{\mathcal A^j}}=e^{2k^+}=1\). This proves the claimed
equivalence.
\end{proof}

\begin{remark}
For \(N=2\), the exclusion of the equal-sign branch was established by
Tomimatsu and Kihara~\cite{TomimatsuKihara1982}. For general \(N\),
Yamazaki~\cite{Yamazaki1983}, using a different parametrization and
different determinant formulas, obtained conditions equivalent to
\eqref{eq:HsquareEquilibriumConditions} and
\eqref{eq:DEquilibriumConditions}. He also displayed the reduced
equilibrium conditions corresponding to
\eqref{eq:mainDconditions}--\eqref{eq:mainHconditions}. His treatment,
however, did not prove that the equal-sign alternatives allowed by the
squared conditions are impossible, and hence did not establish that the
alternating branch is forced. Theorem~\ref{thm:branchSelection} supplies
this missing branch-selection argument for arbitrary \(N\).
\end{remark}

\section{Thermodynamics of the solution}

Following Section~3 of Neugebauer and Hennig~\cite{NeugebauerHennig2009},
we associate algebraic horizon quantities with each interval
\begin{equation}
 \mathcal H^j=(K_{2j},K_{2j-1}),
 \qquad
 \ell_j:=K_{2j-1}-K_{2j}>0,
 \qquad
 j\in\{1,\ldots,N\}.
 \label{eq:horizonLengthThermo}
\end{equation}
They will be denoted by \(\kappa_j,\Omega_j,A_j,\mathfrak M_j\), and
\(J_j\), and are called the candidate surface gravity, angular
velocity, area, Komar mass, and angular momentum, respectively.  The
word ``candidate'' is essential: the determinant construction by itself
does not establish that the intervals \(\mathcal H^j\) extend to regular
Killing horizons, nor does it provide the spacetime regularity needed
for the usual geometric and Komar-integral interpretations of these
quantities.

Assume throughout this section that the endpoint values
\begin{equation}
 f_m:=f(0,K_m),
 \qquad m\in\{1,\ldots,2N\},
 \label{eq:endpointValuesThermo}
\end{equation}
exist and that all determinants occurring in the denominators below are
nonzero.  We use the algebraic definitions
\begin{align}
 \kappa_j+i\Omega_j
 &:=\frac12
 \left(f^{\mathcal A^{j-1}}\right)_{,\zeta}(K_{2j-1}),
 \label{eq:kappaOmegaDefinition}\\
 \kappa_jA_j&:=2\pi\ell_j,
 \label{eq:areaDefinition}\\
 \mathfrak M_j
 &:=\frac{i}{4\Omega_j}(f_{2j-1}-f_{2j}),
 \label{eq:massDefinition}\\
 J_j
 &:=\frac1{2\Omega_j}
 \left(\mathfrak M_j-\frac{\ell_j}{2}\right).
 \label{eq:JDefinition}
\end{align}
Equations \eqref{eq:areaDefinition} and \eqref{eq:JDefinition} imply
immediately the Smarr relation
\begin{equation}
 \mathfrak M_j
 =2\Omega_jJ_j+\frac{\kappa_jA_j}{4\pi}.
 \label{eq:SmarrCandidateRelation}
\end{equation}
The dimensionless ratio relevant to the subextremal inequality is
\begin{equation}
 p_j:=\frac{8\pi J_j}{A_j}.
 \label{eq:pjDefinition}
\end{equation}

For later use, let
\begin{equation}
 \mathcal V_{N-1}^{(j)}
 :=
 \bigl(
 c^{\mathcal H^j}K^{N-2},dK^{N-2},\ldots,
 c^{\mathcal H^j}K,dK,c^{\mathcal H^j},d
 \bigr)
 \label{eq:reducedRowFamilyThermo}
\end{equation}
and define the reduced determinant
\begin{equation}
 H_{N-1}^{(j)}
 :=
 [\mathcal V_{N-1}^{(j)}]_{\widehat{2j-1,2j}}.
 \label{eq:reducedHThermo}
\end{equation}
Thus the columns numbered \(2j-1\) and \(2j\) are deleted before the
determinant is taken.  For \(N=1\), we use the convention
\begin{equation}
 H_0^{(1)}:=1.
 \label{eq:H0Convention}
\end{equation}

\begin{proposition}[Candidate horizon quantities]
\label{prop:candidateHorizonQuantities}
For every \(j\in\{1,\ldots,N\}\), the quantities defined above satisfy
\begin{align}
 \kappa_j+i\Omega_j
 &=\frac{H_N^{\mathcal A^{j-1}}
 \bigl(H_N^{\mathcal A^{j-1}}-H_N^{\mathcal H^j}\bigr)}
 {2\alpha_{2j-1}\Lambda_{2j-1}E_{2j-1}^2},
 \label{eq:kappaOmegaDeterminant}\\
 \kappa_j
 &=\frac{\bigl(H_N^{\mathcal A^{j-1}}\bigr)^2}
 {2\alpha_{2j-1}\Lambda_{2j-1}E_{2j-1}^2},
 \label{eq:kappaDeterminant}\\
 \Omega_j
 &=i\,\frac{H_N^{\mathcal A^{j-1}}H_N^{\mathcal H^j}}
 {2\alpha_{2j-1}\Lambda_{2j-1}E_{2j-1}^2},
 \label{eq:OmegaDeterminant}\\
 f_{2j-1}-f_{2j}
 &=-2\,\frac{H_N^{\mathcal H^j}H_{N-1}^{(j)}}
 {E_{2j-1}E_{2j}},
 \label{eq:endpointDifferenceThermo}\\
 A_j
 &=4\pi\ell_j\,
 \frac{\alpha_{2j-1}\Lambda_{2j-1}E_{2j-1}^2}
 {\bigl(H_N^{\mathcal A^{j-1}}\bigr)^2},
 \label{eq:areaDeterminant}\\
 \mathfrak M_j
 &=-\alpha_{2j-1}\Lambda_{2j-1}
 \frac{E_{2j-1}H_{N-1}^{(j)}}
 {H_N^{\mathcal A^{j-1}}E_{2j}},
 \label{eq:massDeterminant}\\
 p_j
 &=\frac{i}{H_N^{\mathcal H^j}E_{2j}}
 \left[
 H_N^{\mathcal A^{j-1}}E_{2j}
 +2\alpha_{2j-1}\frac{\Lambda_{2j-1}}{\ell_j}
 E_{2j-1}H_{N-1}^{(j)}
 \right].
 \label{eq:pjDeterminant}
\end{align}
In addition,
\begin{equation}
 \frac{\Omega_j}{\kappa_j}
 =i\,\frac{H_N^{\mathcal H^j}}{H_N^{\mathcal A^{j-1}}},
 \label{eq:OmegaOverKappa}
\end{equation}
and the boundary values of \(a\) obey
\begin{equation}
 a^{\mathcal A^{j-1}}-a^{\mathcal H^j}
 =\frac1{\Omega_j}.
 \label{eq:aJumpThroughOmega}
\end{equation}
\end{proposition}

\begin{proof}[Indication of the calculation]
Differentiate the boundary quotient of Proposition~4.1 at
\(K_{2j-1}\).  After the quotient rule is applied, the two products of
determinants have \(2N-2\) reduced rows in common.  The Pl\"ucker
relation of Lemma~4.4 reduces their difference to the product involving
\(E_{2j-1}\), \(\det(X_{2j-1})\), and a second copy of
\(E_{2j-1}\).  Since
\(\det(X_{2j-1})=-\Lambda_{2j-1}\), this gives
\eqref{eq:kappaOmegaDeterminant}.  The conjugation identities from
Section~4 separate its real and imaginary parts and yield
\eqref{eq:kappaDeterminant}--\eqref{eq:OmegaDeterminant}.

The difference of the two endpoint quotients is handled by the same
bordered-minor Pl\"ucker calculation.  Deleting the endpoint columns
leaves precisely the determinant \(H_{N-1}^{(j)}\), which proves
\eqref{eq:endpointDifferenceThermo}.  Equations
\eqref{eq:areaDeterminant} and \eqref{eq:massDeterminant} then follow by
substitution into \eqref{eq:areaDefinition} and \eqref{eq:massDefinition}.
Substitution of \eqref{eq:JDefinition} into \eqref{eq:pjDefinition},
followed by cancellation of the common determinant factors, gives
\eqref{eq:pjDeterminant}.

Finally, division of \eqref{eq:OmegaDeterminant} by
\eqref{eq:kappaDeterminant} proves \eqref{eq:OmegaOverKappa}.
Combining \eqref{eq:OmegaDeterminant} with the boundary formula for
\(a\) in Proposition~4.10 gives \eqref{eq:aJumpThroughOmega}.
\end{proof}

The sign of the candidate area is already fixed algebraically.  Indeed,
for the odd endpoint \(2j-1\), equation \eqref{eq:alphaEkphaseRemoved}
and the conjugation identity \eqref{eq:HsquareThroughU0} give
\begin{equation}
 \alpha_{2j-1}E_{2j-1}^2=(-1)^NAR_{2j-1}^2,
 \qquad
 \bigl(H_N^{\mathcal A^{j-1}}\bigr)^2
 =(-1)^NA\left|H_N^{\mathcal A^{j-1}}\right|^2.
 \label{eq:areaPositiveFactors}
\end{equation}
Here the second identity uses the fact that \(\mathcal A^{j-1}=I_{2j-2}\),
so that the product of the boundary signs is unchanged from the upper
axis.  Consequently, \eqref{eq:areaDeterminant} becomes
\begin{equation}
 A_j
 =4\pi\ell_j\Lambda_{2j-1}
 \frac{R_{2j-1}^2}
 {\left|H_N^{\mathcal A^{j-1}}\right|^2}>0.
 \label{eq:areaStrictlyPositive}
\end{equation}
This is a purely algebraic positivity statement: it does not prove that
\(A_j\) is the area of a regular horizon section.  It does, however,
show that the determinant expression has the indispensable positivity
property expected of such an area and therefore supports its status as
a natural candidate.

\section{Examples for \texorpdfstring{\(N=1,2,3\)}{N=1,2,3}}

\subsection{The case \texorpdfstring{\(N=1\)}{N=1}}

Let
\[
 K_1>K_2,\qquad |\alpha_1|=|\alpha_2|=1,\qquad
 \ell:=K_1-K_2>0.
\]
Put
\[
 r_i=\sqrt{(K_i-\zeta)^2+\rho^2},\qquad i\in\{1,2\},
\]
and recall that
\[
 \mu=\lambda(K+iz),\qquad \mu^2=(K-\zeta)^2+\rho^2.
\]
For \(N=1\), the free spectral column and the fixed columns are
\[
 V_1=\begin{pmatrix}K\\ \mu\\ 1\end{pmatrix},
 \qquad
 C_i=\begin{pmatrix}K_i\\ \alpha_i r_i\\ 1\end{pmatrix},
 \qquad i\in\{1,2\}.
\]
Thus
\begin{equation}
 \mathcal M_1=
 \begin{pmatrix}
 K&K_1&K_2\\
 \mu&\alpha_1r_1&\alpha_2r_2\\
 1&1&1
 \end{pmatrix}.
 \label{eq:NOneM}
\end{equation}

Let \(q_1,q_2,q_3\) be the signed cofactors of the entries in the first
column. Then
\[
 \det\mathcal M_1=Kq_1+\mu q_2+q_3,
\]
where
\begin{align}
 q_1&=\alpha_1r_1-\alpha_2r_2,\\
 q_2&=-\ell,\\
 q_3&=K_1\alpha_2r_2-K_2\alpha_1r_1.
\end{align}
The determinant denominator is
\begin{equation}
 \Delta_1=
 \begin{vmatrix}
 1&K_1&K_2\\
 -1&\alpha_1r_1&\alpha_2r_2\\
 0&1&1
 \end{vmatrix}
 =q_1-q_2=q_1+\ell.
 \label{eq:NOneDelta}
\end{equation}
Consequently,
\begin{equation}
 \chi_1
 =\frac{Kq_1+\mu q_2+q_3}{K(q_1-q_2)}
 =\frac{Kq_1-\ell\mu+q_3}{K(q_1+\ell)}.
 \label{eq:NOneChi}
\end{equation}
With \(\psi_1:=\chi_1^\sharp\), the normalized matrix candidate is
\[
 \Phi_1=
 \begin{pmatrix}
 \psi_1(\lambda)&\psi_1(-\lambda)\\
 \chi_1(\lambda)&-\chi_1(-\lambda)
 \end{pmatrix}.
\]
At lower spectral infinity, \(\mu/K\to-1\), and hence
\(\chi_1\to1\). At upper spectral infinity, \(\mu/K\to1\), and the
Ernst candidate is
\begin{equation}
 f=\frac{q_1+q_2}{q_1-q_2}
 =\frac{q_1-\ell}{q_1+\ell}
 =\frac{\alpha_1r_1-\alpha_2r_2-\ell}
        {\alpha_1r_1-\alpha_2r_2+\ell}.
 \label{eq:NOneErnst}
\end{equation}

\medskip
\noindent\textit{Zeros of the denominator.}
Introduce
\[
 x:=\frac{r_1+r_2}{\ell},\qquad
 y:=\frac{r_1-r_2}{\ell}.
\]
The triangle inequalities give
\[
 x\geq1,\qquad -1\leq y\leq1,
\]
while
\[
 r_1=\frac{\ell}{2}(x+y),\qquad
 r_2=\frac{\ell}{2}(x-y).
\]
It follows that \(\Delta_1=0\) if and only if
\begin{equation}
 (\alpha_1-\alpha_2)x+(\alpha_1+\alpha_2)y=-2
 \label{eq:NOneDenominatorZeroCondition}
\end{equation}
for some
\[
 (x,y)\in[1,\infty)\times[-1,1].
\]
If \(\Delta_1=0\), then \(q_1=-\ell\), whereas the numerator of \(f\)
equals \(q_1-\ell=-2\ell\ne0\). Thus every zero of \(\Delta_1\) is a
genuine pole of \(f\); no removable \(0/0\) singularity occurs.

\medskip
\noindent\textit{Boundary values of \(f\).}
The three open boundary intervals are
\[
 \mathcal A^+=(K_1,\infty),\qquad
 \mathcal H^1=(K_2,K_1),\qquad
 \mathcal A^-=(-\infty,K_2).
\]
The corresponding pairs \(u_i^I=\alpha_i\lambda_i^I\) are
\[
 (u_1^+,u_2^+)=(-\alpha_1,-\alpha_2),
\]
\[
 (u_1^{\mathcal H^1},u_2^{\mathcal H^1})
 =(\alpha_1,-\alpha_2),
\]
and
\[
 (u_1^-,u_2^-)=(\alpha_1,\alpha_2).
\]
For \(\zeta\in I\), one has
\[
 \alpha_i r_i=u_i^I(K_i-\zeta).
\]
Hence
\[
 f^I(\zeta)=\frac{\widetilde P_1^I(\zeta)}
                  {\widetilde Q_1^I(\zeta)},
\]
where
\begin{align}
 \widetilde P_1^I(\zeta)
 &=-H_1^I\zeta+u_1^IK_1-u_2^IK_2-\ell,\\
 \widetilde Q_1^I(\zeta)
 &=-H_1^I\zeta+u_1^IK_1-u_2^IK_2+\ell
 =:\Delta_1^I(\zeta),
\end{align}
and
\[
 H_1^I=[c^I,d]=u_1^I-u_2^I.
\]
Explicitly,
\begin{equation}
 H_1^+=\alpha_2-\alpha_1,\qquad
 H_1^{\mathcal H^1}=\alpha_1+\alpha_2,\qquad
 H_1^-=\alpha_1-\alpha_2=-H_1^+.
 \label{eq:NOneHs}
\end{equation}

\medskip
\noindent\textit{Boundary values of \(e^{2U}\) and \(e^{2k}\).}
The general determinant formula reduces to
\begin{equation}
 e^{2U^I}(\zeta)
 =\frac{(K_1-\zeta)(K_2-\zeta)|H_1^I|^2}
        {|\Delta_1^I(\zeta)|^2}.
 \label{eq:NOneBoundaryU}
\end{equation}
Thus, whenever the displayed numerator and denominator are nonzero,
\[
 e^{2U^+}>0,\qquad
 e^{2U^{\mathcal H^1}}<0,\qquad
 e^{2U^-}>0.
\]
If
\[
 H_1^{\mathcal H^1}=\alpha_1+\alpha_2=0,
\]
then
\[
 e^{2U^{\mathcal H^1}}=0.
\]

Assume
\[
 H_1^+\ne0,\qquad\text{namely}\qquad
 \alpha_1-\alpha_2\ne0,
\]
and normalize \(e^{2k^+}=1\). Then
\[
 e^{2k^I}=\left(\frac{H_1^I}{H_1^+}\right)^2.
\]
In particular, \(e^{2k^-}=1\). If
\[
 H_1^{\mathcal H^1}=\alpha_1+\alpha_2\ne0,
\]
then
\[
 e^{2k^{\mathcal H^1}}
 =\left(\frac{\alpha_1+\alpha_2}{\alpha_2-\alpha_1}\right)^2<0.
\]
If \(H_1^{\mathcal H^1}=0\), then \(e^{2k^{\mathcal H^1}}=0\).

\medskip
\noindent\textit{Boundary values of \(a\).}
For \(N=1\),
\[
 D_1^I=[v^IK,v^I]
 =\ell(1+u_1^I)(1+u_2^I),
\qquad
 v^I=(1+u_1^I,1+u_2^I).
\]
Whenever \(H_1^I\ne0\),
\[
 a^I+C=-i\frac{D_1^I}{H_1^I}.
\]
We fix \(C\) by the upper-axis normalization \(a^+=0\). Since
\[
 D_1^+=\ell(1-\alpha_1)(1-\alpha_2),
\]
this gives
\[
 C=-i\ell\frac{(1-\alpha_1)(1-\alpha_2)}
                    {\alpha_2-\alpha_1}.
\]
Consequently,
\[
 a^I=-i\left(\frac{D_1^I}{H_1^I}
             -\frac{D_1^+}{H_1^+}\right).
\]
For the lower axis this simplifies to
\begin{equation}
 a^-=-2i\ell\frac{1+\alpha_1\alpha_2}
                         {\alpha_1-\alpha_2}.
 \label{eq:NOneAMinus}
\end{equation}
If \(H_1^{\mathcal H^1}=\alpha_1+\alpha_2\ne0\), then
\[
 a^{\mathcal H^1}
 =-2i\ell\frac{\alpha_1(1-\alpha_2)^2}
 {(\alpha_1-\alpha_2)(\alpha_1+\alpha_2)}.
\]

\medskip
\noindent\textit{The exceptional case \(H_1^{\mathcal H^1}=0\).}
Suppose
\[
 \alpha_1+\alpha_2=0,\qquad
 \alpha_1=\alpha,\qquad
 \alpha_2=-\alpha,\qquad |\alpha|=1.
\]
For \(\zeta\in\mathcal H^1=(K_2,K_1)\), one has
\[
 r_1+r_2\longrightarrow\ell,\qquad
 q_1\longrightarrow\alpha\ell,\qquad
 \Delta_1(\zeta)=q_1+\ell\longrightarrow\ell(1+\alpha)
\]
as \(\rho\downarrow0\). If \(\alpha\ne-1\), the Ernst potential has
the finite boundary value
\[
 f^{\mathcal H^1}=\frac{\alpha-1}{\alpha+1}.
\]

The cancelled formula involving \(D_1^{\mathcal H^1}/H_1^{\mathcal H^1}\)
cannot be used. Returning to the un-cancelled reconstruction formula and
imposing \(a^+=0\), one obtains
\begin{align}
 a(\rho,\zeta)
 ={}&-i\ell(\alpha-\alpha^{-1})
 \frac{(\zeta-K_2)r_1-(K_1-\zeta)r_2}
 {(r_1+r_2)^2-\ell^2}
 \notag\\
 &+\frac{i\ell}{2}(\alpha-\alpha^{-1}).
 \label{eq:NOneExceptionalAInterior}
\end{align}
Both the numerator and denominator in the quotient vanish to order
\(\rho^2\). Their leading terms give
\[
 \frac{(\zeta-K_2)r_1-(K_1-\zeta)r_2}
 {(r_1+r_2)^2-\ell^2}
 =\frac{2\zeta-K_1-K_2}{2\ell}+O(\rho^2).
\]
Consequently,
\begin{equation}
 a^{\mathcal H^1}(\zeta)
 =i(\alpha-\alpha^{-1})(K_1-\zeta),
 \qquad
 a^-=i\ell(\alpha-\alpha^{-1}).
 \label{eq:NOneExceptionalBoundaryA}
\end{equation}
For generic \(\alpha\), the first expression is not constant and
\(a^-\ne a^+=0\). It is constant precisely when \(\alpha=\pm1\).

For \(\alpha=1\),
\[
 \alpha_1=1,\qquad \alpha_2=-1,\qquad
 f=\frac{r_1+r_2-\ell}{r_1+r_2+\ell},
\]
and
\[
 a^{\mathcal H^1}=a^-=a^+=0.
\]
This is the usual Schwarzschild Ernst potential with mass \(m=\ell/2\).

For \(\alpha=-1\),
\[
 f=\frac{r_1+r_2+\ell}{r_1+r_2-\ell},
\]
the reciprocal of the usual Schwarzschild Ernst potential. Since
\(r_1+r_2=\ell\) on the whole interval \(\mathcal H^1\), its denominator
vanishes there whereas its numerator does not. Thus \(f\) has a genuine
pole along the entire interval.

\medskip
\noindent\textit{The equilibrium condition.}
With \(a^+=0\) and \(e^{2k^+}=1\), the lower-axis elementary-flatness
condition is automatic because \(H_1^-=-H_1^+\). The only nontrivial
equilibrium equation is
\[
 a^-=a^+=0.
\]
Under \(H_1^+\ne0\), equation \eqref{eq:NOneAMinus} shows that this is
equivalent to
\begin{equation}
 \alpha_1\alpha_2=-1.
 \label{eq:NOneEquilibrium}
\end{equation}
If one additionally requires the middle interval to carry a constant
value of \(a\), then the exceptional calculation above leaves
\(\alpha_1=\pm1\). Of these choices, only
\[
 \alpha_1=1,\qquad \alpha_2=-1
\]
has a finite Ernst potential on \(\mathcal H^1\); it gives the usual
positive-mass Schwarzschild solution.

\medskip
\noindent\textit{The degeneration \(H_1^+=0\).}
Suppose
\[
 H_1^+=0,\qquad
 \alpha_1=\alpha_2=\alpha,\qquad |\alpha|=1.
\]
The determinant formula becomes
\begin{equation}
 f(\rho,\zeta)
 =\frac{\alpha(r_1-r_2)-\ell}{\alpha(r_1-r_2)+\ell}.
 \label{eq:NOneHplusZeroF}
\end{equation}
On \(\mathcal A^+\), one has \(r_1-r_2=-\ell\). Thus, provided
\(\alpha\ne1\),
\begin{equation}
 f^+(\zeta)=-\frac{1+\alpha}{1-\alpha},
 \label{eq:NOneHplusZeroAxisF}
\end{equation}
which is independent of \(\zeta\). Nevertheless, the corresponding
Ernst potential is not locally constant.

More importantly, it does not satisfy the usual asymptotic-flatness
normalization. Introduce polar coordinates in the orbit half-plane by
\[
 \rho=R\sin\theta,\qquad \zeta=R\cos\theta.
\]
As \(R\to\infty\),
\[
 r_1-r_2=-\ell\cos\theta+O(R^{-1}),
\]
and hence
\begin{equation}
 f(\rho,\zeta)
 \longrightarrow
 -\frac{1+\alpha\cos\theta}{1-\alpha\cos\theta}.
 \label{eq:NOneHplusZeroInfinity}
\end{equation}
The limiting value generally depends on the direction \(\theta\) and
cannot be identically equal to one. This example shows that constant
data on an open part of the axis do not determine an Ernst potential
among arbitrary local solutions. Any corresponding uniqueness theorem
must impose additional global hypotheses, such as asymptotic flatness,
suitable regularity, and the analytic or spectral assumptions required
by the relevant inverse-scattering uniqueness result.

\medskip
\noindent\textit{Simplification by equilibrium conditions.}
We now impose the equilibrium condition
\[
 \alpha_1\alpha_2=-1.
\]
Write
\[
 \alpha_1=\alpha=e^{i\theta},
 \qquad
 \alpha_2=-\alpha^{-1}=-e^{-i\theta},
 \qquad
 \theta\in\mathbb R.
\]
The Ernst potential simplifies to
\begin{equation}
 f(\rho,\zeta)
 =
 \frac{\alpha r_1+\alpha^{-1}r_2-\ell}
      {\alpha r_1+\alpha^{-1}r_2+\ell}.
 \label{eq:NOneEquilibriumFAlpha}
\end{equation}
Equivalently,
\begin{equation}
 f(\rho,\zeta)
 =
 \frac{e^{i\theta}r_1+e^{-i\theta}r_2-\ell}
      {e^{i\theta}r_1+e^{-i\theta}r_2+\ell}.
 \label{eq:NOneEquilibriumFTheta}
\end{equation}

The denominator-zero condition
\eqref{eq:NOneDenominatorZeroCondition} becomes
\[
 (\alpha+\alpha^{-1})x+(\alpha-\alpha^{-1})y=-2.
\]
Since
\[
 \alpha+\alpha^{-1}=2\cos\theta,
 \qquad
 \alpha-\alpha^{-1}=2i\sin\theta,
\]
this is equivalent to the two real equations
\begin{equation}
 \cos\theta\,x=-1,
 \qquad
 \sin\theta\,y=0,
 \label{eq:NOneEquilibriumDenominatorZeros}
\end{equation}
where
\[
 x\geq1,
 \qquad
 -1\leq y\leq1.
\]

If \(\cos\theta>0\), the first equation in
\eqref{eq:NOneEquilibriumDenominatorZeros} has no solution, and the
determinant denominator is nowhere zero. If \(\cos\theta=0\), it again
has no solution, but
\[
 H_1^+=-2\cos\theta=0,
\]
so this is the non-asymptotically-flat degeneration discussed above.

If
\[
 -1<\cos\theta<0,
\]
the denominator has exactly one zero, determined by
\[
 x=-\frac1{\cos\theta},
 \qquad
 y=0.
\]
In the original orbit-space coordinates, this point is
\[
 \zeta=\frac{K_1+K_2}{2},
 \qquad
 \rho=\frac{\ell}{2}|\tan\theta|.
\]
It is a genuine pole of \(f\).

Finally, if
\[
 \alpha=-1,
 \qquad\text{equivalently}\qquad
 \theta=\pi\pmod{2\pi},
\]
the zero condition reduces to \(x=1\). This corresponds to the entire
closed interval
\[
 \rho=0,
 \qquad
 K_2\leq\zeta\leq K_1,
\]
along which \(f\) has a genuine pole.

Consequently, the family satisfying
\[
 \operatorname{Re}\alpha=\cos\theta>0
\]
is free of zeros of the determinant denominator. This is precisely the
Kerr family in the present parametrization. The special value
\[
 \alpha=1,
 \qquad\text{equivalently}\qquad
 \theta=0,
\]
gives the Schwarzschild solution.

\medskip
\noindent\textit{Candidate thermodynamic quantities.}
We continue to impose the equilibrium condition and use the
parametrization
\[
 \alpha_1=\alpha=e^{i\theta},
 \qquad
 \alpha_2=-\alpha^{-1}=-e^{-i\theta}.
\]
The determinants entering Proposition~\ref{prop:candidateHorizonQuantities}
reduce to
\[
 H_1^+=-\alpha-\alpha^{-1}=-2\cos\theta,
 \qquad
 H_1^{\mathcal H^1}=\alpha-\alpha^{-1}=2i\sin\theta,
\]
and
\[
 \Lambda_1=\ell,
 \qquad
 E_1=1+\alpha^{-1},
 \qquad
 E_2=1+\alpha,
 \qquad
 H_0^{(1)}=1.
\]
Assume first that
\[
 \cos\theta\ne0,
 \qquad
 \sin\theta\ne0,
\]
so that the definitions containing \(\Omega^{-1}\) can be used
directly.  Proposition~\ref{prop:candidateHorizonQuantities} gives
\[
 \kappa+i\Omega
 =\frac{H_1^+(H_1^+-H_1^{\mathcal H^1})}
 {2\alpha\ell E_1^2}
 =\frac{\alpha\cos\theta}{\ell(1+\cos\theta)}.
\]
Consequently,
\begin{equation}
 \kappa=\frac{\cos^2\theta}{\ell(1+\cos\theta)},
 \qquad
 \Omega=\frac{\sin\theta\cos\theta}
 {\ell(1+\cos\theta)}.
 \label{eq:NOneKappaOmega}
\end{equation}

The endpoint values can be read directly from
\eqref{eq:NOneEquilibriumFAlpha}.  At \(\zeta=K_1\), one has
\(r_1=0\) and \(r_2=\ell\), whereas at \(\zeta=K_2\), one has
\(r_1=\ell\) and \(r_2=0\).  Hence
\[
 f_1=\frac{1-\alpha}{1+\alpha}
 =-i\tan\frac\theta2,
 \qquad
 f_2=\frac{\alpha-1}{\alpha+1}
 =i\tan\frac\theta2,
\]
and therefore
\[
 f_1-f_2=-2i\tan\frac\theta2.
\]
The definitions in Section~6 now yield
\begin{align}
 \mathfrak M
 &=\frac{i}{4\Omega}(f_1-f_2)
 =\frac{\ell}{2\cos\theta},
 \label{eq:NOneCandidateMass}\\
 A
 &=\frac{2\pi\ell}{\kappa}
 =2\pi\ell^2\frac{1+\cos\theta}{\cos^2\theta},
 \label{eq:NOneCandidateArea}\\
 J
 &=\frac1{2\Omega}\left(\mathfrak M-\frac\ell2\right)
 =\frac{\ell^2\sin\theta}{4\cos^2\theta},
 \label{eq:NOneCandidateAngularMomentum}\\
 p
 &=\frac{8\pi J}{A}
 =\frac{\sin\theta}{1+\cos\theta}
 =\tan\frac\theta2.
 \label{eq:NOneCandidateP}
\end{align}
The boundary identity \eqref{eq:aJumpThroughOmega}, together with
\(a^+=0\), also gives
\begin{equation}
 a^{\mathcal H^1}
 =-\frac1\Omega
 =-\ell\frac{1+\cos\theta}{\sin\theta\cos\theta}.
 \label{eq:NOneCandidateHorizonA}
\end{equation}
The last formula applies only when \(\Omega\ne0\).

Since
\[
 \ell=2\mathfrak M\cos\theta,
\]
the same quantities can equivalently be written as
\begin{equation}
 \kappa=\frac{\cos\theta}
 {2\mathfrak M(1+\cos\theta)},
 \qquad
 \Omega=\frac{\sin\theta}
 {2\mathfrak M(1+\cos\theta)},
 \label{eq:NOneKerrKappaOmega}
\end{equation}
and
\begin{equation}
 A=8\pi\mathfrak M^2(1+\cos\theta),
 \qquad
 J=\mathfrak M^2\sin\theta.
 \label{eq:NOneKerrAreaJ}
\end{equation}
These are the Kerr relations in the present parametrization.  In
particular,
\[
 \mathfrak M>0
 \quad\Longleftrightarrow\quad
 \cos\theta>0,
\]
while
\[
 \frac{8\pi|J|}{A}
 =\left|\tan\frac\theta2\right|.
\]
It follows that
\begin{equation}
 \mathfrak M>0
 \quad\Longleftrightarrow\quad
 8\pi|J|<A
 \quad\Longleftrightarrow\quad
 \operatorname{Re}\alpha=\cos\theta>0.
 \label{eq:NOneEquivalentAdmissibilityConditions}
\end{equation}
This is exactly the parameter range in which the determinant
denominator was shown above to have no zeros.

At the special value
\[
 \theta=0,
 \qquad
 \alpha=1,
\]
one has \(\Omega=0\), so the definitions containing
\(\Omega^{-1}\) cannot be applied directly.  Their continuous limits,
together with the formulas not involving \(\Omega^{-1}\), give
\begin{equation}
 \kappa=\frac1{2\ell},
 \qquad
 \Omega=0,
 \qquad
 A=4\pi\ell^2,
 \qquad
 \mathfrak M=\frac\ell2,
 \qquad
 J=0,
 \qquad
 p=0.
 \label{eq:NOneSchwarzschildThermodynamics}
\end{equation}
These are precisely the Schwarzschild values for mass
\(m=\ell/2\).

\subsection{The case \texorpdfstring{\(N=2\)}{N=2}}

Let
\[
 K_1>K_2>K_3>K_4,
 \qquad
 |\alpha_i|=1,
 \qquad
 i\in\{1,2,3,4\}.
\]
Assume that the boundary denominator has no zeros and that
\[
 H_2^+\ne0.
\]
We use the normalizations
\[
 e^{2k^+}=1,
 \qquad
 a^+=0.
\]
By Section~5, the equilibrium and branch-selection conditions are
equivalent to
\[
 D_2^{\mathcal A^1}=-D_2^+,
 \qquad
 D_2^-=D_2^+,
\]
and
\[
 H_2^{\mathcal A^1}=-H_2^+,
 \qquad
 H_2^-=H_2^+.
\]
In particular, \(H_2\) is nonzero on every open axis interval.  For
\(N=2\), these conditions reduce to
\begin{equation}
 \alpha_1\alpha_2+\alpha_3\alpha_4=0,
 \label{eq:NTwoFirstEquilibrium}
\end{equation}
and
\begin{equation}
 \frac{(1-\alpha_4)^2}{\alpha_4}w^2
 =\frac{(1-\alpha_3)^2}{\alpha_3},
 \qquad
 \frac{(1+\alpha_2)^2}{\alpha_2}{w'}^2
 =\frac{(1+\alpha_1)^2}{\alpha_1},
 \label{eq:NTwoRemainingEquilibrium}
\end{equation}
where
\begin{equation}
 w:=\sqrt{\frac{K_{14}K_{24}}{K_{13}K_{23}}},
 \qquad
 w':=\sqrt{\frac{K_{23}K_{24}}{K_{13}K_{14}}},
 \qquad
 K_{ij}:=K_i-K_j.
 \label{eq:NTwoWWPrime}
\end{equation}
For strictly extended horizons,
\[
 w>1,
 \qquad
 0<w'<1.
\]

Following the explicit parametrization of the double-Kerr equilibrium
conditions obtained by Manko and Ruiz~\cite{MankoRuiz2001}, introduce
\[
 \alpha=e^{i\phi},
 \qquad
 \varepsilon\in\{+1,-1\}.
\]
Then
\begin{equation}
 \alpha_1
 =\frac{w'\alpha^2+i\varepsilon\alpha}
 {w'-i\varepsilon\alpha},
 \qquad
 \alpha_2
 =\frac{\alpha^2+iw'\varepsilon\alpha}
 {1-iw'\varepsilon\alpha},
 \label{eq:NTwoAlpha12}
\end{equation}
and
\begin{equation}
 \alpha_3
 =\frac{w\alpha^2-\alpha}{w-\alpha},
 \qquad
 \alpha_4
 =\frac{\alpha^2-w\alpha}{1-w\alpha},
 \label{eq:NTwoAlpha34}
\end{equation}
where
\[
 \alpha=\sqrt{-\alpha_1\alpha_2}
 =\sqrt{\alpha_3\alpha_4},
 \qquad
 |\alpha|=1.
\]
Thus, apart from an overall translation and scaling of the \(K_i\),
the family satisfying the equilibrium conditions is parametrized by
\[
 w,\qquad w',\qquad\phi,\qquad\varepsilon.
\]
The preceding equilibrium analysis was developed in
\cite{MankoRuizSanabria2000,MankoRuiz2001}.

Assume in addition that
\[
 H_2^{\mathcal H^1}\ne0,
 \qquad
 H_2^{\mathcal H^2}\ne0,
\]
so that the candidate quantities containing \(\Omega_j^{-1}\) are
directly defined.  Put
\[
 p_j:=\frac{8\pi J_j}{A_j},
 \qquad
 j\in\{1,2\},
\]
and
\[
 \mathit\Phi:=\cos\phi+\varepsilon\sin\phi.
\]
The formulas of Section~6 become
\begin{equation}
 p_1=\varepsilon\frac{1+\mathit\Phi w'}
 {w'(\mathit\Phi+w')},
 \qquad
 p_2=\varepsilon\frac{w(w-\mathit\Phi)}
 {1-w\mathit\Phi}.
 \label{eq:NTwoPValues}
\end{equation}
The two strict subextremal inequalities
\[
 |p_1|<1,
 \qquad
 |p_2|<1
\]
would require, respectively,
\begin{equation}
 {w'}^2+2\mathit\Phi w'+1<0,
 \qquad
 w^2-2\mathit\Phi w+1<0.
 \label{eq:NTwoSubextremalConditions}
\end{equation}
The first inequality implies \(\mathit\Phi w'<0\), whereas the second
implies \(\mathit\Phi w>0\).  Since \(w,w'>0\), these requirements are
incompatible.  Consequently, within the non-degenerate family satisfying
the equilibrium conditions, the two inequalities \(|p_1|<1\) and \(|p_2|<1\) cannot hold
simultaneously.  This is the algebraic contradiction used by Hennig
and Neugebauer in the non-existence proof for two subextremal
constituents.  The exceptional cases
\[
 H_2^{\mathcal H^j}=0
\]
are not covered by the formulas containing \(\Omega_j^{-1}\) and
require a separate limiting analysis.

\subsection{The case \texorpdfstring{\(N=3\)}{N=3}}

Let
\[
 K_1>K_2>K_3>K_4>K_5>K_6,
 \qquad |\alpha_i|=1,\qquad i\in\{1,\ldots,6\}.
\]
Write
\[
 \mathcal A^0=\mathcal A^+,\qquad
 \mathcal A^1,\qquad \mathcal A^2,\qquad
 \mathcal A^3=\mathcal A^-.
\]
Assume
\[
 \Delta_3(0,\zeta)\ne0\quad(\zeta\in\mathbb R),
 \qquad H_3^+\ne0,
\]
and normalize \(a^+=0\) and \(e^{2k^+}=1\).  For an open boundary
interval \(I\), recall
\[
 c^I=(u_1^I,\ldots,u_6^I),\qquad
 d=(1,\ldots,1),\qquad v^I=c^I+d,
\]
and put \(K=\operatorname{diag}(K_1,\ldots,K_6)\).  Then
\begin{equation}
 H_3^I=[c^IK^2,dK^2,c^IK,dK,c^I,d],
 \qquad
 D_3^I=[v^IK^3,v^IK^2,v^IK,v^I,dK,d].
 \label{eq:NThreeHDDefinitions}
\end{equation}
The axis sign patterns are
\[
\begin{array}{c|rrrrrr}
 I&u_1^I&u_2^I&u_3^I&u_4^I&u_5^I&u_6^I\\ \hline
 \mathcal A^0&-\alpha_1&-\alpha_2&-\alpha_3&-\alpha_4&-\alpha_5&-\alpha_6\\
 \mathcal A^1& \alpha_1& \alpha_2&-\alpha_3&-\alpha_4&-\alpha_5&-\alpha_6\\
 \mathcal A^2& \alpha_1& \alpha_2& \alpha_3& \alpha_4&-\alpha_5&-\alpha_6\\
 \mathcal A^3& \alpha_1& \alpha_2& \alpha_3& \alpha_4& \alpha_5& \alpha_6
\end{array}
\]
Theorem~\ref{thm:branchSelection} gives the five equilibrium equations
\begin{equation}
\begin{aligned}
 D_3^{\mathcal A^0}+D_3^{\mathcal A^1}&=0,&
 D_3^{\mathcal A^1}+D_3^{\mathcal A^2}&=0,&
 D_3^{\mathcal A^2}+D_3^{\mathcal A^3}&=0,\\
 H_3^{\mathcal A^0}+H_3^{\mathcal A^1}&=0,&
 H_3^{\mathcal A^1}+H_3^{\mathcal A^2}&=0.&
\end{aligned}
\label{eq:NThreeFiveDeterminantEquations}
\end{equation}
The missing lower-axis \(H_3\)-equation follows automatically from
\(H_3^-=(-1)^3H_3^+\).

We now expand the five equations.  For
\(S\subset\{1,\ldots,6\}\), define
\[
 V_S:=\prod_{\substack{i,j\in S\\i<j}}K_{ij},
 \qquad K_{ij}:=K_i-K_j.
\]
For a three-element set \(\{i,j,k\}\), put
\begin{equation}
 W_{ijk}:=V_{\{i,j,k\}}
 V_{\{1,\ldots,6\}\setminus\{i,j,k\}},
 \qquad A_{ijk}:=\alpha_i\alpha_j\alpha_k.
 \label{eq:NThreeWAdefinitions}
\end{equation}
For a four-element set \(S=\{i,j,k,l\}\), let
\(S^c=\{m,n\}\), where \(m<n\), and define
\begin{equation}
 Z_{ijkl}:=(-1)^{i+j+k+l}V_{\{i,j,k,l\}}K_{mn}.
 \label{eq:NThreeZDefinition}
\end{equation}
Thus every \(W_{ijk}\) and \(Z_{ijkl}\) is an explicit product of
the differences \(K_{pq}\).

The two \(H_3\)-equations are
\begin{equation}
\begin{aligned}
0={}&W_{123}(A_{123}-A_{456})
-W_{124}(A_{124}-A_{356})\\
&+W_{125}(A_{125}-A_{346})
-W_{126}(A_{126}-A_{345}),
\end{aligned}
\label{eq:NThreeExpandedHOne}
\end{equation}
and
\begin{equation}
\begin{aligned}
0={}&-W_{134}(A_{134}-A_{256})
+W_{234}(A_{234}-A_{156})\\
&+W_{345}(A_{345}-A_{126})
-W_{346}(A_{346}-A_{125}).
\end{aligned}
\label{eq:NThreeExpandedHTwo}
\end{equation}

The equation \(D_3^{\mathcal A^0}+D_3^{\mathcal A^1}=0\) is
\begin{align}
0={}&Z_{3456}(1-\alpha_3)(1-\alpha_4)(1-\alpha_5)(1-\alpha_6)\notag\\
&+Z_{1234}(1+\alpha_1\alpha_2)(1-\alpha_3)(1-\alpha_4)\notag\\
&+Z_{1235}(1+\alpha_1\alpha_2)(1-\alpha_3)(1-\alpha_5)\notag\\
&+Z_{1236}(1+\alpha_1\alpha_2)(1-\alpha_3)(1-\alpha_6)\notag\\
&+Z_{1245}(1+\alpha_1\alpha_2)(1-\alpha_4)(1-\alpha_5)\notag\\
&+Z_{1246}(1+\alpha_1\alpha_2)(1-\alpha_4)(1-\alpha_6)\notag\\
&+Z_{1256}(1+\alpha_1\alpha_2)(1-\alpha_5)(1-\alpha_6)\notag\\
&+Z_{1345}(1-\alpha_3)(1-\alpha_4)(1-\alpha_5)\notag\\
&+Z_{1346}(1-\alpha_3)(1-\alpha_4)(1-\alpha_6)\notag\\
&+Z_{1356}(1-\alpha_3)(1-\alpha_5)(1-\alpha_6)\notag\\
&+Z_{1456}(1-\alpha_4)(1-\alpha_5)(1-\alpha_6)\notag\\
&+Z_{2345}(1-\alpha_3)(1-\alpha_4)(1-\alpha_5)\notag\\
&+Z_{2346}(1-\alpha_3)(1-\alpha_4)(1-\alpha_6)\notag\\
&+Z_{2356}(1-\alpha_3)(1-\alpha_5)(1-\alpha_6)\notag\\
&+Z_{2456}(1-\alpha_4)(1-\alpha_5)(1-\alpha_6).
\label{eq:NThreeExpandedDOne}
\end{align}

The equation \(D_3^{\mathcal A^1}+D_3^{\mathcal A^2}=0\) is
\begin{align}
0={}&Z_{1256}(1+\alpha_1)(1+\alpha_2)(1-\alpha_5)(1-\alpha_6)\notag\\
&+Z_{1234}(1+\alpha_3\alpha_4)(1+\alpha_1)(1+\alpha_2)\notag\\
&+Z_{1345}(1+\alpha_3\alpha_4)(1+\alpha_1)(1-\alpha_5)\notag\\
&+Z_{1346}(1+\alpha_3\alpha_4)(1+\alpha_1)(1-\alpha_6)\notag\\
&+Z_{2345}(1+\alpha_3\alpha_4)(1+\alpha_2)(1-\alpha_5)\notag\\
&+Z_{2346}(1+\alpha_3\alpha_4)(1+\alpha_2)(1-\alpha_6)\notag\\
&+Z_{3456}(1+\alpha_3\alpha_4)(1-\alpha_5)(1-\alpha_6)\notag\\
&+Z_{1235}(1+\alpha_1)(1+\alpha_2)(1-\alpha_5)\notag\\
&+Z_{1236}(1+\alpha_1)(1+\alpha_2)(1-\alpha_6)\notag\\
&+Z_{1356}(1+\alpha_1)(1-\alpha_5)(1-\alpha_6)\notag\\
&+Z_{2356}(1+\alpha_2)(1-\alpha_5)(1-\alpha_6)\notag\\
&+Z_{1245}(1+\alpha_1)(1+\alpha_2)(1-\alpha_5)\notag\\
&+Z_{1246}(1+\alpha_1)(1+\alpha_2)(1-\alpha_6)\notag\\
&+Z_{1456}(1+\alpha_1)(1-\alpha_5)(1-\alpha_6)\notag\\
&+Z_{2456}(1+\alpha_2)(1-\alpha_5)(1-\alpha_6).
\label{eq:NThreeExpandedDTwo}
\end{align}

\clearpage
The equation \(D_3^{\mathcal A^2}+D_3^{\mathcal A^3}=0\) is
\begin{align}
0={}&Z_{1234}(1+\alpha_1)(1+\alpha_2)(1+\alpha_3)(1+\alpha_4)\notag\\
&+Z_{1256}(1+\alpha_5\alpha_6)(1+\alpha_1)(1+\alpha_2)\notag\\
&+Z_{1356}(1+\alpha_5\alpha_6)(1+\alpha_1)(1+\alpha_3)\notag\\
&+Z_{1456}(1+\alpha_5\alpha_6)(1+\alpha_1)(1+\alpha_4)\notag\\
&+Z_{2356}(1+\alpha_5\alpha_6)(1+\alpha_2)(1+\alpha_3)\notag\\
&+Z_{2456}(1+\alpha_5\alpha_6)(1+\alpha_2)(1+\alpha_4)\notag\\
&+Z_{3456}(1+\alpha_5\alpha_6)(1+\alpha_3)(1+\alpha_4)\notag\\
&+Z_{1235}(1+\alpha_1)(1+\alpha_2)(1+\alpha_3)\notag\\
&+Z_{1245}(1+\alpha_1)(1+\alpha_2)(1+\alpha_4)\notag\\
&+Z_{1345}(1+\alpha_1)(1+\alpha_3)(1+\alpha_4)\notag\\
&+Z_{2345}(1+\alpha_2)(1+\alpha_3)(1+\alpha_4)\notag\\
&+Z_{1236}(1+\alpha_1)(1+\alpha_2)(1+\alpha_3)\notag\\
&+Z_{1246}(1+\alpha_1)(1+\alpha_2)(1+\alpha_4)\notag\\
&+Z_{1346}(1+\alpha_1)(1+\alpha_3)(1+\alpha_4)\notag\\
&+Z_{2346}(1+\alpha_2)(1+\alpha_3)(1+\alpha_4).
\label{eq:NThreeExpandedDThree}
\end{align}
Equations~\eqref{eq:NThreeExpandedHOne},
\eqref{eq:NThreeExpandedHTwo}, and
\eqref{eq:NThreeExpandedDOne}--\eqref{eq:NThreeExpandedDThree} are the
complete branch-independent equilibrium equations in terms of
\(K_i\) and \(\alpha_i\).  By the conjugation identities of Section~4,
each apparently complex equation represents one real condition.  No
assumption on \(H_3^{\mathcal H^j}\) is needed at this stage.

\begin{remark}
The same subset expansion gives, in principle, completely explicit
equilibrium equations for every \(N\).  Their number grows linearly
with \(N\), but the number of terms in their expansion in \(K_i\) and
\(\alpha_i\) grows combinatorially; consequently, the unreduced
formulas rapidly become unwieldy.

The more substantial problem is to find coordinates adapted to these
equations.  For \(N=1\), after the inessential length scale and the
normalizations have been separated, the equilibrium condition leaves
a single phase parameter.  For \(N=2\), the Manko--Ruiz
parametrization provides a particularly effective set of parameters
in which the equilibrium conditions can be solved explicitly.  To the
authors' knowledge, no comparably useful parametrization is presently
available for \(N=3\), or for general \(N\).  The five equations above
therefore give an explicit algebraic starting point for the search for
such a parametrization in the first open case.
\end{remark}

\section{Some conjectures}

We conclude with three algebraic conjectures suggested by the preceding
determinant formulas.  They concern only the functions and candidate
quantities defined in this paper.  In particular, their formulation does
not presuppose that the determinant solution extends to a globally regular
spacetime or that the intervals \(\mathcal H^j\) are regular Killing
horizons.

The first conjecture concerns the relation between the asymptotic mass and
the candidate horizon masses.  Define the candidate ADM mass by the upper
axis expansion
\begin{equation}
 f^+(\zeta)
 =1-\frac{2M_{\mathrm{ADM}}}{\zeta}
  +O(\zeta^{-2})
 \qquad (\zeta\to+\infty).
 \label{eq:CandidateADMDefinition}
\end{equation}

The coefficient in this expansion has a direct determinant expression.
Recall that on the upper axis
\[
 c^+=(-\alpha_1,\ldots,-\alpha_{2N}),
 \qquad
 d=(1,\ldots,1),
\]
and define
\begin{equation}
 \mathcal G_N^+
 :=
 [dK^N,dK^{N-1},c^+K^{N-2},dK^{N-2},
 \ldots,c^+K,dK,c^+,d].
 \label{eq:ADMnumeratorDeterminant}
\end{equation}
Thus \(\mathcal G_N^+\) is obtained from \(H_N^+\) by replacing its
first row \(c^+K^{N-1}\) by \(dK^N\).  Whenever \(H_N^+\ne0\),
\begin{equation}
 M_{\mathrm{ADM}}
 =-\frac{\mathcal G_N^+}{H_N^+}.
 \label{eq:CandidateADMDeterminantFormula}
\end{equation}
Indeed, Proposition~\ref{prop:boundaryErnst} gives
\[
 f^+(\zeta)
 =\frac{q_1^+(\zeta)+q_2^+(\zeta)}
        {q_1^+(\zeta)-q_2^+(\zeta)},
\]
and the common leading coefficient of its numerator and denominator is
\begin{equation}
 [\zeta^N]\bigl(q_1^++q_2^+\bigr)
 =
 [\zeta^N]\bigl(q_1^+-q_2^+\bigr)
 =(-1)^NH_N^+.
 \label{eq:ADMcommonLeadingCoefficient}
\end{equation}
Moreover, \(q_2^+=-B_N^+\), where
\[
 B_N^+
 =[dK^N,dK^{N-1},c_1^+K^{N-2},dK^{N-2},
 \ldots,c_1^+,d]
\]
and \(c_1^+=c^+K-\zeta c^+\).  Selecting the term
\(-\zeta c^+K^m\) from each of the \(N-1\) rows containing
\(c_1^+K^m\) yields
\begin{equation}
 [\zeta^{N-1}]q_2^+=(-1)^N\mathcal G_N^+.
 \label{eq:ADMq2SubleadingCoefficient}
\end{equation}
Consequently,
\[
 \begin{aligned}
 f^+(\zeta)
 &=1+\frac{2q_2^+(\zeta)}{q_1^+(\zeta)-q_2^+(\zeta)}\\
 &=1+\frac{2\mathcal G_N^+}{H_N^+}\frac1\zeta
   +O(\zeta^{-2}),
 \end{aligned}
\]
which proves \eqref{eq:CandidateADMDeterminantFormula} by comparison with
\eqref{eq:CandidateADMDefinition}.  Before the equilibrium conditions are
imposed, this coefficient need not be real; it is therefore understood
here in the same algebraic candidate sense as the other quantities in
Section~6.

\begin{conjecture}[Candidate mass decomposition]
For determinant data satisfying the equilibrium conditions, the candidate
ADM mass is the sum of the candidate Komar masses:
\begin{equation}
 M_{\mathrm{ADM}}
 =\sum_{j=1}^{N}\mathfrak M_j.
 \label{eq:CandidateMassDecomposition}
\end{equation}
\end{conjecture}

For \(N=1\), identity~\eqref{eq:CandidateMassDecomposition} follows
directly from the formulas in Section~7.1.  For \(N=2\), it is the mass
identity obtained by Hennig and Neugebauer from the asymptotic expansion of
the equilibrium double-Kerr-NUT family
\cite{NeugebauerHennig2009}.  We have also verified the identity for a
finite collection of \(N=3\) parameter sets satisfying the five equilibrium
equations of Section~7.3.  No proof for general \(N\) is presently known.

The second conjecture concerns zeros of the determinant denominator in the
closed Weyl half-plane.

\begin{conjecture}[Zeros of the determinant denominator]
Let \(N\geq2\), and suppose that the determinant data are non-degenerate and
satisfy the equilibrium conditions.  Then \(\Delta_N\) has at least one
zero in the closed Weyl half-plane
\begin{equation}
 \{(\rho,\zeta)\in\mathbb R^2:\rho\geq0\}.
 \label{eq:ClosedWeylHalfPlane}
\end{equation}
\end{conjecture}

For \(N=2\), this assertion was not proved by Hennig and Neugebauer.  They
examined numerically a large number of parameter values satisfying the
special equilibrium relations and found zeros of the determinant
denominator in every case considered.  They further checked that the
numerator does not vanish at the same points, so that these zeros produce
singular rings outside the candidate horizons.  Their numerical
investigations led them to conjecture that every equilibrium
double-Kerr-NUT candidate suffers from such singular rings; see
Section~5.2 of \cite{NeugebauerHennig2009}.  The conjecture above extends
this \(N=2\) conjecture to arbitrary \(N\geq2\).

For \(N=3\), we have tested finitely many parameter sets satisfying the
equilibrium equations.  In every case examined, \(\Delta_3\) has a zero in
the closed Weyl half-plane.  These calculations provide evidence, but not a
proof.

The restriction \(N\geq2\) is essential.  Indeed, the \(N=1\) calculation
of Section~7.1 shows that the positive-mass Kerr family is characterized in
the present parametrization by
\begin{equation}
 \operatorname{Re}\alpha=\cos\theta>0.
 \label{eq:ConjecturesPositiveKerrRange}
\end{equation}
For precisely this parameter range, \(\Delta_1\) is nowhere zero in the
closed Weyl half-plane.

The third conjecture is the general counterpart of the obstruction used in
the two-horizon non-existence proof.

\begin{conjecture}[Failure of simultaneous sub-extremality]
Let \(N\geq2\), and suppose that the determinant data satisfy the
equilibrium conditions and that all candidate ratios
\begin{equation}
 p_j=\frac{8\pi J_j}{A_j},
 \qquad j=1,\ldots,N,
 \label{eq:ConjecturesCandidateRatios}
\end{equation}
are defined.  Then the inequalities
\begin{equation}
 |p_j|<1,
 \qquad j=1,\ldots,N,
 \label{eq:SimultaneousSubextremality}
\end{equation}
cannot all hold simultaneously.
\end{conjecture}

For \(N=2\), this follows from the explicit equilibrium parametrization and
the formulas reviewed in Section~7.2.  For \(N=3\), finite exact and
numerical tests have so far produced no parameter set for which all three
inequalities \(|p_1|<1\), \(|p_2|<1\), and \(|p_3|<1\) hold.  Again, this is
evidence rather than a proof.

\begin{remark}[The supplementary branches for \(N=2\)]
Recall the supplementary discussion at the end of Section~4.  For the
two-horizon family satisfying the equilibrium conditions, we obtained the
factorizations
\[
 H_2^{\mathcal H^1}
 =(w'+\mathit\Phi)Q_\varepsilon^{(1)}
\]
and
\[
 H_2^{\mathcal H^2}
 =(1-w\mathit\Phi)Q_\varepsilon^{(2)};
\]
see \eqref{eq:H1balancedfactorization} and
\eqref{eq:H2balancedfactorization}.  The vanishing of
\(H_2^{\mathcal H^j}\) therefore has more than one possible algebraic
source.  We supplemented the equilibrium parametrization by imposing the
following additional requirements: whenever
\[
 H_2^{\mathcal H^1}=0,
\]
require
\[
 w'+\mathit\Phi=0;
\]
and whenever
\[
 H_2^{\mathcal H^2}=0,
\]
require
\[
 1-w\mathit\Phi=0.
\]
These are sufficient supplementary requirements selected in order to
obtain the constant boundary values
\[
 a|_{\mathcal H^1}=0,
 \qquad
 a|_{\mathcal H^2}=0,
\]
respectively; see \eqref{eq:H1supplementarycondition} and
\eqref{eq:H2zeroa}.  They are not consequences of
\(H_2^{\mathcal H^j}=0\), and their necessity is not claimed.

Away from the two branches selected by these supplementary requirements,
both expressions in \eqref{eq:NTwoPValues} are directly defined.  The
argument of Hennig and Neugebauer then shows that \(|p_1|<1\) and
\(|p_2|<1\) cannot hold simultaneously.

It remains to consider the selected supplementary branches.  Their defining
factors cannot vanish simultaneously.  Indeed,
\(w'+\mathit\Phi=0\) gives \(\mathit\Phi=-w'<0\), whereas
\(1-w\mathit\Phi=0\) gives \(\mathit\Phi=1/w>0\).  Equivalently,
simultaneous vanishing would require \(ww'=-1\), which is impossible for
\[
 w>1,
 \qquad
 0<w'<1.
\]

If \(w'+\mathit\Phi=0\), the standard expression for \(p_1\) is not
directly defined, but the expression for \(p_2\) remains regular and gives
\[
 p_2
 =\varepsilon\frac{w(w+w')}{1+ww'}.
\]
Since
\[
 w(w+w')-(1+ww')=w^2-1>0,
\]
we obtain \(|p_2|>1\).  Similarly, if \(1-w\mathit\Phi=0\), the standard
expression for \(p_2\) is not directly defined, whereas
\[
 p_1
 =\varepsilon\frac{w+w'}{w'(1+ww')}.
\]
Now
\[
 (w+w')-w'(1+ww')=w(1-{w'}^2)>0,
\]
and hence \(|p_1|>1\).

Consequently, on either branch selected by the supplementary requirements,
the non-exceptional component already violates the strict sub-extremal
inequality.  The \(N=2\) obstruction therefore does not require the
candidate quantity associated with the exceptional component to be defined
by a limiting procedure.  This extends the argument of Hennig and
Neugebauer to the two cases in which one of the standard algebraic
expressions for \(p_j\) is not directly defined.
\end{remark}

The same \(N=1\) Kerr family in
\eqref{eq:ConjecturesPositiveKerrRange} explains why this conjecture must
also begin with \(N=2\).  Section~7.1 gives
\begin{equation}
 p_1
 =\frac{\sin\theta}{1+\cos\theta}
 =\tan\frac{\theta}{2},
 \label{eq:ConjecturesNOneP}
\end{equation}
and proves
\begin{equation}
 \mathfrak M>0
 \quad\Longleftrightarrow\quad
 |p_1|<1
 \quad\Longleftrightarrow\quad
 \operatorname{Re}\alpha=\cos\theta>0.
 \label{eq:ConjecturesNOneEquivalence}
\end{equation}
Thus one and the same positive-mass Kerr family simultaneously satisfies
\begin{equation}
 \Delta_1\neq0
 \quad\hbox{throughout the closed Weyl half-plane},
 \qquad
 |p_1|<1.
 \label{eq:ConjecturesNOneCounterexamples}
\end{equation}

We finally indicate a possible conditional relation between the last two
conjectures.  This is the natural \(N\)-component extension of the mechanism
used for the two-Kerr problem in \cite{ChruscielEtAlTwoKerr}.  Suppose that
the candidate horizons are regular Killing horizons and that the associated
initial hypersurface satisfies the global regularity assumptions required
there.  The area--angular-momentum inequality (see, for example,
\cite{DainReiris2011}) then gives
\begin{equation}
 A_j\geq 8\pi|J_j|
 \qquad (j=1,\ldots,N).
 \label{eq:DainReirisConditionalInequality}
\end{equation}
After exclusion of the extreme equality case, one obtains
\begin{equation}
 |p_j|<1
 \qquad (j=1,\ldots,N)
 \label{eq:DainReirisConditionalStrictInequality}
\end{equation}
simultaneously, contrary to Conjecture~8.3.  Thus, if Conjecture~8.3 holds,
some global regularity assumption must fail.  The zeros of \(\Delta_N\)
predicted by Conjecture~8.2 provide a natural algebraic candidate for this
failure, although this last identification remains unproved.

\section*{Acknowledgements}

The author would like to express his sincere gratitude to Hans-Joachim Hein
and Gustav Holzegel for many helpful and stimulating discussions, as well as
for their continued encouragement and support. He would also like to thank
Raimar Wulkenhaar and Gilles Jordan Pedieu for valuable discussions and helpful
comments related to this work.

\medskip
\begingroup
\small
\noindent Mathematisches Institut, Universit\"at M\"unster, Einsteinstr.~62,
48149 M\"unster, Germany.\\
E-mail: \href{mailto:dr.daheng.min@uni-muenster.de}{dr.daheng.min@uni-muenster.de}
\endgroup

\end{document}